\documentclass[acmsmall,screen]{acmart-primary/acmart}

\usepackage{mydef}

\usepackage{xcolor}
\usepackage{subcaption}
\usepackage{hyperref}
\usepackage{bookmark}
\usepackage[capitalise]{cleveref}
\AtBeginDocument{
\hypersetup{
  colorlinks  = true,
  linkcolor   = blue,        
  citecolor   = teal,        
  urlcolor    = cyan,        
  filecolor   = magenta,
}
}

\usepackage{algorithm}
\usepackage{algpseudocode}

\usepackage{listings}
\usepackage{tcolorbox}
\tcbuselibrary{skins,breakable}
\usepackage{tikz}

\definecolor{mGreen}{rgb}{0,0.6,0}
\definecolor{mGray}{rgb}{0.5,0.5,0.5}
\definecolor{mPurple}{rgb}{0.58,0,0.82}
\definecolor{backgroundColour}{rgb}{0.95,0.95,0.95}

\definecolor{codegreen}{rgb}{0,0.6,0}
\definecolor{codegray}{rgb}{0.5,0.5,0.5}
\definecolor{codepurple}{rgb}{0.58,0,0.82}
\definecolor{backcolour}{rgb}{0.95,0.95,0.95}

\lstdefinestyle{mystyle}{
    commentstyle=\color{codegreen},
    keywordstyle=\color{magenta},
    numberstyle=\tiny\color{codegray},
    stringstyle=\color{codepurple},
    basicstyle=\ttfamily\footnotesize,
    breakatwhitespace=false,         
    breaklines=true,                 
    captionpos=b,                    
    keepspaces=true,                 
    numbers=left,
    firstnumber=1,
    numbersep=8pt,
    xleftmargin=20pt,
    showspaces=false,                
    showstringspaces=false,
    showtabs=false,                  
    tabsize=2,
    language=C++
}

\tcolorboxenvironment{lstlisting}{
    enhanced,
    colback=backcolour,
    colframe=backcolour,
    boxrule=0pt,
    arc=3pt,
    left=2pt, right=2pt, top=2pt, bottom=2pt,
    breakable,
}

\AtBeginDocument{%
  }

\setcopyright{none}
\renewcommand\footnotetextcopyrightpermission[1]{}

\acmConference[]{}{}{}
\acmPrice{}
\acmYear{}
\acmDOI{}
\acmISBN{}
\copyrightyear{}
\acmVolume{}
\acmNumber{}
\acmArticle{}

\makeatletter
\def\@acmArticle{}
\def\ACM@history{}
\def\@received{}
\makeatother
\AtBeginDocument{%
  \fancyfoot{}%
  \fancyfoot[C]{\thepage}%
  \fancypagestyle{firstpagestyle}{%
    \fancyhf{}%
    \fancyfoot[C]{\thepage}%
  }%
}

\begin{document}

\title{MultiWave: A computational laboratory for adaptive numerical methods approximating hyperbolic balance laws}


\author{Adrian Kolb}
\email{kolb@eddy.rwth-aachen.de}
\affiliation{%
	\institution{Institut für Geometrie und Praktische Mathematik, RWTH Aachen University}
	\streetaddress{Im S\"usterfeld 2}
	\city{Aachen}
	\country{Germany}
	\postcode{52076}
}


\author{Aleksey Sikstel}
\email{sikstel@acom.rwth-aachen.de}
\affiliation{%
	\institution{Institut für Geometrie und Praktische Mathematik, RWTH Aachen University}
	\streetaddress{Kreuzherrenstrasse 2}
	\city{Aachen}
	\country{Germany}
	\postcode{52062}
}

\renewcommand{\shortauthors}{Kolb et al.}

\begin{abstract}
	The \textsl{MultiWave} C++-framework for adaptive numerical methods approximating hyperbolic balance laws~\revASSS{is presented}. \textsl{MultiWave} has been designed as a computational laboratory where new mathematical concepts can be quickly implemented and tested.
	\revAK{Starting from the mathematical background and proceeding to the low-level implementation details,  the realisation of a discontinuous Galerkin method with multiresolution-based grid adaptation \revASSS{is demonstrated}.}
	The design choices, in particular regarding the modularity that allows one to extend the code reusing existing infrastructure, \revASSS{is discussed}.
	\revAK{Scaling studies on a distributed-memory machine show that the framework retains its efficiency up to large rank counts, and that multiresolution-based adaptivity reduces the cost per time step by orders of magnitude compared to the uniform-grid computation at \revakk{the order of accuracy} of the \revASSS{full} reference discretisation.}
\end{abstract}




\maketitle

\section{Introduction}

Systems of hyperbolic balance laws constitute an important class of first-order partial differential equations (PDEs). Their solutions model diverse phenomena in fields such as continuum physics, biomedical sciences, manufacturing, logistics, and traffic flow~\cite{bessonov2016methods,  liu2024distributed, othman2022pde,  colombo2011modelling, garavello2006traffic, colombo2003hyperbolic, dafermosHyperbolicConservationLaws2016}.
\begin{figure*}[htbp]
	\centering
	\includegraphics{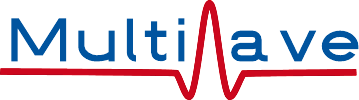}
\end{figure*}

Given the space-time domain \(\Omega_T := [0,T] \times \Omega\) where \(\Omega\subset \mathbb{R}^d\) is a \(d\)-dimensional spatial domain and \(T>0\) the final time, we consider  systems of PDEs in the form
\revASS{\begin{equation}
		\label{eq:general-pde}
		\frac{\partial\vec{u}}{\partial t} + \nabla\cdot \vec{F}(\vec{u})  + \vec{S}(\vec{x}, \vec{u}) = \vec{0}.
	\end{equation}
	The vector of unknown functions \revASSS{depending on} time and space reads \(\vec{u} \,\colon\, \Omega_T \to D\subset\R^m\), the vector of flux functions is \(\vec{F}(\vec{u})=\left(\vec{f}_1(\vec{u}), \ldots, \vec{f}_d(\vec{u})  \right)^T\) with  \(\vec{f}_i \,\colon\, D \to \R^{m}\) and \(\vec{S}\, \colon\, \Omega\times D \to \mathbb{R}^m\) denotes the source term where $D\subset \mathbb{R}^m$ is the set of admissible values. The system is endowed with initial and boundary conditions
	\begin{align*}
		\begin{split}
			 & \vec{u}(0, \vec{x}) = \vec{u}_0(\vec{x}) \text{ for } \vec{x} \in \Omega,                                                  \\
			 & \vec{u}(t, \vec{x}) = \vec{b}(t, \vec{x}, \vec{u}^-(t, \vec{x})) \text{ for } t>0, \vec{x} \in \partial\Omega_{\text{bc}},
		\end{split}
	\end{align*}
	where $\vec{b}\,\colon\, [0,T]\times \partial\Omega_{\text{bc}}\times D \to D$ and $\partial\Omega_{\text{bc}} \subset \partial\Omega$ is the inflow boundary.
	The system of balance laws~\eqref{eq:general-pde} is called hyperbolic if the matrix 
\begin{align}
\label{eq:flux-jacobian}
  \mathbf{A}(\vec{u}, \vec{n})  := \sum_{j=1}^d  \vec{n}_j\cdot  \vec{f}_j' (\vec{u})
\end{align}
has real eigenvalues and is diagonalisable in any  direction $\vec{n}\in \mathbb{R}^d \setminus \{\vec{0}\}$ and $\vec{u} \in D$.
}

\revASSS{Numerical simulations} of these problems require high-precision approximations whose accuracy is, ideally, rigorously controlled. \revASSS{Typically, their} solutions exhibit a heterogeneous structure and might consist of smooth areas with steep and flat gradients, vortices and discontinuities such as shock or contact waves~\cite{liu2021shock, bressanHyperbolicSystemsConservation2000}. Consequently, the development of accurate adaptive numerical methods based on a sound theoretical foundation is vital for applications. The \revASS{study} of such methods strongly benefits from numerical experiments, in particular from empirical error analysis and convergence studies.

\revASS{Implementation of numerical methods for hyperbolic \revASSS{PDEs} has a long history. Among the earliest widely adopted frameworks are finite volume codes build around structured or block-structured meshes. Codes like \textsl{Athena++} \cite{Stone2020} and \textsl{OpenFOAM} \cite{greenshields2025} address specific application domains: magnetohydrodynamics and general computational fluid dynamics, respectively. They offer robust, production-grade solvers for single-physics problems. Their design, however, reflects the assumptions of their target applications, i.e. fixed spatial discretisation order, domain-specific flux treatments, and limited support for  heterogeneous coupling that arises in multi-material or multi-phase systems.}

\revASS{The advent of high-order discontinuous Galerkin (DG) methods spurred a second wave of frameworks targeting compressible flow at high polynomial degree. \textsl{FLEXI}~\cite{flexi} and \textsl{HORSES3D}~\cite{ferrer2023horses} are \revASSS{representatives} of this class. They deliver high-order accuracy on hexahedral meshes with efficient tensor-product implementations  \revASSS{taylored} to scale on modern HPC architectures. \textsl{Nektar++}~\cite{cantwell2015nektar} broadens the scope to mixed element types and a wider range of PDE classes within a spectral/$hp$-element framework.}

\revASS{More recent frameworks have moved toward greater algorithmic flexibility and a more general treatment of hyperbolic problems. \textsl{ExaHyPE} \cite{reinarz2020exahype} provides an engine for systems of hyperbolic PDEs discretised by ADER-DG schemes with cell-local limiting and dynamic AMR, targeting exascale platforms via a task-based parallelism model. \textsl{Trixi.jl} \cite{ranocha2022adaptive}, written in Julia, has attracted considerable attention for its entropy-stable DG discretisations, composable solver architecture, and support for a broad range of equation systems. 
}

\revASS{The \textsl{Clawpack} software family \cite{clawpack} occupies a unique position in this landscape. Build around LeVeque's wave-propagation algorithms, it provides a unifying framework in which the Riemann solver is the central computational primitive, directly reflecting the mathematical structure of the equations. Extensions including \textsl{AMRClaw} and \textsl{GeoClaw} add adaptive mesh refinement and geophysical applications, while \textsl{PyClaw} provides a high-level Python interface. \textsl{Clawpack}'s philosophy of centring the discretisation on Riemann solutions at cell interfaces makes it conceptually closest to the approach taken in the present work. Nevertheless, it is based on finite volume methods of at most second order, provides no support for $hp$-adaptivity, and was not designed for the systematic treatment of non-conservative products or the coupling of distinct hyperbolic systems across material interfaces.}

\revASS{The present work introduces \textsl{MultiWave}\footnote{\url{https://www.igpm.rwth-aachen.de/forschung/multiwave}}, a framework designed  around a fundamentally different philosophy: rather than a production solver optimised for a fixed discretisation strategy, \textsl{MultiWave} is  a laboratory for numerical methods, in the spirit of \textsl{Clawpack} but operating at a deeper level of algorithmic abstraction. To this end, the main design decision was to expose the full computational hierarchy as independently exchangeable components: from the PDE system and Riemann solver down through the DG discretisation and basis representation, the refinement strategy and error estimator, the inter-domain coupling mechanism, and ultimately to the underlying grid data structures and elementary cell operations. Nothing in this hierarchy is hardwired. This is realised through a systematic use of C++ template metaprogramming, where each layer of the hierarchy is encoded as a compile-time policy type rather than a runtime abstraction. The consequence is that replacing e.g.~a Roe solver with a relaxation solver, switching from a modal to a nodal DG basis or coupling two physically distinct subsystems across an interface requires only a change at the template parameter level, with no modification to the surrounding framework and no runtime overhead. Beyond enabling systematic comparative studies, this design substantially lowers the cost of experimentation: new PDE systems, solvers, and refinement strategies can be prototyped and integrated into existing  pipelines with minimal friction, reducing the gap between a theoretical concept and a working numerical implementation drastically.}

\revASS{Remarkably, the flexibility does not come at the cost of performance. Because the entire hierarchy is resolved at compile time, the compiler has full visibility into the instantiated configuration, enabling aggressive inlining and optimisation that would be prohibited by virtual dispatch or runtime polymorphism.  \textsl{MultiWave} is built on a distributed-memory MPI architecture with dynamic load balancing across the adaptive mesh. This enables the framework to solve multi-dimensional problems and coupled multi-physics configurations, and to conduct convergence studies at high polynomial degree. In addition, parallel scalability is a prerequisite for the kinds of large-scale computations the framework targets.}

\revASS{The adaptive DG method serves as the primary discretisation workhorse of the framework. However, the template architecture is equally capable of implementing other schemes, such as CWENO or generalised Riemann problem methods, reusing the  parallelised adaptive infrastructure. Refinement decisions are driven by multiresolution analysis (MRA) or rigorous a posteriori error estimators. Path-conservative Riemann solvers are a native concept in the framework, and the template architecture ensures that solver design, discretisation, and adaptivity interact as first-class citizens at every level. This makes \textsl{MultiWave} a unique environment for systematically studying adaptivity and solver design across a broad range of hyperbolic PDE problems. The result is a codebase that functions simultaneously as a high-performance research instrument and a foundation for application-specific extensions.}

\revAK{
	The contributions of this work are the following:
	\begin{itemize}
		\item The entire adaptive DG method, from the PDE system and Riemann solver down to the grid data structures, is implemented as compile-time template policies. Components are exchanged at the template parameter level, without virtual dispatch (\cref{sec:multiwave}).
		\item For demonstration purposes, three method variants are realised within the framework: a second-order spatial derivative term, a path-conservative solver for non-conservative products, and a well-balanced shallow water scheme. Each is obtained by supplying new policy types, leaving the existing components unchanged (\cref{sec:extens-new-feat}).
		\item The distributed-memory design uses a Morton space-filling curve with subtree-atomic partitioning and dynamic load rebalancing. It reaches $90.6\,\%$ strong-scaling efficiency on 12\,288 MPI ranks (\cref{sec:performance}).
		\item Grid adaptivity is driven by multiresolution analysis with multiwavelets and a level-dependent threshold hierarchy. On adaptive grids, the framework retains a weak-scaling efficiency above $55\,\%$ from 384 to 12\,288 ranks, with the compute-bound time stepping at $91\,\%$ (\cref{sec:performance}).
	\end{itemize}
}



\revAK{Modularity is the connecting theme of this work, as the mathematical hierarchy of the method maps directly to the code, enabling each component to be studied and extended in isolation.} This work is organised as follows: we briefly introduce the mathematical background in \cref{sec:adapt-schemes}, \revAK{describe how the mathematical hierarchy maps to independently exchangeable code components} in \cref{sec:multiwave}, examine performance on distributed-memory machines in \cref{sec:performance}, and present numerical results in \cref{sec:appl-fram}.

\section{Adaptive numerical schemes for hyperbolic balance laws}
\label{sec:adapt-schemes}

In this section, we first introduce  the discretisation of systems of hyperbolic balance laws by means of  discontinuous Galerkin (DG) methods~\cite{cockburnRungeKuttaDiscontinuous1998}. This approach follows closely the mathematical analysis framework and, due to its locality, is well-suited for  dynamic local grids ($h$-adaptivity) and discretisation order variation ($p$-adaptivity). Secondly, we describe the \revASSS{multiresolution-based} adaptivity approach that is at the core of \textsl{MultiWave} and also serves as a foundation for further adaptivity methods.

\subsection{Discontinuous Galerkin methods}
\label{sec:disc-galerk-meth}


\revASSS{Since hyperbolic balance laws may develop discontinuities in finite time}, the system \revASSS{may} not  possess classical solutions which makes weak formulations mandatory.
\revASSS{The weak formulation is the starting point for deriving the DG method introduced in~\cite{cockburnTVBRungeKuttaLocal1989}, in which the trial and test spaces are chosen to be the same finite-dimensional space of piecewise polynomials.  }
\revASSS{To be more specific}, the spatial domain $\Omega$ is partitioned by a finite number of disjoint open cells $V_{\lambda}$ that are numbered by an index set $\mathcal{I}_h$, i.e.
\begin{equation*}
	\overline{\Omega} = \bigcup_{\lambda\in \mathcal{I}_h} \overline{V_{\lambda}},\text{ where } V_\lambda \cap V_\mu = \emptyset  \text{ for } \lambda \neq \mu \in \mathcal{I}_h.
\end{equation*}
The set $\mathcal{G}_h : = \{V_\lambda\}_{\lambda \in \mathcal{I}_h}$ is called the computational grid. On this grid the DG-space $\mathcal{F}_{h,p} $ of order $p$ is defined as
\begin{align*}
	\mathcal{F}_{h,p} := \left\{ f \in L^2(\Omega) \,\colon\, f|_{V_{\lambda}} \in \Pi _{p-1}(V_{\lambda}),\, \lambda \in \mathcal{I}_h \right\},
\end{align*}
where $\Pi_{p-1}(V)$ denotes the space of polynomials of \revAK{tensor} degree less than $p$. For the DG-space $\mathcal{F}_{h,p}$  we require a set of basis functions  $\Phi_h := \{\varphi_{\lambda, \vec{i}}\}_{\lambda \in \mathcal{I}_h, \, \vec{i} \in \mathcal P}$ with the index set \(\mathcal{P}:=\{ \vec{i}\in\N^d_0 \,\colon\, \|\vec{i}\|_{\infty} \leq p-1\}\). The elements of the basis are assumed to be orthogonal with respect to the standard $L^2(\Omega)$-inner product and compactly supported, i.e.
\begin{align*}
	\langle \varphi_{\lambda', \vec{i}'}, \varphi_{\lambda, \vec{i}} \rangle_{\Omega} = \delta_{\lambda',\lambda}\delta_{\vec{i}',\vec{i}}, \quad   \text{supp } \varphi_{\lambda,\vec{i}}  = V_{\lambda},
\end{align*}
for all $\lambda, \lambda' \in \mathcal{I}_h$ and $\vec{i}, \vec{i}'\in \mathcal{P}$.

The weak solution of the system of balance laws~(\ref{eq:general-pde}) is approximated in the DG-space by an expansion of the basis $\Phi_h$ as follows
\begin{equation}
	\label{eq:dg-space-approx}
	\vec{u}(t, \,\cdot\,) \approx \vec{u}_h(t, \,\cdot\,) := \sum_{\lambda\in \mathcal{I}_h} \sum_{\vec{i}\in \mathcal{P}} \vec{u}_{\lambda,\vec{i}}(t) \varphi_{\lambda,\vec{i}}(\,\cdot\,) \in (\mathcal{F}_{h,p})^m.
\end{equation}
Due to the orthogonality of the basis $\Phi_h$, the coefficients $\vec{u}_{\lambda,\vec{i}}(t)$ are recovered as
\begin{align*}
	\vec{u}_{\lambda,\vec{i}}(t) \coloneq \langle \vec{u}_h(t, \,\cdot\,), \varphi_{\lambda,\vec{i}} \rangle_{V_{\lambda}}.
\end{align*}

\revASSS{In order to handle fluxes at cell interfaces where the the numerical solution $\vec{u}_h$, defined in~\eqref{eq:dg-space-approx},  may be discontinuous,  the fluxes \revASS{$\vec{F}\cdot\vec{n}$} in the direction $\vec{n}$ are approximated by numerical fluxes $\widehat{\vec{f}}$. Hence, following~\cite{cockburnTVBRungeKuttaLocal1989}, one obtains the semi-discrete system of equations for every $\lambda \in \mathcal{I}_h$:
	\begin{equation}
		\label{eq:weak-approx-cl-semi-discr}
		\frac{d \vec{u}_{\lambda, \vec{i}}}{dt}
		= \int_{V_{\lambda}} \sum_{j=1}^d \vec{f}_j(\vec{u}_h) \cdot\partial_{\vec{x}_j} \varphi_{\lambda,\vec{i}} d\vec{x}
		-\! \int_{\partial V_{\lambda}} \widehat{\vec{f}} \left(\vec{u}_h^+, \vec{u}_h^-, \vec{n}_\lambda\right) \varphi_{\lambda,\vec{i}} dS
		-\! \int_{V_{\lambda}} \vec{S}(\vec{x}, \vec{u}_h) \varphi_{\lambda,\vec{i}} d\vec{x} \eqcolon \mathcal{R}_{DG}(\vec{u}_h),
	\end{equation}
}
\revASSS{Here,} the numerical flux $\widehat{\vec{f}}$ in direction $\vec{n}_\lambda$ is evaluated at the \revAK{outer} value $\vec{u}_h^+$ and the \revAK{inner} value $\vec{u}_h^-$ of $\vec{u}_h$ at the boundary of $V_\lambda$, i.e.~ $\vec{u}_h^{\pm}(\vec{x}) = \lim_{\varepsilon \to 0} \vec{u}_h(\vec{x} \pm \varepsilon \vec{n}_{\lambda})$ \revASSS{and} $\vec{n}_\lambda$ denotes the outward unit normal vector to the cell boundary $\partial V_\lambda$.

The semi-discrete system~(\ref{eq:weak-approx-cl-semi-discr}) is a system of ordinary differential equations (ODEs) \revASSS{$\vec{u}_{\lambda, \vec{i}}'(t) = \mathcal{R}_{DG}(\vec{u}_{h}(t))$}.
This system is solved by strong stability preserving Runge--Kutta methods~\cite{gottlieb2001strong} in compliance with time-step bounds given by the Courant--Friedrichs--Lewy (CFL) condition~\cite{chalmersRobustCFLCondition2020}. In order to guarantee convergence in the mean of the DG-method the discretisation needs to be stabilised by a local projection limiter, see e.g.~\cite{cockburnTVBRungeKuttaLocal1989}. \revAK{The limiter eliminates spurious oscillations caused by the Gibbs phenomenon by locally reducing the discretisation order to linear polynomials, thereby restoring monotonicity.}

\begin{rem}
	The DG method can be adapted to mixed systems, i.e.~balance laws~\eqref{eq:general-pde} with additional differential operators like viscosity or non-conservative products, see~\cite{bassi1997high} and~\cite{dal1995definition}, respectively. For the sake of brevity, we restrict ourselves here to hyperbolic balance laws. The following section on discretisation explains how \textsl{MultiWave} is extended to incorporate such operators.
\end{rem}

\subsection{Multiresolution analysis and \revakk{grid adaptation}}
\label{sec:grid-adaptivity}

The idea of \revAK{multiresolution-based} grid adaptation methods is to consider the discretisation on an adaptive grid as a perturbation of the discretisation on the full reference grid. This approach allows \revAK{us} to control the induced perturbation error by keeping significant contributions only.  Thereby, the accuracy of the reference solution is asymptotically maintained at a fraction of the cost that would be necessary for a computation on the full reference grid.
This approach has been successfully applied in the context of hyperbolic balance laws \cite{Gerhard2017,Gerhard2015,bramkamp2004adaptive,dahmenMultiresolutionSchemesConservation2001, hovhannisyanAdaptiveMultiresolutionDiscontinuous2014, Caviedes-Voullieme2020, gerhardHighOrderDiscontinuousGalerkin2015, gerhard2016adaptive, gerhardWaveletFreeApproachMultiresolutionBased2021, harten1995multiresolution} using the concept of MRA  introduced in~\cite{mallat1989multiresolution}.

\revAK{The} MRA \revakk{is determined by a sequence $\mathcal{S}\coloneq \left\{ \mathcal{S}_{\ell} \right\}_{\ell\in \mathbb{N}_0}$ of nested and closed subspaces $\mathcal{S}_{\ell}$ of $L^2(\Omega)$} \revAK{on a hierarchy} of nested grids $\mathcal{G}_{\ell} \coloneq \left\{ V_{\lambda}  \right\} _{\lambda \in \mathcal{I}_{\ell}}$ where $\ell\in\mathbb{N}_0$ denotes the refinement level. \revakk{This MRA is assumed to be dense in $L^2(\Omega)$.} \revAK{Two grids at successive refinement levels $\ell$ and $\ell+1$ are called nested if each cell $V_\lambda, \lambda\in\mathcal I_\ell$, on level $\ell$ is composed of cells at level $\ell +1$, i.e., $V_\lambda = \bigcup_{\mu\in\mathcal M_\lambda}V_\mu$ for the refinement set $\mathcal M_\lambda \subset \mathcal I_{\ell+1}$ of cell $V_\lambda$}. To simplify the notation, we restrict ourselves to Cartesian dyadic grids, although \revAK{the} MRA can be applied, for instance, on triangulations as well~\cite{Yu1999}.
For a precise formulation of the details in the MRA we refer to the  literature  cited above  and the references therein.

\subsubsection*{\textbf{Multiscale decomposition}:}

\begin{figure}[!htb]
	\centering
	\begin{tikzpicture}
		\node[anchor=south west, inner sep=0] (img) at (0,0)
		{\includegraphics[width=0.4\textwidth]{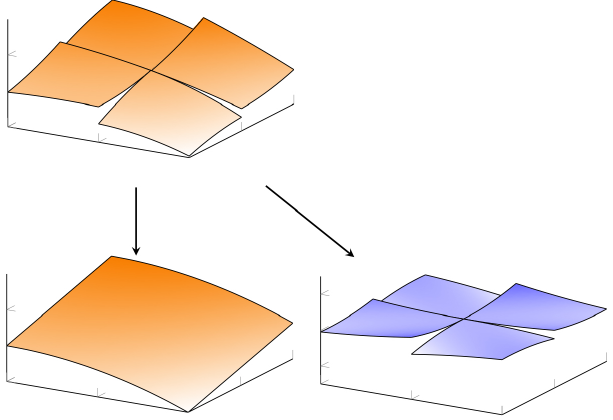}};
		\begin{scope}[x={(img.south east)},y={(img.north west)}]
			\node[anchor=east] at (0.02, 0.83) {$\vec{u}^{\ell+1}_\lambda$};
			\node[anchor=east] at (0.02, 0.24) {$\vec{u}^\ell_\lambda$};
			\node[anchor=west] at (0.98, 0.24) {$d^\ell_\lambda$};
			\node[anchor=east] at (0.20, 0.50) {$P_{\mathcal{S}_\ell}$};
			\node[anchor=west] at (0.47, 0.52) {$P_{\mathcal{W}_\ell}$};
		\end{scope}
	\end{tikzpicture}
	\caption{Illustration of local two-scale decomposition $\vec u^{\ell+1}_\lambda = \vec u^\ell_\lambda + \vec d^\ell_\lambda$, cf. \cite{Gerhard2017}.}
	\label{fig:mra-decomposition}
\end{figure}

The MRA is employed to assess the difference of the solution between two refinement levels by introducing the orthogonal complement space
\begin{align*}
	\revAS{	\mathcal{W}_{\ell} \coloneq \left\{ d\in \mathcal{S}_{\ell+1} \,\colon\, \langle d, u \rangle_{\Omega}  = 0, u\in \mathcal{S}_{\ell} \right\}},
\end{align*}
such that $\mathcal{S}_{\ell+1} = \mathcal{S}_{\ell} \oplus \mathcal{W}_{\ell}$ for any $\ell\in \mathbb{N}_0$. In other words, the space $\mathcal{W}_{\ell}$ contains the information that is lost when projecting a solution $u^{\ell+1}\in S_{\ell+1}$ from level $\ell+1$ onto level $\ell$. By recursively applying the two-scale decomposition, one obtains the \revAK{multiscale} representation $\mathcal{S}_{\ell}= \mathcal{S}_0 \oplus \mathcal{W}_0 \oplus \mathcal{W}_1 \oplus \cdots \oplus \mathcal{W}_{\ell-1}$, i.e. any function $u \in L^2(\Omega)$ can be decomposed into its different scales $u = u^0 + \sum_{\ell\in\mathbb{N}_0} d^{\ell}$ with

\begin{equation}
	\label{eq:orth-projections}
	u^\ell \coloneq P_{S_{\ell}}(u) =  P_{S_{\ell}}(u^{\ell+1}), \quad d^{\ell} \coloneq P_{\mathcal{W}_{\ell}}(u) = P_{\mathcal{W}_{\ell}}(u^{\ell+1}), \quad \ell\in \mathbb{N}_0,
\end{equation}
where $P_{\mathcal A}\,\colon\, L^2(\Omega) \to \mathcal{A}$ denotes the $L^2$-projection to a closed subspace $\mathcal{A}\subset L^2(\Omega)$.
The two-scale decomposition, $u^{\ell+1} = u^{\ell} + d^{\ell} $ is illustrated in \cref{fig:mra-decomposition}.

Since the DG-solution is in the space of piece-wise polynomials, the projections $P_{\mathcal{S}_{\ell}}$ and $P_{\mathcal{W}_{\ell}}$ can be localised on each element of the grid by setting $u^{\ell}_{\lambda} \coloneq u^{\ell} \mychi{V_{\lambda}}$ and $d^{\ell}_{\lambda} \coloneq d^{\ell} \mychi{V_{\lambda}}$ where $\ell\in \mathbb{N}_0$, $V_\lambda \in \mathcal{G}_{\ell}$ is a grid cell and $\mychi{A}$ is the characteristic function of a set $A$. With this notation the localised multi-scale decomposition reads
\begin{equation}\label{eq:loc-ms-decomposition}
	u = \sum_{\lambda\in \mathcal I_0} u^0_{\lambda} + \sum_{\ell\in \mathbb{N}_0} \sum_{\lambda \in \mathcal I_{\ell}} d^{\ell}_\lambda,
\end{equation}
where $d^{\ell}_\lambda\in \mathcal W_\ell$ are the local contributions of the details on each level.

In order to interpret these contributions, i.e. the norms of the details, in relation to the smoothness of functions $u\in H^p(V_{\lambda})$, we follow~\cite{Gerhard2017} and apply Whitney's theorem together with the Cauchy--Schwarz inequality to obtain
\begin{equation}
	\label{eq:cancellation-prop}
	\| d_{\lambda}^{\ell}\|_{L^2(V_\lambda)} \leq \operatorname{diam}(V_{\lambda})^p \sum_{\| \boldsymbol{\alpha} \|_1 = p } \frac{1}{\boldsymbol{\alpha}!} \|D^{\boldsymbol\alpha}u\|_{L^2(V_\lambda)}.
\end{equation}
This inequality is \revAK{essential} as it implies that the norm of local contributions decays  with the grid size (i.e.~with increasing $\ell$) at order $p$ as long as all weak derivatives of $u$ up to $p$-th order are bounded. \revAK{Thus, the norm of the details} $\|d_{\lambda}^{\ell}\|_{L^2(V_\lambda)}$, is associated with the smoothness of $u$ and can be employed to drive \revakk{data compression} and consequently grid adaptivity.

\subsubsection*{\textbf{Thresholding}:}

\revakk{The inequality \eqref{eq:cancellation-prop} allows for} a sparse approximation $u^{L,\varepsilon}$ of the solution $u^L$ at a fixed maximal refinement level $L\in\mathbb N$. To this end, we employ the technique of hard thresholding, which consists in discarding non-significant local contributions according to a hierarchy of local threshold values $\varepsilon_{\lambda, L} = 2^{\ell-L}\varepsilon_L$ with a global threshold value $\varepsilon_L > 0$. A local contribution is called \emph{significant} if
\begin{equation}
	\label{eq:significant-detail}
	\| d_{\lambda}^{\ell}\|_{L^2(V_\lambda)}  > \varepsilon_{\lambda, L}.
\end{equation}
\revAK{A sparse approximation is then obtained by discarding non-significant details and projecting the local DG solution onto a coarser grid.}
\revAK{Then,} the sparse approximation $u^{L,\varepsilon} \in \mathcal{S}_L$ of $u^L$ is defined as
\begin{equation}
	\label{eq:sparse-approx}
	u^{L,\varepsilon} \coloneq \sum_{\lambda \in \mathcal I_0}u^0_{\lambda} + \sum_{\ell=0}^{L-1}\sum _{\lambda \in \mathcal{D}_{\varepsilon} \cap \mathcal{I}_{\ell}} d^{\ell}_{\lambda},
\end{equation}
where the index set $\mathcal{D}_{\varepsilon}$ is defined as the smallest set containing the indices of significant contributions
\begin{equation}
	\label{eq:significant-indices}
	\mathcal{D}_\varepsilon \supset \left\{ \lambda \in \bigcup_{\ell=0}^{L-1} \mathcal I_{\ell} \,\colon\, \|d^{\ell}_{\lambda} \|_{L^2(V_{\lambda})} > \varepsilon_{\lambda, L}  \right\}
\end{equation}
and being a tree, i.e., $\lambda\in\mathcal D_\varepsilon$ \revAK{implies for all} $\mu\in \mathcal D_\varepsilon$ with $V_\lambda \subset V_\mu$.
Given $\sum_{\ell=0}^{L-1} \max_{\lambda \in \mathcal I_{\ell}} \varepsilon_{\lambda, L} \leq \varepsilon_{\text{max}}$ the \revAK{approximation} error is bounded by
\begin{equation}
	\label{eq:error-sparse-approx}
	\| u^L - u^{L,\varepsilon}\|_{L^q(\Omega) } \leq |\Omega|^{1/q} C \varepsilon_{\text{max}}, \text{ for } q=1,2.
\end{equation}
Furthermore, the thresholding procedure is $L^2$-stable in the sense of $\|u^{L,\varepsilon}\|_{L^2(\Omega)} \leq \|u^{L}\|_{L^2(\Omega)}$ and conservative, i.e., $\int_{V_\lambda} u^L(\vec x)d\vec x = \int_{V_\lambda} u^{L,\varepsilon}(\vec x)d\vec x$. For more details we refer to \cite[Sec. 3]{gerhardWaveletFreeApproachMultiresolutionBased2021} and the references therein.

\begin{rem}
	\revAK{
		The thresholding procedure depends only on the solution data, not on the underlying PDE. It can therefore be used to compress data to the sparse grid within the approximation error bound~\eqref{eq:error-sparse-approx}. This has been applied, e.g., to the compression of climate reanalysis data, see~\cite{Khosrawi2026}.
	}
\end{rem}

\subsubsection*{\textbf{Prediction}:}
\revAK{To ensure that the solution of the DG scheme is appropriately resolved at both the old time step $t_n$ \emph{and} the new time step $t_{n+1}$, we need to perform grid refinement before updating the solution.}
Following Harten~\cite{harten1995multiresolution}, we employ the prediction strategy in a slightly modified form, see~\cite[Sec. 5]{gerhardWaveletFreeApproachMultiresolutionBased2021}:
\begin{enumerate}
	\item significant details remain significant: $\mathcal D^{n}_\varepsilon \subset \mathcal D^{n+1}$,
	\item neighbours of significant details become significant:
	      \begin{equation*}
		      \lambda \in \mathcal D^n_\varepsilon \Rightarrow
		      \{\tilde\lambda: V_{\tilde\lambda} \text{ is neighbour of } V_\lambda\}\subset \mathcal D^{n+1},
	      \end{equation*}
	\item significant details may become significant on a higher refinement level:
	      \begin{equation*}
		      \|d^\ell_\lambda\|_{L^2(V_\lambda)} > 2^{p+1} \varepsilon_{\lambda, L} \Rightarrow \mathcal M_\lambda \subset \mathcal D^{n+1},
	      \end{equation*}
	      with refinement set $\mathcal M_\lambda \subset \mathcal I_{\ell+1}$ of cell $V_\lambda$ and
	\item $\mathcal D^{n+1} $ is a tree.
\end{enumerate}
Property~(2) is justified by the CFL condition, which ensures that information in one cell propagates at most into its direct neighbouring cells during one time step. Property~(3) accounts for the development of sharp gradients within a cell. Although the reliability of the prediction, i.e. the property that no significant cells are discarded after performing a time step, has not been proven in general, the strategy has been successfully applied to a variety of problems, see~\cite{Gerhard2017} and references therein.

\subsubsection*{\revakk{\textbf{Grid adaptation}}:}
\revakk{To perform grid adaptation, the tree structure of the significant details~\eqref{eq:significant-indices} has to be translated onto an adaptive grid representing the sparse approximation~\eqref{eq:sparse-approx} in terms of single-scale coefficients. To this end, we define an adaptive grid $\mathcal G_{L,\varepsilon} := \cup_{\lambda \in \mathcal I_\varepsilon} V_\lambda$ with index set
	\begin{align*}
\mathcal I_{L,\varepsilon} := \left( \mathcal I_0 \setminus \mathcal D_\varepsilon \right) \cup \left\{
		\mu \in \bigcup_{\ell = 1}^L \mathcal I_\ell : \exists \lambda \in \mathcal D_\varepsilon : \mu \in \mathcal M_\lambda \text{ and } \mu \not\in \mathcal D_\varepsilon
		\right\},
		\intertext{such that}
		\overline\Omega = \overline{\bigcup_{\lambda \in \mathcal I_\varepsilon} V_\lambda} \quad\text{ and }\quad V_\lambda \cap V_\mu = \emptyset \quad\text{ for all }\quad \lambda \neq \mu \in \mathcal I_\varepsilon.
	\end{align*}
On the adaptive grid $\mathcal G_\varepsilon$ we denote a DG space $\mathcal S_{L,\varepsilon} \subset \mathcal S_L$. The sparse approximation~\eqref{eq:sparse-approx} is then represented in $\mathcal S_{L,\varepsilon}$ by the inverse multiscale transformation, yielding
	\begin{align}
		\label{eq:sparse-approx-dg}
		u^{L,\varepsilon} = \sum_{\lambda \in \mathcal I_0}u^0_{\lambda} + \sum_{\ell=0}^{L-1}\sum _{\lambda \in \mathcal{D}_{\varepsilon} \cap \mathcal{I}_{\ell}} d^{\ell}_{\lambda} =
\sum_{\ell=0}^L \sum_{\lambda \in \mathcal I_\ell \cap \mathcal I_{L,\varepsilon}} u^\ell_\lambda.
	\end{align}}

\subsubsection*{\textbf{Adaptive RKDG scheme}:}
The prediction step, the RKDG time step, and the \revAK{thresholding} step constitute the complete adaptive RKDG scheme, summarised in \cref{alg:general-adaptive-method}. For the initialisation step \textsc{adapt\_init}, a bottom-up strategy is employed, in which the initial data are projected level by level and cells are retained wherever the detail coefficients are significant, see~\cite{Gerhard2017}.

\begin{algorithm}
	\caption{Time marching of the \revAK{adaptivity method}}\label{alg:general-adaptive-method}
	\begin{algorithmic}
		\Require $\vec{u}_0\,\colon\, D \to \mathbb{R}^m$,  $L > 0$, $T > 0$
		\State $t \gets 0$
		\State $\vec{u}_h^{L,\varepsilon} \gets \textsc{adapt\_init}(\vec{u}_0)$
		\While{$t < T$}
		\State $\vec{u}_h^{L,\varepsilon} \gets$ \textsc{predict}($\vec{u}_h^{L,\varepsilon}$) \Comment{Refinement}
		\State $t \gets \textsc{step}(\vec{u}_h^{L,\varepsilon})$
		\State $\vec{u}_h^{L,\varepsilon} \gets$\textsc{coarsen}($\vec{u}_h^{L,\varepsilon}$) \Comment{Threshold}
		\State $\vec{u}_h^{L,\varepsilon} \gets$\textsc{limiter}($\vec{u}_h^{L,\varepsilon}$)
		\EndWhile
	\end{algorithmic}
\end{algorithm}

\begin{rem}
	\label{rem:threshold-choice}
	To preserve the accuracy of the reference scheme, i.e.,
	\begin{equation*}
		\|u - u^{L,\varepsilon}\| \approx \|u-u^L\| = \mathcal O(h_L^\beta)
	\end{equation*}
	for $L\to\infty$ and order $\beta$, the global threshold value $\varepsilon_L > 0$ must be chosen appropriately. In \cite{Gerhard2015}, a heuristic strategy using
	\begin{align}
		\label{eq:local-threshold-value}
		\varepsilon_L = C_\text{thr}h^\gamma_L
	\end{align}
	with some $C_\text{thr} > 0$ was proposed. For non-smooth solutions, setting $\gamma =1 $ and choosing $C_\text{thr}$ as the expected amplitude of the solution typically yields satisfactory numerical results, cf. \cite{Gerhard2017} and references therein.
\end{rem}

\begin{rem}
	Since limiting is required solely around discontinuities, a robust adaptivity method is important to maintain high discretisation order in smooth parts of the solution. This is achieved by applying the limiting procedure only in regions marked for highest grid resolution.
\end{rem}

\begin{rem}
	Besides providing an efficient and robust perturbation $h$-adaptivity method, MRA forms the basis for further analysis. For instance, MRA plays a key role in a posteriori error estimates derived in~\cite{giesselmann2025posteriori} where it is used to identify discontinuities.  Moreover, data structures of the MRA implementation are useful building blocks for design and validation of various adaptivity approaches, such as $hp$-adaptivity.
\end{rem}
In the following section on implementation, we will show how new features can be quickly incorporated into the \textsl{MultiWave}-ecosystem.

\section{Library design}
\label{sec:multiwave}

\begin{figure}[htbp]
	\centering
	\includegraphics[width=0.95\textwidth]{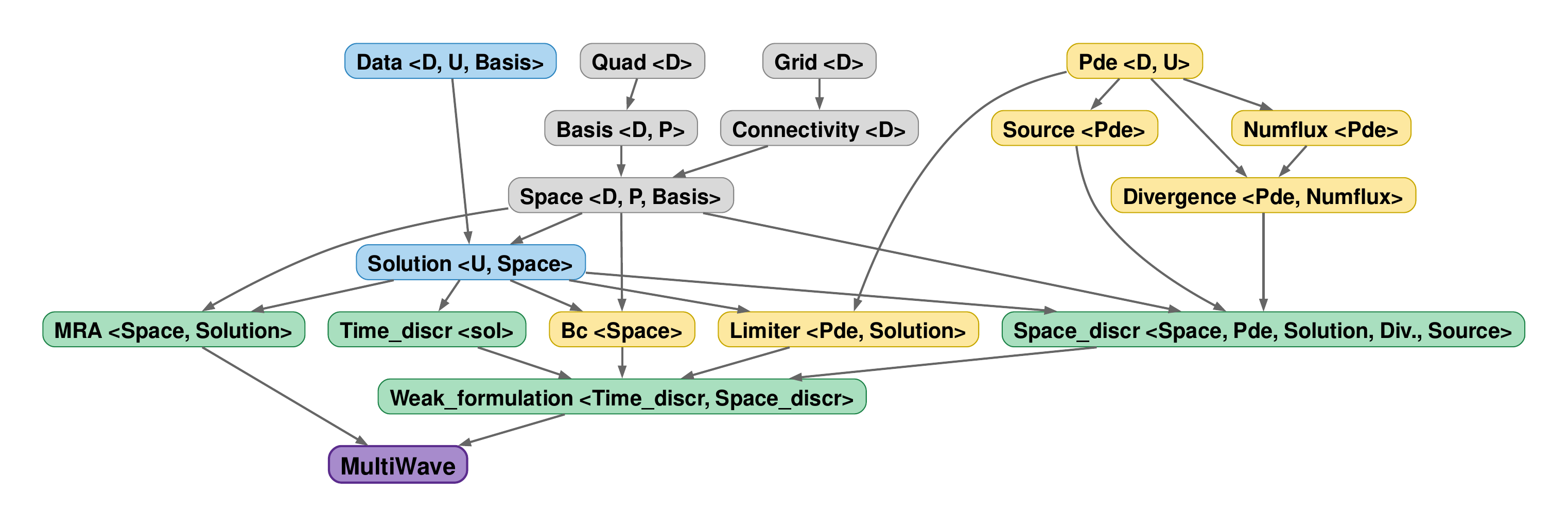}
	\caption{Dependency graph of the \textsl{MultiWave} library's template class hierarchy. Arrows indicate template parameter dependencies between components.}
	\label{fig:mw-structure}
\end{figure}

\textsl{MultiWave} is implemented in modern C++ and makes heavy use of object-oriented programing, template
programming and first-class functions. The code structure follows the mathematical
hierarchy of the numerical method, with each ingredient corresponding to one
independently exchangeable software component. \revAK{Two considerations motivate this choice.}
First, the library serves as a laboratory for mathematical research: it provides fast and
readable prototyping of theoretical ideas, with successful concepts subsequently admitted into
the code base.
\revAK{Every component of the hierarchy is a compile-time template parameter, so integrating
	a new PDE system, Riemann solver, or refinement strategy into an existing pipeline
	requires nothing more than a new class and a recompile.}
To this end, abstraction facilities of C++ are used as much as possible.
The overall structure of the library and the template dependencies between its components are
illustrated in \cref{fig:mw-structure}.
Second, high performance is crucial, in particular for simulations in multiple space dimensions.
\revAK{Because the dependency chain is resolved at compile time, the compiler has full
	visibility into the instantiated configuration and can apply inlining and specialisation that
	runtime polymorphism would prohibit, at no runtime cost.}

\revAK{\cref{fig:mw-structure} organises the library into four layers, colour-coded by role.
	The grey nodes at the base form the geometry layer: \texttt{Quad} and \texttt{Grid} define the
	mesh, \texttt{Connectivity} encodes the cell-neighbour graph, and \texttt{Basis} together with
	\texttt{Space} fix the approximation space.
	The yellow nodes encode the physics: \texttt{Pde} supplies the PDE system and is the parent
	of \texttt{Numflux}, which discretises the Riemann problem at cell interfaces,
	\texttt{Divergence}, which assembles the discrete flux operator, a \texttt{Source} term, a slope
	\texttt{Limiter}, and boundary conditions \texttt{Bc}.
	The blue nodes \texttt{Data} and \texttt{Solution} manage the solution state.
	At the top, the green algorithm layer assembles these into a complete method:
	\texttt{Space\_discr} builds the spatial operator from geometry, physics, and solution;
	\texttt{Time\_discr} and \texttt{Weak\_formulation} close the time loop; \texttt{MRA} drives
	adaptive mesh refinement.
	The \texttt{MultiWave} entry point (purple) is the only symbol the user
	instantiates directly.}

\revAK{This layered structure is the concrete realisation of the framework's modularity. Each node represents one independently replaceable policy, connected only through the template interfaces shown. Replacing the PDE system, switching the Riemann solver, or adopting a different time-stepping strategy requires only a new class at the corresponding node. All remaining layers recompile unchanged.}

\revAK{\cref{fig:mw-laplace,fig:mw-pathcons,fig:mw-swe} show three such experiments. In each case only the orange-highlighted nodes are exchanged, while the rest of the graph stays unchanged.}

\revAK{The flexibility extends beyond the choice of algorithmic components to the underlying data representation itself. As part of ongoing work on high-dimensional stochastic problems, the standard DG coefficients are replaced by tensors. Since the framework interacts with coefficients exclusively through overloaded operators, this substitution requires only a minimal modification to the data layer, leaving all algorithmic components unchanged.}


\subsection{Structure of the discretised weak formulation}
\label{sec:struct-discr-weak}
In this section we focus on the implementation details of the modal Runge--Kutta DG (RKDG) method~\eqref{eq:weak-approx-cl-semi-discr}.  Our ansatz makes as few assumptions as possible  and accurately follows  the mathematical structure of the \revASSS{semi-discrete weak formulation~(\ref{eq:weak-approx-cl-semi-discr})}. This  allows for high code reusability of the framework as it provides building blocks for the construction of a variety of adaptive numerical methods.

\textsl{MultiWave} implements the  \revASSS{semi-discrete weak formulation~(\ref{eq:weak-approx-cl-semi-discr})} in time-marching fashion, cf. \cref{alg:general-adaptive-method}, evolving data stored in space-time slabs:
\begin{lstlisting}[style=mystyle]%,float, label={lst:mycode}]
template <typename TSol>
struct Step_slab {
  double t, dt; // current time and current timestep size
  std::shared_ptr<TSol> u_old, u_new;
};
\end{lstlisting}
The split time and space discretisation~\eqref{eq:weak-approx-cl-semi-discr} is provided in the \texttt{Split\_semi} class.
\begin{lstlisting}[style=mystyle]
template <typename TOptions,
          template <typename> typename TTime_discr,
          template <typename> typename TSpace_discr>
class Split_semi{
  /* .. type definitions .. */

  struct Post_proc{
    TOptions::Limiter limiter;
    void apply(auto &u){
      limiter.apply(u);
    }
  };

  TTime_discr<TOptions> temporal;
  TSpace_discr<TOptions> spatial;
  TOptions::Bc bc;
  Post_proc post;

  void step(Step_slab<Solution> slab) {
    temporal.step(slab,
      [&](double t, auto& u, auto& Rs) {  // RHS lambda
        spatial.rhs(u, Rs);
        bc.apply(t, u, Rs);
      },
      [&](auto& u) {                   // post-step lambda
        post.apply(u);                 // limiter
      });
  }
};  
\end{lstlisting}
The space-time slab is updated by the method \texttt{step} that calls the ODE-solver \texttt{temporal.step} with the right-hand side \texttt{spatial.rhs}. In addition, a post-processing lambda-function is  passed as the third argument and is called after each RK-stage. Consecutive calls of \texttt{Split\_semi::step} evolve the solution in the space-time slab, swapping \texttt{u\_old} and \texttt{u\_new} and then updating \texttt{u\_new} until the final time is reached in the main loop of the code, cf. \cref{alg:general-adaptive-method}.
The \texttt{TOptions} template parameter expects a class that stores all compile-time settings  set by the user and will be specified in what follows.

For the temporal discretisations standard ODE-solvers are employed. In our case strong stability preserving RK methods~\cite{gottlieb2001strong} are wrapped in the \texttt{SSPRK} class. The ODE-solvers are applied to evolve each solution coefficient $u_{\lambda, \vec{i}}(t)$ while bulk-processing the entire solution map. Moreover, the \texttt{SSPRK} class provides caching of the stage-storage maps (that are subject to change in each time-step due to adaptivity) to minimise memory allocations.  The \texttt{SSPRK} class is passed as the \texttt{TTime\_discr} parameter to \texttt{Split\_semi}. The spatial  right-hand side operator in~\eqref{eq:weak-approx-cl-semi-discr}, $\mathcal{R}_{DG}$,  is implemented in the method \texttt{rhs} of the \texttt{Modal\_dg} class. In order to construct an RKDG discretisation scheme, \texttt{Modal\_dg} instantiates the template  parameter \texttt{TSpace\_discr} in \texttt{Split\_semi}.

Next, we elaborate on the implementation of the spatial discretisation \texttt{Modal\_dg}. The operator $\mathcal{R}_{DG}$ is composed of $L^2$-inner products of nonlinear expressions of the solution coefficients and comprises two types of discretisations of terms in the \revASSS{semi-discrete system~(\ref{eq:weak-approx-cl-semi-discr})}:  expressions that \revASSS{approximate terms containing} the solution $\vec{u}_h$ only and expressions that \revASSS{approximate} spatial first-order  partial derivatives of $\vec{u}_h$. We denote the class that implements the first type  by \texttt{Source} and the class for the second type  by \texttt{Divergence}.  In order to evaluate \revASSS{$\mathcal{R}_{DG}(\vec{u}_h )$} for a single cell $V_\lambda$, two types of $L^2$-inner products arising from the Gauss theorem have to be approximated: volume integrals over the cell $V_\lambda$ and surface integrals over $\partial V_\lambda$. Since the DG-approximation allows jumps on the cell interfaces, the surface integrals depend on values of both cells adjacent to the surface. Thus, both classes \texttt{Divergence} and \texttt{Source} are required to implement two methods: \texttt{vol\_coeff} that evaluates the quadrature of volume integrals for all coefficients in the cell and \texttt{surf\_coeff} that evaluates the quadrature of surface integrals. Although typically there are no surface terms in the discretisation of the source, we keep this structure consistent for all components of \texttt{Modal\_dg}. In addition to the surface integrals on the inner cell edges, boundary conditions are implemented in the \texttt{bc} class's method \texttt{apply}. \Cref{alg:dg-rhs} summarises the discretisation of $\mathcal{R}_{DG}$.
\begin{algorithm}[H]
	\caption{\texttt{Modal\_dg::rhs(slab)}}\label{alg:dg-rhs}
	\begin{algorithmic}[1]
		\State  $\mathbf{R} \gets $ empty map $\lambda \mapsto \vec{u}_\lambda$
		\For{each cell $V_{\lambda}$}
		\State $\mathbf{R}_\lambda \gets \texttt{Divergence::vol\_coeff}(\vec{u}_\lambda)
			+ \texttt{Source::vol\_coeff}(\vec{u}_\lambda)$
		\EndFor
		\For{each interior edge $e$ between cells $V_{\lambda_1}$, $V_{\lambda_2}$}
		\State $(\vec{w}_{L}, \vec{w}_{R}) \gets
			\texttt{Div::surf\_coeff}(e, \vec{u}_{\lambda_1}, \vec{u}_{\lambda_2})
			+ \texttt{Source::surf\_coeff}(e, \vec{u}_{\lambda_1}, \vec{u}_{\lambda_2})$
		\State $\mathbf{R}_{\lambda_1} \mathrel{-}= \vec{w}_{L} $
		\State $\mathbf{R}_{\lambda_2} \mathrel{+}= \vec{w}_{R}$
		\EndFor
		\State $\texttt{bc.apply}{(t,\,\vec{u},\,\mathbf{R})}$
		\State \Return $\mathbf{R}$
	\end{algorithmic}
\end{algorithm}
For the sake of flexibility and performance, we employ the \emph{curiously recurring template pattern} (CRTP) and derive the \texttt{Modal\_dg} class from template parameters \texttt{Divergence} and \texttt{Source}.
\begin{lstlisting}[style=mystyle]
  template<typename TOptions,
           typename TDivergence,
           typename TSource = Zero_or_source<TOptions, Source_DG>>
  struct Modal_dg
      : public TDivergence,  public TSource
  {
    void rhs(const auto& solution, auto& Rs);
  };
\end{lstlisting}
While the \texttt{TDivergence} parameter is required,  the parameter \texttt{TSource}  is  instantiated by the \texttt{Zero\_or\_source} constant expression function as either \texttt{Zero} or as \texttt{Source\_dg} depending on the settings in \texttt{TOptions}. The \texttt{Zero}  class implements the \texttt{vol\_coeff} and \texttt{surf\_coeff} methods as empty functions that cause zero runtime costs. This approach allows the user to extend the \texttt{Modal\_dg} class by adding discretised operators that can be switched on and off on demand.

The template parameter \texttt{TSource}  is controlled by the availability of a source term implementation inside the \texttt{Pde} class (specified in the \texttt{TOptions} instantiation) and the source policy setting.
\begin{lstlisting}[style=mystyle]
  struct MyOptions {
    static constexpr int DIM = 2;               // two dimensions in space
    static constexpr int P_DIM = 3;             // order P=3 polynomials

    using Pde            = Euler<DIM>;      
    using Num_flux       = LLF;             // local Lax--Friedrichs num. flux
    using Source_policy  = Zero_policy;     // no source discretisation

    // ... types derived from user defined aliases
  };
\end{lstlisting}
The \texttt{TDivergence} parameter expects the \texttt{Pde} class to implement the \texttt{inv\_flux} method (i.e. the flux function $\vec{f}$). Since the flux function $\vec{f}$ is approximated by a numerical flux $\widehat{\vec{f}}$ at the cell interface, the \texttt{TDivergence} class requires the \texttt{Num\_flux} alias to be defined in the \texttt{TOptions} as a non-\texttt{Zero} class. For instance, in the above  code \texttt{Num\_flux} is set to a class implementing the local Lax--Friedrichs flux. The handling of \texttt{TSource} (and any additional base classes of \texttt{Modal\_DG}) is similar. However, the user may set the discretisation policy in the \texttt{TOptions} parameter to \texttt{Zero}. This (temporarily) turns off the \texttt{TSource} computation despite the presence of the \texttt{source} method in the \texttt{PDE} class. Thus, the user can convert a balance law approximation into a conservation law by changing one line of code, reusing existing \texttt{PDE} implementations. The compile-time switch between real implementation and a zero operation is implemented in the following way:
\begin{lstlisting}[style=mystyle]
  template<typename  Opt, template <typename> typename Impl>
  using Zero_or_source = std::conditional_t<
      has_source_policy<Opt> && Opt::Pde::has_source,
      Impl<Opt>,    // Source_dg discretisation
      Zero_s<Opt>   // empty zero-op
    >;
\end{lstlisting}
Finally, if the selected \texttt{PDE}-class does not declare the required method while the corresponding policy is set non-\texttt{Zero}, a \texttt{static\_assert} is triggered at compile time. The discretisation policy validation rules are summarised in \cref{tab:policy-validation}.
\begin{table}[!h]
	\begin{center}
		\begin{tabular}{cc|l}
			Method implemented in PDE & Policy set         & Result                       \\\hline
			\checkmark                & \checkmark         & real implementation          \\
			\checkmark                & \text{\sffamily X} & Zero (silently skipped)      \\
			\text{\sffamily X}        & \checkmark         & \texttt{static\_assert} FAIL \\
			\text{\sffamily X}        & \text{\sffamily X} & Zero                         \\
		\end{tabular}
	\end{center}
	\caption{Discretisation policy validation rules.}
	\label{tab:policy-validation}
\end{table}
Thus, the user has some flexibility when setting up a numerical method while the class template logic implementation remains hidden. In \revASSS{the subsequent} section we show how to add a second-order spatial derivative term to an existing DG method.

\subsection{Discretisation of the divergence-terms}
\label{sec:discr-diverg-terms}

The \texttt{Divergence\_dg} class implements the evaluation of the inner products
\revASSS{\[
		\int_{V_{\lambda}}\sum_{j=1}^d  \vec{f}_j(\vec{u}_h) \cdot\partial_{\vec{x}_j} \varphi_{\lambda,\vec{i}} d\vec{x}\quad  \text{and}\quad  \int_{\partial V_{\lambda}} \widehat{\vec{f}} \left(\vec{u}_h^+, \vec{u}_h^-, \vec{n}_\lambda\right) \varphi_{\lambda,\vec{i}} dS,
	\]
}
that occur in the semi-discrete weak formulation~\eqref{eq:weak-approx-cl-semi-discr}. It is passed as the \texttt{TDivergence} template  parameter to the \texttt{Modal\_dg} class. Internally it utilises  \texttt{inner\_prod} methods of the \texttt{DG\_space} class as follows
\begin{lstlisting}[style=mystyle]
result = 0;
const auto indices = std::views::iota(0u, N); // N - number of quadrature points
for (int d = 0; d < DIM; ++d){ 
  std::for_each(std::execution::unseq,
                indices.begin(), indices.end(), [&](unsigned int j) {
    buffer.vol_fluxes[j] =
        flux(j, buffer.vol_values[j], positive_normals[d].vec);
  });
  
  result += space.inner_prod(buffer.vol_fluxes,
                             buffer.basis_vol_der[d], lambda) / dx[d];
}
\end{lstlisting}
A mutable \texttt{buffer} that stores the values of the flux function at each quadrature point as well as the values of the gradient of the basis $\nabla\varphi_{\lambda, \vec{i}}$ is maintained in the class. Note that the computation is performed on a reference cell $[0,1]^d$ such that the evaluation of the gradient of the basis is required only once. Moreover, during the computation of the flux values, the compiler is permitted to generate SIMD instructions via the execution policy \texttt{std::execution::unseq}.

The evaluation of the surface inner products is performed in a similar manner:
\begin{lstlisting}[style=mystyle]
  up_L = space.inner_prod(buf.surf_fluxes_L, buf.basis_surf_L, edge); // w_L
  up_R = space.inner_prod(buf.surf_fluxes_R, buf.basis_surf_R, edge); // w_R
\end{lstlisting}
At return the resulting coefficients are added to the cell coefficients with signs according to the flux direction, cf. \cref{alg:dg-rhs}, lines 7 and 8.

Finally, the last component that is missing is the actual computation of the inner products.
To this end,  we approximate the volume inner products, \revASS{$\langle \vec{F},\varphi_{\lambda, i} \rangle_{V_{\lambda}}  = \int_V\sum_{j=1}^d  \vec{f}_j\,\partial_{\vec{x}_j}\varphi_{\lambda, i} d\vec{x}$}, by quadrature rules, typically using Gauss quadrature with nodes \revASS{$\vec{\zeta}^{\text{ref}}_q$} and weights $w_q$ where $q\in \{1, \ldots, N_q\}$ on the reference cell $[0,1]^d$. Both the nodes and the weights are computed once and then transformed to each physical cell $V_{\lambda}$ by an affine transformation $\phi_{\lambda}\, \colon\, \mathbb{R}^d \to \mathbb{R}^d$. By the change of variables of the integral the integrand is multiplied by $|J_{\lambda}|$, the determinant of the Jacobian matrix of $\phi_{\lambda}$.
\revASS{\[
		\langle \vec{F},\varphi_{\lambda, i} \rangle_{V_{\lambda}}
		\approx \sum_{j=1}^d\sum_{q=1}^{N_q} \vec{f}_j(\vec{u}_h(\vec{\zeta}_q))\, \partial_{\vec{x}_j}\varphi_{\lambda, i}(\vec{x}_q)\;w_q\;|J_{\lambda}|,
		\qquad \vec{\zeta}_q \coloneq \phi_{\lambda}( \vec{\zeta}^{\text{ref}}_{q}),
	\]}
\begin{lstlisting}[style=mystyle, mathescape=true]
  template<typename Tf_elems, integrable_val Tv>
    requires requires(Tf_elems x, Tv y) { x * y; }
  auto
  Space<D,P,TBasis>::inner_prod(
      const auto& f_vec,
      const auto& g_vec,
      const auto& lambda) const{                     // cell index                       
    return dot(f_vec * g_vec,                        // (1) $\vec{f}_q \cdot \varphi_q(\vec{x}_q)$
                   quadrature.weights)               // (2) dot with $w_q$  
                     * quadrature.deref_determinant( // (3) multiply by $|J_{\lambda}|$
                         grid.dom(lambda));          // cell geometry
    }
\end{lstlisting}
The \texttt{integrable\_val} concept ensures that the instantiation of \texttt{Tv} allows basic vector space properties, i.e. multiplication by scalars and summation of its objects:
\begin{lstlisting}[style=mystyle]
  template <typename T>
  concept integrable_val = requires(T y) {
    { y + y } -> std::convertible_to<T>;
    { 1.0 * y } -> std::convertible_to<T>;
  };
\end{lstlisting}

Inner products on surfaces $\Gamma = \overline{V_{\lambda_1}} \cap \overline{V_{\lambda_2}}$, i.e. over an edge $e$,
\revASS{\[
\langle \widehat{\vec{f}}, \varphi_j \rangle_e \coloneq  \int_{\Gamma} \widehat{\vec{f}}\,\varphi_j dS
\approx \sum_{q=1}^{N_q^s} \widehat{\vec{f}}(\vec{u}_h(\vec{\zeta}_q^s))\;\varphi_j(x_q^s)\;w_q^s\;|J_e^s|,
\qquad \vec{\zeta}^s_q \coloneq \phi_{\Gamma}(\vec{\zeta}^{\text{ref}}_q),
\]}
are handled in a similar way, except for \revASS{$\vec{\zeta}_q^{\text{ref}}\in \partial([0,1]^d)$} and the spatial dimension one, where the surface is a single point and, thus, \texttt{inner\_prod} is set to the product of the first two elements. On edges of cells on different levels (i.e. adjacent to hanging nodes), the values on the edge of the finer cell are obtained from pre-computed values at quadrature points while the values of the coarser cell are computed online by directly evaluating the basis function on the sub-cell edge.
\begin{lstlisting}[style=mystyle]
  template<typename Tf_elems, integrable_val Tv>
    requires requires(Tf_elems x, Tv y) { x * y; }
  auto
  Space<D,P,TBasis>::inner_prod(
      const auto& f_vec,
      const auto& g_vec,
      const Edge& e) const{
    if constexpr (D == 1) {
      return (f_vec[0] * g_vec[0]);        // 1-D: single point, no weights
    } else {
      // take the finer cell's domain for the Jacobian
      const auto& fine_cell =
          left(e).level < right(e).level ? right(e) : left(e);
      return dot(f_vec * g_vec, quadrature.surface.weights)
               * quadrature.surf_deref_determinant(grid.dom(fine_cell),
                                                   dir(e)));
   }
 }
\end{lstlisting}

\begin{rem}
	Due to the modular structure of the discrete weak formulation~\eqref{eq:weak-approx-cl-semi-discr} implementation, \textsl{MultiWave} is not restricted to this particular DG method. The DG method is our workhorse due to its general formulation that resembles the functional analytic theory. However, further approximation methods like CWENO~\cite{craveroCWENOUniformlyAccurate2018} or Lax-Wendroff-type schemes~\cite{du2018hermite} can be easily added.
\end{rem}

\subsection{Multiresolution analysis}
In this section we describe the implementation of the multiresolution-based grid adaptation in \textsl{MultiWave}. \revAK{The central class is \texttt{Multiscale}, which provides the multiscale transformation and its inverse as well as the coarsening and refinement logic.} As in \cref{sec:struct-discr-weak}, the design closely mirrors the mathematical structure of the MRA introduced in \cref{sec:grid-adaptivity}. For further details, we refer to \cite{gerhardWaveletFreeApproachMultiresolutionBased2021, gerhard2016adaptive}

To incorporate multiresolution-based grid adaptation into the DG solver, a refinement step is performed before each time step to anticipate the propagation of significant features, followed by a coarsening step to discard cells with non-significant detail coefficients. An adaptation step is summarised in \cref{alg:adaptation}. \revAK{The multiscale transformation~\eqref{eq:loc-ms-decomposition} decomposes} the solution into its coarsest-level representation and the associated detail coefficients $\mathcal{D}$ across all refinement levels. The detail coefficients $\mathcal{D}$ are then passed to either Harten's prediction strategy~\cite{harten1995multiresolution} for refinement or hard thresholding for coarsening, cf.~\cref{sec:grid-adaptivity}.
Both strategies determine which cells should remain active, but rather than modifying $\mathcal{D}$ directly, significant cells are collected into an auxiliary set $\mathcal T$. This provides a unified interface for both states: in the coarsening step, $\mathcal T$ contains only the significant cells identified by hard thresholding \eqref{eq:significant-detail}, while in the refinement step it contains all currently active details and is extended by Harten's prediction strategy. Once $\mathcal T$ is determined, the tree property is enforced by \texttt{generate\_tree}, which adds any missing ancestors to ensure a valid adaptive index set. The detail coefficients $\mathcal D$ are then synchronised with $\mathcal T$ via \texttt{sync\_d\_with\_t}, initialising details of newly included cells to zero and discarding the details of cells no longer in $\mathcal T$. The resulting set $\mathcal D$ constitutes the sparse approximation~\eqref{eq:sparse-approx}, from which the inverse multiscale transformation reconstructs the solution on the new adaptive grid \revakk{$\mathcal G_{L,\varepsilon}$, whose index set $\mathcal I_{L,\varepsilon}$ collects the leaves of $\mathcal T$, yielding the single-scale representation~\eqref{eq:sparse-approx-dg}, cf.~\cref{sec:grid-adaptivity}}. \revAK{The implementation of both the multiscale transformation and its inverse is detailed in the following section.}

\begin{algorithm}[H]
	\caption{\texttt{multiscale::adaptation(sol, state)}}
	\label{alg:adaptation}
	\begin{algorithmic}[1]
		\Require solution $\mathcal{U}$, state $\in \{\texttt{refinement},\, \texttt{coarsening}\}$
		\State $\textsc{mst}(\mathcal{U},\, 0,\, L)$ \Comment{multiscale transformation~\eqref{eq:loc-ms-decomposition}}
		\If{state $= \texttt{refinement}$}
		\State $\mathcal{T} \gets \mathcal{D}$
		\State $\textsc{predict\_multiscale}(\mathcal{U},\, 0,\, L)$
		\ElsIf{state $= \texttt{coarsening}$}
		\State $\textsc{hard\_thresholding}(0,\, L)$
		\EndIf
		\State $\textsc{generate\_tree}(0,\, L)$
		\State $\textsc{sync\_d\_with\_t}(0,\, L)$
		\State $\textsc{inverse\_mst}(\mathcal{U},\, 0,\, L)$
		\State $\textsc{update\_neighbours}(\mathcal{U})$
	\end{algorithmic}
\end{algorithm}

\begin{rem}
	\revakk{Apart from the prediction step, which relies on Harten's strategy~\cite{harten1995multiresolution,gerhardWaveletFreeApproachMultiresolutionBased2021} and thus on the propagation properties of the underlying PDE, the grid adaptation in \cref{alg:adaptation} operates on the solution space only.}
	As in \cref{sec:struct-discr-weak}, the adaptation method can be set in \texttt{MyOptions} by
	\begin{lstlisting}[style=mystyle]
  struct MyOptions {
    // ...
    using MRA  = Multiscale<U, TSpace>;     
    // ...
  };
  \end{lstlisting}
	For more advanced adaptation strategies tailored to specific applications, e.g.  shallow water equations \cite{Caviedes-Voullieme2020, Gerhard2015}, goal-based grid adaptation \cite{herty2024multiresolution} or wavelet-free grid adaptation \cite{gerhardWaveletFreeApproachMultiresolutionBased2021,Gerhard2017}, the \texttt{Mra} type can be set accordingly.
\end{rem}

\subsection{Two-scale transformation}
\label{sec:two-scale-trafo}
\revAK{The core operations in \cref{alg:adaptation} are the multiscale transformation $\textsc{mst}$ and its inverse $\textsc{inverse\_mst}$.} Both perform a local two-scale transformation applied level wise. We now derive these local operations and their matrix representation.

We express the space $\mathcal S_\ell$ in terms of its scaling functions $S_\ell = \mathrm{span}_{i\in \mathcal{I}^S_\ell} \phi_i$ and the space $\mathcal{W}_\ell$ in terms of multiwavelets $\mathcal{W}_\ell = \mathrm{span}_{i\in \mathcal{I}^W_\ell}\psi_i$, where $\mathcal{I}^S_\ell$ and $\mathcal{I}^W_\ell$ are index sets of the global degrees of freedom of $\mathcal{S}_\ell$ and $\mathcal{W}_\ell$, respectively.
To perform the multiscale transformation \eqref{eq:loc-ms-decomposition}, we need to express the basis functions of the DG space $\mathcal{S}_\ell$ and of its orthogonal complement space $\mathcal{W}_\ell$ in terms of scaling functions on level $\ell + 1$. Due to the locality of the basis functions it holds,
\begin{equation}
	\phi_{\lambda,i} = \sum_{\mu \in \mathcal M_\lambda}\sum_{k\in \mathcal{P}} \langle \phi_{\lambda,i},\phi_{\mu, k}\rangle_{V_\lambda} \phi_{\mu,k},\quad
	\psi_{\lambda,j} = \sum_{\mu \in \mathcal M_\lambda}\sum_{k\in \mathcal{P}} \langle \psi_{\lambda,j},\phi_{\mu, k}\rangle_{V_\lambda} \phi_{\mu,k},
	\label{eq:local-basis}
\end{equation}
where $\mathcal M_\lambda\subset \mathcal I_{\ell+1}$ is the refinement set of cell $V_\lambda$. The inner products in \eqref{eq:local-basis} are represented as matrices $\vec M_{\lambda, 0} \in \mathbb R^{|\mathcal{P}|\times |\mathcal{P}|\, |\mathcal M_\lambda|}$ and $\vec M_{\lambda, 1} \in \mathbb R^{|\mathcal{P}^*|\times |\mathcal{P}|\, |\mathcal M_\lambda|}$, respectively, where $|\mathcal{P}^*| := |\mathcal{P}|\,(|\mathcal M_\lambda| - 1)$ is the local degree of freedom for the complement space. Thus, the scaling functions and multiwavelets can be expressed in terms of functions on a finer level $\vec{\Phi}_\lambda = \vec M_{\lambda,0}\vec\Phi_{\mathcal M_\lambda}$ and $\vec\Psi_{\lambda} = \vec M_{\lambda,1}\vec\Phi_{\mathcal M_\lambda}$, and we define the mask-matrix
\begin{equation}
	\vec M_\lambda := \begin{pmatrix}
		\vec M_{\lambda,0} \\
		\vec M_{\lambda,1}
	\end{pmatrix}
	\in \mathbb R^{|\mathcal{P}|\,|\mathcal M_\lambda| \times |\mathcal{P}|\,|\mathcal M_\lambda| }.
	\label{eq:maskmatrix}
\end{equation}
Since the scaling functions and multiwavelets are orthogonal, the mask-matrix \eqref{eq:maskmatrix} is orthogonal. Thus, the scaling functions on a finer level can be represented as the sum of scaling functions and multiwavelets on the next coarser level
\begin{equation}
	\vec \Phi_{\mathcal M_\lambda} = \vec M^T_{\lambda,0}\vec\Phi_\lambda + \vec M^T_{\lambda,1}\vec \Psi_\lambda.
	\label{eq:ms-decomposition}
\end{equation}
Using \eqref{eq:ms-decomposition}, the local projections of \eqref{eq:orth-projections} are computed as

\begin{align*}
	u^\ell_\lambda = P_{\mathcal{S}_\ell}(u^{\ell+1}_\lambda) = (\vec M_{\lambda, 0}\vec u_{\mathcal M_\lambda})\cdot\vec\Phi_\lambda,\quad
	d^\ell_\lambda = P_{\mathcal{W}_\ell}(u^{\ell+1}_\lambda) = (\vec M_{\lambda, 1}\vec u_{\mathcal M_\lambda})\cdot\vec\Psi_\lambda.
\end{align*}
\revakk{Finally}, the DG coefficients of the two-scale transformation and the inverse two-scale transformation are determined by
\begin{align*}
	\begin{pmatrix}
		\vec u_\lambda \\
		\vec d_\lambda
	\end{pmatrix} = \vec M^T_\lambda \vec u_{\mathcal M_\lambda},\quad
	\vec u_{\mathcal M_\lambda} = \vec M_\lambda \begin{pmatrix}
		\vec u_\lambda \\
		\vec d_\lambda
	\end{pmatrix}.
\end{align*}
Given children $\vec u_{\mathcal M_\lambda}$ of the cell $V_\lambda$, we perform the two-scale transformation to obtain the detail and scaling coefficients on the coarser refinement level:

\begin{lstlisting}[style=mystyle, escapeinside={<@}{@>}]
template <int U, typename TSpace>
void Multiscale<U, TSpace>::two_scale(
    const std::array<StaticMatrix<U, NR_BASES>, NR_CHILDREN> &children, 
    auto &v, auto &d) const {
  //   children_flat[:, s*NR_BASES .. (s+1)*NR_BASES-1] = children[s]
  StaticMatrix<U, NR_BASES * NR_CHILDREN> children_flat;
  for (auto s = 0u; s < NR_CHILDREN; ++s)
    submatrix(children_flat, 0, s * NR_BASES, U, NR_BASES) =
        children[s];

  const StaticMatrix<U, NR_BASES * NR_CHILDREN> wavelet_flat =
      children_flat * transpose(mst_mat_ptr->fwd_transform);

  // col block 0 <@$\to$@> v, blocks 1..NR_CHILDREN-1 <@$\to$@> d[0..NR_CHILDREN-2]
  v = submatrix(wavelet_flat, 0, 0, U, NR_BASES);
  for (auto e = 1u; e < NR_CHILDREN; ++e)
    d[e - 1] = submatrix(wavelet_flat, 0, e * NR_BASES, U, NR_BASES);
}
\end{lstlisting}

Analogously, we obtain children of a cell $V_\lambda$ by performing the inverse two-scale transformation:

\begin{lstlisting}[style=mystyle, escapeinside={<@}{@>}]
template <int U, typename TSpace>
void multiscale<U, TSpace>::inv_two_scale(
    std::array<StaticMatrix<U, NR_BASES>, NUM_CHILDREN> &children, const auto &v,
    const auto &d) const {
  //   wavelet_flat[:, 0 .. NR_BASES-1]              = v
  //   wavelet_flat[:, (e+1)*NR_BASES .. (e+2)*NR_BASES-1] = d[e]
  StaticMatrix<U, NR_BASES * NR_CHILDREN> wavelet_flat;
  submatrix(wavelet_flat, 0, 0, U, NR_BASES) = v;

  for (auto e = 1u; e < NUM_CHILDREN; ++e)
    submatrix(wavelet_flat, 0, e * NR_BASES, U, NR_BASES) = d[e - 1];

  const StaticMatrix<U, NR_BASES * NR_CHILDREN> children_flat =
      wavelet_flat * transpose(mst_mat_ptr->inv_transform);

  for (auto s = 0u; s < NUM_CHILDREN; ++s)
    children[s] =
        submatrix(children_flat, 0, s * NR_BASES, U, NR_BASES);
}
\end{lstlisting}

\begin{rem}
	\revAK{Note that in \texttt{two\_scale} and \texttt{inv\_two\_scale}} we perform the transformations simultaneously for \texttt{U} variables. The mask-matrix and its inverse are given by \texttt{fwd\_transform} and \texttt{inv\_transform}. Choosing the Legendre polynomials as basis, these matrices can be pre-computed for all cells, since the orthogonality of the inner products \eqref{eq:local-basis} is invariant under affine transformation.
\end{rem}

\subsection{Adding new features}
\label{sec:extens-new-feat}
One of the objectives of \textsl{MultiWave} is to provide the user with flexibility to add or modify its features on any level of abstraction.
In the previous two sections we have shown how templates, \revAK{in particular the CRTP pattern}, are employed to generate a tight class structure that avoids any computational overhead.
The main rationale of this design choice, however, is the ability of the user to replace any component in \cref{fig:mw-structure} by a custom implementation (whether inheriting existing classes or implementing new ones) without touching the code base.
This allows for prototyping new adaptive numerical methods and easy incorporation of successful prototypes. In this section we present three examples of how this is achieved in practice.

\subsubsection*{\textbf{Adding \revASS{second-order} operators to RKDG}}
\begin{figure}[htbp]
	\centering
	\includegraphics[width=\textwidth]{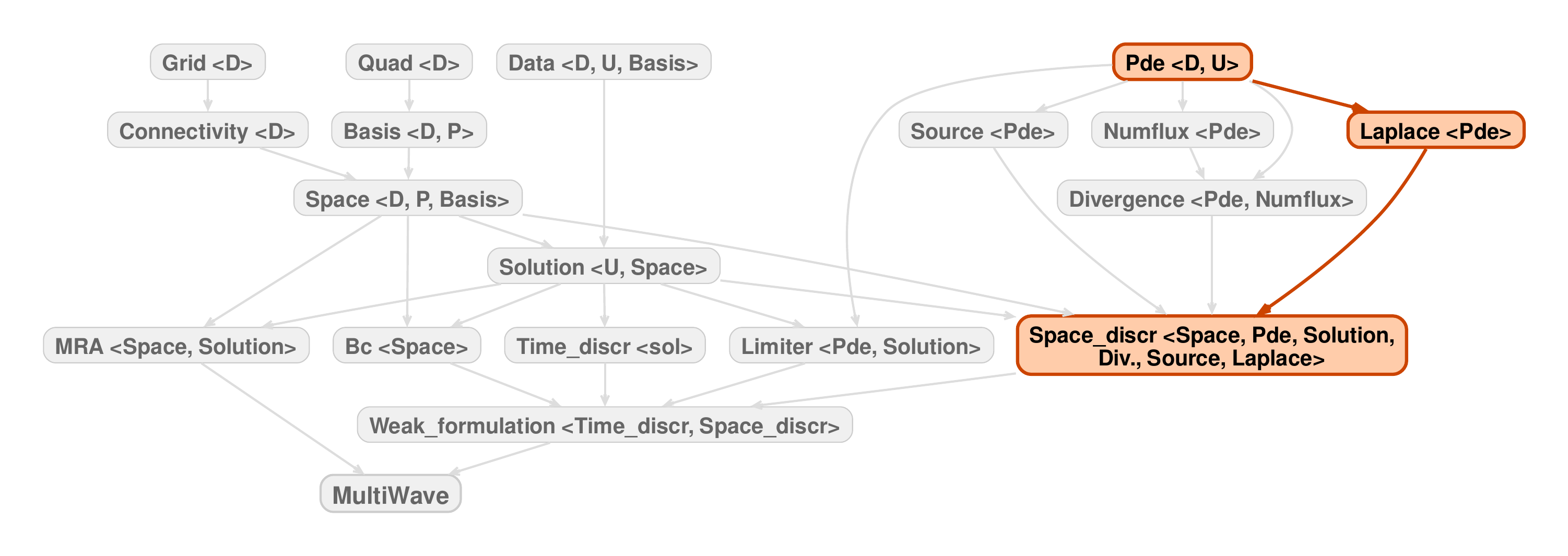}
	\caption{Only the \texttt{Pde} class and the new \texttt{Laplace} node (highlighted) need to be modified; all other parts of the hierarchy remains unchanged.}
	\label{fig:mw-laplace}
\end{figure}

The RKDG method can be applied to solve convection-dominated balance laws
\revASS{\begin{align*}
		\frac{\partial\vec{u}}{\partial t} + \nabla\cdot \vec{F}(\vec{u}) + {\color{red}\nabla\cdot \vec{H}(\vec{u}, \nabla \vec{u})}  + \vec{S}(\vec{x}, \vec{u}) = \vec{0}.
	\end{align*}
	with $\vec{H}\,\colon\, D \times \mathbb{R}^{m \times d} \to \mathbb{R}^{m \times d}$, by stabilising the viscous fluxes $\vec{H}$}, for instance, using the Bassi--Rebay method (BR2)~\cite{bassi1997high}. For details we also refer to~\cite{Gerhard2017} and the literature cited therein. For our purposes it is sufficient to extend the \cref{alg:dg-rhs} by adding
\begin{align*}
	\mathbf{R}_\lambda \gets \texttt{Divergence::vol\_coeff}(\lambda, \vec{u})
	 & + {\color{red}\texttt{Laplace::vol\_coeff}(\lambda, \vec{u})} \\
	 & + \texttt{Source::vol\_coeff}(\lambda, \vec{u})
\end{align*}
and
\begin{align*}
	(\vec{w}_{L}, \vec{w}_{R}) \gets
	\texttt{Divergence::surf\_coeff}(e, \vec{u})
	and
	 & + {\color{red}\texttt{Laplace::surf\_coeff}(e, \vec{u})} \\
	 & + \texttt{Source::surf\_coeff}(e, \vec{u}).
\end{align*}
In order to implement this, first, the \texttt{Pde} class is extended by a method \texttt{visc\_flux} that corresponds to the function $\vec{H}$. For example, the Navier--Stokes equations are obtained when deriving from the existing Euler equations class and implementing the viscous fluxes. Second, a \texttt{Laplace}  class together with the selection mechanism \texttt{Zero\_or\_laplace} is required. The \texttt{Laplace} class is added to the base classes of \texttt{Modal\_dg} providing the \texttt{vol\_coeff} and \texttt{surf\_coeff} functions. Note that all previously implemented systems remain functional, since in that case \texttt{Zero\_or\_laplace} returns the \texttt{Zero} type. The resulting modifications to the class hierarchy are illustrated in \cref{fig:mw-laplace}.
\begin{lstlisting}[style=mystyle, escapeinside={<@}{@>}]
  template<typename TOptions,
           typename TDivergence,
           typename <@{\color{red}\texttt{TLaplace}}@> = Zero_or_laplace<TOptions, Bassi_Rebay>, 
           typename TSource = Zero_or_source<TOptions, Source_dg>>
  struct Modal_dg
      : public TDivergence, public <@{\color{red}\texttt{TLaplace}}@>, public TSource
  {
    void rhs(const auto& solution, auto& Rs);
  };
\end{lstlisting}
Finally, the BR2 scheme itself is implemented as a discretisation policy and set in the compile-time options.
\begin{lstlisting}[style=mystyle, escapeinside={<@}{@>}]
  struct MyOptions {
    static constexpr int DIM = 2;
    static constexpr int P_DIM = 3;

    using Pde            = Navier_Stokes<DIM>;
    using Num_flux       = LLF;           
    using <@{\color{red}\texttt{Laplace\_policy}}@>   = Bassi_Rebay;

   // ...
  };
\end{lstlisting}

\subsubsection*{\textbf{Path-conservative DG}}
In multi-phase flow models, some quantities are not conserved across fluid contact boundaries. In this case, the underlying hyperbolic PDE includes non-conservative terms $\vec{G}$, i.e.
\revASS{\begin{align*}
		\frac{\partial\vec{u}}{\partial t} + \nabla\cdot \vec{F}(\vec{u}) + {\color{red}\sum_{i=1}^d \vec{G}_i(\vec{u}) \partial_{x_i}\vec{u}}  + \vec{S}(\vec{x}, \vec{u}) = \vec{0}.
	\end{align*}}
Following the Dal Maso--LeFloch--Murat theory~\cite{dal1995definition}, the weak formulation for the non-conservative \revAK{product} reads
\[
	\sum_{i=1}^d\int_{V_\lambda} \vec{G}_i(\vec{u_h}) \partial_{x_i}\vec{u}_h \varphi_{\lambda, \vec{i}} d\vec{x} + \int_{\partial V_{\lambda}} \widehat{\vec{g}} \left(\vec{u}_h^+, \vec{u}_h^-, \vec{n}_\lambda\right) \varphi_{\lambda,\vec{i}} dS,
\]
with the path-conservative numerical flux
\[
	\revAK{
		\widehat{\vec{g}}\left(\vec{u}_h^+, \vec{u}_h^-, \vec{n}_\lambda\right)
		= \frac12 \frac{\varphi^-_{\lambda,\vec{i}} + \varphi^+_{\lambda,\vec{i}}}{\varphi^-_{\lambda,\vec{i}}} \int_0^1 \sum_{i=1}^d n_{\lambda,i}\,
		\vec{G}_i\!\left(\Phi(\tau, \vec{u}_h^-, \vec{u}_h^+)\right)
		\partial_{\tau}\Phi(\tau, \vec{u}_h^-, \vec{u}_h^+)\, d\tau
	}
\]

Here, $\Phi\, \colon\, [0,1] \times D \times D \to D$ is \revAK{a} path, i.e.~a homotopy in the first argument $\tau$ that is a part of the underlying physical model.

\begin{figure}[htbp]
	\centering
	\includegraphics[width=\textwidth]{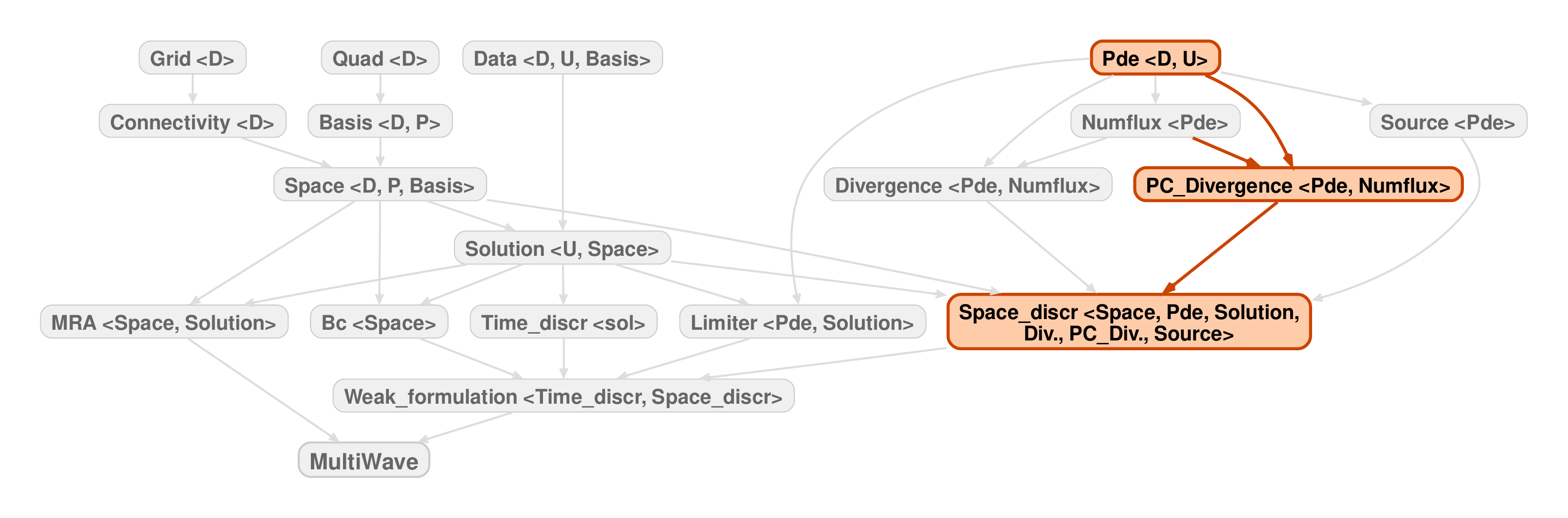}
	\caption{The new \texttt{PC\_Divergence} component (highlighted) is added alongside the existing \texttt{Divergence} term; \texttt{Space\_discr} is extended to accept both.}
	\label{fig:mw-pathcons}
\end{figure}

The non-conservative terms are added by deriving a class \texttt{PC\_Divergence\_dg} from the \texttt{Divergence\_dg} and implementing additional inner products in the \texttt{vol\_coeff} and \texttt{surf\_coeff} methods.  In particular, the path-conservative surface integrals are computed with the normal directed outwards of the cell, i.e. the path-conservative fluxes have different signs than the conservative ones and are implemented as follows:
\begin{lstlisting}[style=mystyle]
  // nc_nf = non-conservative path-integral 
  surf_fluxes_L = nf + 0.5 * nc_nf;
  surf_fluxes_R = nf - 0.5 * nc_nf;
\end{lstlisting}
where \texttt{nf} stores the coefficients from the conservative numerical flux $\widehat{\vec{f}}$, see \cref{fig:mw-pathcons}.

\subsubsection*{\textbf{MRA for shallow water equations}}
The shallow water equations read
\begin{align*}
	\frac{\partial}{\partial t}\begin{pmatrix}
		h \\ h\vec{v}
	\end{pmatrix} + \nabla\cdot
	\begin{pmatrix}
		h \vec v \\ h \vec v \otimes \vec v + \frac{1}{2}gh^2\vec I
	\end{pmatrix} =
	{\color{red}\begin{pmatrix}
		\vec 0 \\ -gh\nabla b
	\end{pmatrix}},
\end{align*}
where $h > 0$ is the water depth, $\vec{v} \in \mathbb{R}^d$ is the velocity vector, $b \in \mathbb{R}$ is the bottom topography, and $g$ is the gravitational constant. Although it is a hyperbolic balance law of the form~\eqref{eq:general-pde}, two important issues have to be addressed. First, we have to ensure that the numerical scheme is \emph{well-balanced}, i.e., steady states are preserved on non-constant bottom topography. Second, positivity of the water height must be preserved, especially at wet-dry interfaces, see~\cite{Gerhard2017, Gerhard2015, Caviedes-Voullieme2020} and references therein.

The well-balancing property can be achieved by modifying the numerical flux~\cite{Xing2014} with
\begin{align}
	\widehat{\vec f}_{\text{wb}}^+(\vec u^+_h, \vec u^-_h, b^+,b^-,\vec n_\lambda) := \widehat{\vec f}(\vec u^+_h, \vec u^-_h, \vec n_\lambda)
	+ {\color{red}\frac{g}{2}\left(\left(h^+\right)^2 - \left(h_m^+\right)^2\right)
	\begin{pmatrix}
		0 \\ \vec n_\lambda
	\end{pmatrix}}, \\
	\intertext{for the right cell, and}
	\widehat{\vec f}_{\text{wb}}^-(\vec u^+_h, \vec u^-_h, b^+,b^-,\vec n_\lambda) := \widehat{\vec f}(\vec u^+_h, \vec u^-_h, \vec n_\lambda)
	+ {\color{red}\frac{g}{2}\left(\left(h^-\right)^2 - \left(h_m^-\right)^2\right)\begin{pmatrix}
		                                                                               0 \\ \vec n_\lambda \end{pmatrix}},~
	\label{eq:well-balance-flux}
\end{align}
for the left cell, where $\vec u_m^\pm := (h^\pm_m, (h_m^\pm\vec v^\pm))^T$ with $h^\pm_m := \max\left(0, h^\pm + b^\pm - \max\left(b^+,b^-\right)\right)$.

\begin{figure}[htbp]
	\centering
	\includegraphics[width=\textwidth]{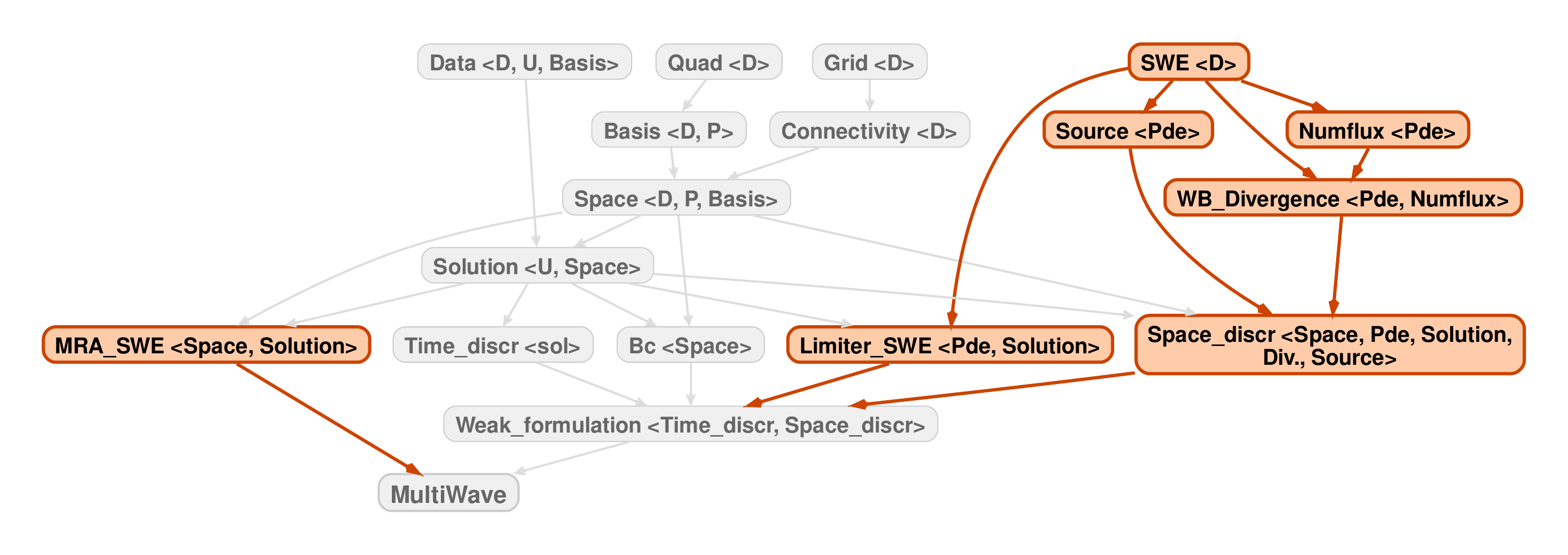}
	\caption{The highlighted components --- \texttt{SWE}, well-balanced numerical flux, \texttt{WB\_Divergence}, \texttt{Source}, \texttt{Limiter\_SWE}, and \texttt{MRA\_SWE} --- are the only parts of the hierarchy that require specialisation.}
	\label{fig:mw-swe}
\end{figure}

To preserve positivity, a positivity-preserving limiter has to be employed, and the multiresolution-based grid adaptation has to be performed on $h + b$ rather than $h$, with additional positivity corrections applied. For the details of these algorithms, we refer to \cite{Gerhard2015, Gerhard2017}. These modifications are implemented in \textsl{MultiWave} (see \cref{fig:mw-swe}) and can be activated by setting the compile-time options

\begin{lstlisting}[style=mystyle, escapeinside={<@}{@>}]
  struct MyOptions {
    static constexpr DIM = 2;
    static constexpr P_DIM = 3;

    using Pde            = SWE<DIM>;        
    using Num_flux       = <@{\color{red}\texttt{LLF\_SWE}}@>;             // well-balanced flux
    using <@{\color{red}\texttt{Source\_policy}}@>    = Source_DG;
    using MRA            = <@{\color{red}\texttt{Multiscale\_SWE<U, TSpace>}}@>;

   // ...
  };
\end{lstlisting}

The following section addresses the data structures and parallelisation strategy that allow \textsl{MultiWave} to scale to large distributed-memory systems.

\section{Performance}
\label{sec:performance}
Many practical problems for hyperbolic balance laws require either high spatial resolution or large computational domains, making performance a critical aspect. In multiple space dimensions, the number of degrees of freedom grows exponentially with the refinement level, and even with the compression provided by the multiresolution-based adaptivity, efficient use of computational resources is indispensable. We address this on two levels: first, through careful selection of data structures and algorithms that exploit cache locality and SIMD vectorisation to maximise single-node throughput; second, through a distributed-memory parallelisation strategy based on MPI that allows the framework to scale to large clusters. In this section, we discuss both aspects in detail.

\subsection{Data structures}

\label{sec:datatypes}
The three most consequential data structure choices in \textsl{MultiWave} concern the per-cell coefficient storage, the representation of the adaptive solution, and the neighbourhood connectivity graph.

\subsubsection*{\textbf{Cell-local data layout.}}
\label{sec:cell-local-data}
The DG coefficients $\vec{u}_\lambda:= (\vec{u}_{\lambda, \vec{i}})_{\vec{i}\in \mathcal P}\in \mathbb R^{m\times |\mathcal{P}|}$ of a single cell $V_\lambda$ form a matrix whose dimensions depend on the number of states and the basis $\Phi_h$. The size of $\Phi_h$, and hence of $\vec{u}_\lambda$, is entirely determined by the choice of basis, polynomial degree $p$, and spatial dimension $d$. Since these are fixed at compile time as template parameters, the size of the coefficient matrix is a compile-time constant. Consequently, the coefficient matrices are stored as static matrices, which avoids heap allocations and enables the compiler to emit fully unrolled SIMD-vectorised loops over the polynomial modes. This \revAK{leads to a significant acceleration of} the small, dense linear algebra operations that dominate the runtime, namely volume integrals, surface integrals, and \revAK{the two-scale transformations.} All linear algebra operations are performed by means of the Blaze library\footnote{\url{https://bitbucket.org/blaze-lib/blaze/src/master/}}~\cite{Iglberger2012, Iglberger2012a}
which is tuned for high-performance applications. Blaze provides implementations of BLAS-like operations tuned for small matrices leveraging SIMD instructions, and additionally simplifies mathematical expressions through smart expression templates, allowing one to write mathematically precise code without sacrificing performance.

\begin{rem}
	To avoid repeated heap allocations in the innermost computational loops, all intermediate quantities required during the evaluation of the weak formulation, including solution values, fluxes, and basis function evaluations at quadrature points, are stored in a persistent, pre-allocated buffer. This buffer is initialised once and reused across all cells and time steps, eliminating dynamic memory overhead and improving cache locality. The same strategy is applied throughout all modules, including the slope limiter, the connectivity, and the MRA.
\end{rem}

\subsubsection*{\textbf{Adaptive solution representation.}}
\label{sec:adaptive-solution-representation}
Each cell $V_\lambda$ on a Cartesian grid is uniquely identified by its \emph{levelmultiindex} (LMI), defined as
$\text{lmi}_\lambda := \bigl(\ell,  (\text{mi}_1,\dots,\text{mi}_d)\bigr) \in \mathbb N_0 \times \mathbb N_0^d$,
where $\ell$ is the refinement level and $(\text{mi}_1,\dots,\text{mi}_d)$ encodes the $\lambda$-th cell position within the computational domain. This representation directly reflects the hierarchical structure of the MRA and enables efficient computation of parent--child relationships. For instance, the parent of $\text{lmi}_\lambda$ is given by
$\bigl(\ell - 1, (\lfloor \text{mi}_1/2\rfloor, \dots, \lfloor\text{mi}_d/2\rfloor)\bigr)$
and its $2^d$ children are
$\bigl(\ell + 1, ( 2\cdot\text{mi}_1 +k_1, \dots, 2\cdot\text{mi}_d + k_d)\bigr)$ for $(k_1,\dots,k_d) \in \{0,1\}^d$. A 2D example is given in \cref{fig:lmi-example}. To reduce memory consumption, each LMI is packed into a single 64-bit unsigned integer by allocating a fixed number of bits for the level and for each spatial index. \revAK{Since both encoding and decoding of the LMI} reduce to bit-wise operations, this representation is cache-friendly and incurs negligible overhead. While this limits the maximum size of the computational grid, the resulting bounds are sufficient for most practical applications, as detailed in \cref{tab:lmi-encoding}. Since the LMI encodes both the refinement level and the cell position, the full geometry of each cell is uniquely defined from the LMI alone, without any additional geometric storage. \revAK{A similar cell-identifier concept for adaptive DG methods on 2D and 3D triangulations has been developed in~\cite{Brix2011a}.}

\begin{rem}
	Since the grid is Cartesian, all cells at level $\ell$ are geometrically congruent, with uniform side length $h_\ell = h_0 / 2^\ell$. An affine map from any cell $V_\lambda$ to the reference cell $[0,1]^d$ is therefore fully determined by its LMI. Reference quadrature points are fixed and shared across all cells; their physical counterparts are recovered on the fly via an affine map, see \cref{sec:discr-diverg-terms}.
	Cell-local DG operators on the reference cell depend on the refinement level only through the Jacobian factor $h_\ell^d$, and can thus be pre-computed once per level and reused across all cells.
\end{rem}

\begin{figure}[h]
	\noindent
	\begin{minipage}[t]{0.48\linewidth}
		\centering
		\begin{tikzpicture}[scale=5.5]
			\draw [thick] (0,0) -- (1,0) -- (1,1) -- (0,1) -- (0,0);
			\draw [thick] (0,0.5) -- (1,0.5);
			\draw [thick] (0.5,0) -- (0.5,1);
			\draw [thick] (0.5,0.75) -- (1,0.75);
			\draw [thick] (0.75,0.5) -- (0.75,1);
			\draw [thick] (0.625,0.5) -- (0.625,0.75);
			\draw [thick] (0.5,0.625) -- (0.75,0.625);
			\node at (0.25, 0.25)     {\textbf\footnotesize$(0,\!(0,\!0)\!)$};
			\node at (0.75, 0.25)     {\textbf\footnotesize$(0,\!(1,\!0)\!)$};
			\node at (0.25, 0.75)     {\textbf\footnotesize$(0,\!(0,\!1)\!)$};
			\node at (0.875, 0.625)   {\textbf\footnotesize$(1,\!(3,\!2)\!)$};
			\node at (0.625, 0.875)   {\textbf\footnotesize$(1,\!(2,\!3)\!)$};
			\node at (0.875, 0.875)   {\textbf\footnotesize$(1,\!(3,\!3)\!)$};
			\node at (0.5625, 0.5625) {\scalebox{0.5}{\textbf\tiny$(2,\!(4,\!4)\!)$}};
			\node at (0.6875, 0.5625) {\scalebox{0.5}{\textbf\tiny$(2,\!(5,\!4)\!)$}};
			\node at (0.5625, 0.6875) {\scalebox{0.5}{\textbf\tiny$(2,\!(4,\!5)\!)$}};
			\node at (0.6875, 0.6875) {\scalebox{0.5}{\textbf\tiny$(2,\!(5,\!5)\!)$}};
		\end{tikzpicture}
		\caption{LMIs on a 2D adaptive Cartesian grid with three refinement levels, showing successive refinement of cell $(0, (1,1))$.}
		\label{fig:lmi-example}
	\end{minipage}
	\hfill
	\begin{minipage}[b]{0.48\linewidth}
		\centering
		\begin{tabular}[t]{l|c|c}
			      & Max.\ level & Max.\ index per direction \\
			\hline
			$d=1$ & 63 (6 bit)  & $2^{58}-1$                \\
			$d=2$ & 63 (6 bit)  & $2^{29}-1$                \\
			$d=3$ & 15 (4 bit)  & $2^{20}-1$                \\
		\end{tabular}
		\captionof{table}{Bit allocation of a 64-bit LMI encoding per spatial dimension.}
		\label{tab:lmi-encoding}
	\end{minipage}
\end{figure}

Due to the dynamic grid adaptation, the solution cannot be represented as a single flat vector efficiently. The MRA requires fast access to parents and children by direct arithmetic operations on the LMI, and the two-scale transforms require rebuilding the solution and detail coefficients level by level, making a flat data layout impractical.

Hash maps with the LMI as key and the coefficient matrix as value are a natural choice in order to guarantee search, insertion and deletion in $\mathcal O(1)$ time. However, standard hash maps are not cache-friendly, leading to frequent cache misses when iterating over the solution.
To overcome this, we use the \texttt{ankerl::unordered\_dense}\footnote{\url{https://github.com/martinus/unordered_dense}} library, which is designed for high-performance applications. Unlike standard hash maps, it stores data contiguously in memory,  providing cache-efficient iteration. Fast lookups and insertions are also guaranteed, and although deletion of individual keys requires two lookups and is therefore slower, this overhead is negligible in practice. Since the two-scale transforms require access to all LMIs at a given refinement level, we store the solution as a vector of \texttt{unordered\_dense} hash maps, one per level. This allows for both fast per-cell operations and efficient iteration over all cells at a fixed refinement level.

\subsubsection*{\textbf{Neighbourhood connectivity.}}
\label{sec:neighbourhood-connec}
Computing the surface integrals in the semi-discrete formulation~\eqref{eq:weak-approx-cl-semi-discr} requires, for each cell $V_\lambda$, the values $\vec{u}_h^-$ of all face-adjacent neighbours. For efficiency reasons, the connectivity of the adaptive grid is stored explicitly, since neighbouring cells cannot be inferred by simple index arithmetic.
Here, the neighbourhood graph is represented as an adjacency list, implemented as an \texttt{ankerl::unordered\_dense} hash map keyed by the LMI. Each entry stores the LMI and outward normal of every face neighbour. This design yields $\mathcal O(1)$ neighbour lookups during flux evaluation.
Storing the outward normals avoids repeated geometric queries in the solver inner loop, at the cost of additional memory proportional to the number of adjacency entries.

A further advantage of this approach is that it does not require the 2:1 balance restriction, which mandates that neighbouring cells differ by at most one refinement level. \textsl{MultiWave} supports arbitrary level differences between neighbouring cells, which is particularly relevant for multiresolution-based grid adaptation in hyperbolic problems, where abrupt level transitions arise naturally near discontinuities. We refer to \cref{sec:discr-diverg-terms} for the treatment of non-conforming interfaces in the modal DG setting.

Since the adaptive grid changes after every time step, the connectivity graph is updated locally during refinement and coarsening, exploiting the tree structure of the MRA hierarchy to avoid a full rebuild. The coarsening procedure is described in \cref{alg:connec-coarsening}: all children $\mathcal{C}_\mu$ of the parent $V_\mu$ are removed from the graph, and their external neighbours are collected into $\mathcal N_\mu$ and assigned to $V_\mu$. The refinement procedure is given in \cref{alg:connec-refinement}: each child $V_\kappa$ inherits the neighbours of $V_\lambda$ that are geometrically face-adjacent to it. For each face direction, a neighbour $V_{\kappa'}$ is assigned to $V_\kappa$ if $V_{\kappa'}$ is coarser than the child level, or if its ancestor at the child level is the direct index-neighbour of $V_\kappa$. Sibling connections between children are added separately. Both operations update the graph locally in $\mathcal O(2^d\cdot k)$, where $k$ is the maximum number of face neighbours per cell. An illustration is given in \cref{fig:connec-update}.

\begin{algorithm}[ht]
	\caption{Coarsening of cell $V_\lambda$}
	\label{alg:connec-coarsening}
	\begin{algorithmic}[1]
		\Require{$V_\lambda \in \mathcal{I}_\ell$, $\ell > 0$, adjacency list $\mathcal{A}$}

		\State $V_{\mu} \gets \textsc{get\_parent}(V_\lambda)$
		\State $\mathcal{C}_{\mu} \gets \textsc{get\_child\_indices}(V_{\mu})$
		\State $\mathcal{N}_{\mu} \gets \emptyset$

		\For{each $V_\kappa \in \mathcal{C}_{\mu}$}
		\For{each $(V_{\kappa'},\, \vec{n}) \in \textsc{get\_neighbours}(\mathcal{A},\, V_\kappa)$}
		\If{$V_{\kappa'} \notin \mathcal{C}_{\mu}$}
		\Comment{No internal sibling connections}
		\State $\mathcal{N}_\mu \gets \mathcal{N}_{\mu} \cup \{(V_{\kappa'},\, \vec{n})\}$
		\EndIf
		\EndFor
		\State $\textsc{remove\_all\_neighbours}(\mathcal{A},\, V_\kappa)$
		\EndFor

		\State $\textsc{add\_neighbours}(\mathcal{A},\, V_{\mu},\, \mathcal{N}_{\mu})$
		\Comment{$V_\mu$ inherits external neighbours}
	\end{algorithmic}
\end{algorithm}

\begin{algorithm}[ht]
	\caption{Refinement of cell $V_\lambda$}
	\label{alg:connec-refinement}
	\begin{algorithmic}[1]
		\Require{$V_\lambda \in \mathcal{I}_\ell$, $\ell < L$, adjacency list $\mathcal{A}$}

		\State $\mathcal{C}_\lambda \gets \textsc{get\_child\_indices}(V_\lambda)$
		\State $\mathcal{N}_\lambda \gets \textsc{get\_neighbours}(\mathcal{A},\, V_\lambda)$
		\State $\mathcal{N}_\kappa \gets \emptyset$ for all $V_\kappa \in \mathcal{C}_\lambda$

		\For{each $\vec{n} \in \{\pm\vec{e}_1, \ldots, \pm\vec{e}_d\}$}
		\For{each $V_\kappa \in \textsc{get\_children\_with\_normal}(V_\lambda,\, \vec{n})$}
		\For{each $(V_{\kappa'},\, \vec{n}) \in \mathcal{N}_\lambda$ with $\vec{n}$-direction}
		\State $V_{\hat\kappa} \gets \textsc{get\_parent}(V_{\kappa'},\; \max(0,\, \ell_{\kappa'} - \ell_\kappa))$
		\If{$\ell_{\kappa'} \leq \ell_\lambda$ \textbf{ or } $\textsc{next\_index}(V_\kappa,\, \vec{n}) = V_{\hat\kappa}$}
		\State $\mathcal N_\kappa \gets \mathcal{N}_\kappa \cup \{(V_{\kappa'},\, \vec{n})\}$
		\EndIf
		\EndFor
		\EndFor
		\EndFor

		\For{each adjacent pair $(V_{\kappa_i},\, V_{\kappa_j}) \in \mathcal{C}_\lambda$} \Comment{Sibling connections}
		\State $\mathcal{N}_{\kappa_i} \gets \mathcal N_{\kappa_i} \cup \{(V_{\kappa_j},\, \vec{n}_{\kappa_i\kappa_j})\}$
		\EndFor

		\For{each $V_\kappa \in \mathcal{C}_\lambda$}
		\State $\textsc{add\_neighbours}(\mathcal{A},\, V_\kappa,\, \mathcal{N}_\kappa)$
		\EndFor

		\State $\textsc{remove\_all\_neighbours}(\mathcal{A},\, V_\lambda)$
	\end{algorithmic}
\end{algorithm}

\begin{figure}[h]
	\centering
	\begin{tikzpicture}[scale=4.0]
		\begin{scope}[shift={(0,0)}]
			\fill[blue!12] (0,0) rectangle (0.5,0.5);    \draw[thick] (0,0) rectangle (0.5,0.5);
			\fill[green!15] (0.5,0) rectangle (1,0.5);   \draw[thick] (0.5,0) rectangle (1,0.5);
			\fill[blue!12] (0.5,0.5) rectangle (0.75,0.75);  \draw[thick] (0.5,0.5) rectangle (0.75,0.75);
			\fill[blue!12] (0.75,0.5) rectangle (1,0.75);    \draw[thick] (0.75,0.5) rectangle (1,0.75);
			\fill[blue!12] (1,0) rectangle (1.125,0.125);    \draw[thick] (1,0) rectangle (1.125,0.125);
			\fill[blue!12] (1,0.125) rectangle (1.125,0.25); \draw[thick] (1,0.125) rectangle (1.125,0.25);
			\fill[blue!12] (1,0.25) rectangle (1.125,0.375); \draw[thick] (1,0.25) rectangle (1.125,0.375);
			\fill[blue!12] (1,0.375) rectangle (1.125,0.5);  \draw[thick] (1,0.375) rectangle (1.125,0.5);
			\node[font=\tiny,fill=blue!12,inner sep=0.5pt]          (a-v1)  at (0.25,0.25)    {$V_1$};
			\node[font=\tiny,fill=green!15,inner sep=0.5pt]         (a-vl)  at (0.75,0.25)    {$V_\lambda$};
			\node[font=\scalebox{0.75}{\tiny},fill=blue!12,inner sep=0.5pt] (a-v2)  at (0.625,0.625)  {$V_2$};
			\node[font=\scalebox{0.75}{\tiny},fill=blue!12,inner sep=0.5pt] (a-v3)  at (0.875,0.625)  {$V_3$};
			\node[font=\scalebox{0.6}{\tiny},fill=blue!12,inner sep=0.3pt]  (a-v4)  at (1.0625,0.0625) {$V_4$};
			\node[font=\scalebox{0.6}{\tiny},fill=blue!12,inner sep=0.3pt]  (a-v5)  at (1.0625,0.1875) {$V_5$};
			\node[font=\scalebox{0.6}{\tiny},fill=blue!12,inner sep=0.3pt]  (a-v6)  at (1.0625,0.3125) {$V_6$};
			\node[font=\scalebox{0.6}{\tiny},fill=blue!12,inner sep=0.3pt]  (a-v7)  at (1.0625,0.4375) {$V_7$};
			\draw[<->,gray!60,thick] (a-v1) -- (a-vl);
			\draw[<->,gray!60,thick] (a-v2) to[bend right=15] (a-vl);
			\draw[<->,gray!60,thick] (a-v3) to[bend left=15]  (a-vl);
			\draw[<->,gray!60,thick] (0.95,0.0625) -- (1.04,0.0625);
			\draw[<->,gray!60,thick] (0.95,0.1875) -- (1.04,0.1875);
			\draw[<->,gray!60,thick] (0.95,0.3125) -- (1.04,0.3125);
			\draw[<->,gray!60,thick] (0.95,0.4375) -- (1.04,0.4375);
			\node[font=\scriptsize] at (0.625,-0.1) {(a)};
		\end{scope}

		\begin{scope}[shift={(1.45,0)}]
			\fill[blue!12] (0,0) rectangle (0.5,0.5);    \draw[thick] (0,0) rectangle (0.5,0.5);
			\fill[green!15] (0.5,0) rectangle (0.75,0.25);   \draw[thick] (0.5,0) rectangle (0.75,0.25);
			\fill[green!15] (0.75,0) rectangle (1,0.25);     \draw[thick] (0.75,0) rectangle (1,0.25);
			\fill[green!15] (0.5,0.25) rectangle (0.75,0.5); \draw[thick] (0.5,0.25) rectangle (0.75,0.5);
			\fill[green!15] (0.75,0.25) rectangle (1,0.5);   \draw[thick] (0.75,0.25) rectangle (1,0.5);
			\fill[blue!12] (0.5,0.5) rectangle (0.75,0.75);  \draw[thick] (0.5,0.5) rectangle (0.75,0.75);
			\fill[blue!12] (0.75,0.5) rectangle (1,0.75);    \draw[thick] (0.75,0.5) rectangle (1,0.75);
			\fill[blue!12] (1,0) rectangle (1.125,0.125);    \draw[thick] (1,0) rectangle (1.125,0.125);
			\fill[blue!12] (1,0.125) rectangle (1.125,0.25); \draw[thick] (1,0.125) rectangle (1.125,0.25);
			\fill[blue!12] (1,0.25) rectangle (1.125,0.375); \draw[thick] (1,0.25) rectangle (1.125,0.375);
			\fill[blue!12] (1,0.375) rectangle (1.125,0.5);  \draw[thick] (1,0.375) rectangle (1.125,0.5);
			\node[font=\tiny,fill=blue!12,inner sep=0.5pt]          (b-v1)  at (0.25,0.25)    {$V_1$};
			\node[font=\scalebox{0.75}{\tiny},fill=green!15,inner sep=0.5pt] (b-k1) at (0.625,0.125)  {$V_{\kappa_1}$};
			\node[font=\scalebox{0.75}{\tiny},fill=green!15,inner sep=0.5pt] (b-k2) at (0.875,0.125)  {$V_{\kappa_2}$};
			\node[font=\scalebox{0.75}{\tiny},fill=green!15,inner sep=0.5pt] (b-k3) at (0.625,0.375)  {$V_{\kappa_3}$};
			\node[font=\scalebox{0.75}{\tiny},fill=green!15,inner sep=0.5pt] (b-k4) at (0.875,0.375)  {$V_{\kappa_4}$};
			\node[font=\scalebox{0.75}{\tiny},fill=blue!12,inner sep=0.5pt]  (b-v2)  at (0.625,0.625)  {$V_2$};
			\node[font=\scalebox{0.75}{\tiny},fill=blue!12,inner sep=0.5pt]  (b-v3)  at (0.875,0.625)  {$V_3$};
			\node[font=\scalebox{0.6}{\tiny},fill=blue!12,inner sep=0.3pt]   (b-v4)  at (1.0625,0.0625) {$V_4$};
			\node[font=\scalebox{0.6}{\tiny},fill=blue!12,inner sep=0.3pt]   (b-v5)  at (1.0625,0.1875) {$V_5$};
			\node[font=\scalebox{0.6}{\tiny},fill=blue!12,inner sep=0.3pt]   (b-v6)  at (1.0625,0.3125) {$V_6$};
			\node[font=\scalebox{0.6}{\tiny},fill=blue!12,inner sep=0.3pt]   (b-v7)  at (1.0625,0.4375) {$V_7$};
			\draw[<->,gray!60,thick] (b-v1) to[bend right=10] (b-k1);
			\draw[<->,gray!60,thick] (b-v1) to[bend left=10]  (b-k3);
			\draw[<->,gray!60,thick] (b-v2) -- (b-k3);
			\draw[<->,gray!60,thick] (b-v3) -- (b-k4);
			\draw[<->,gray!60,thick] (0.95,0.0625) -- (1.04,0.0625);
			\draw[<->,gray!60,thick] (0.95,0.1875) -- (1.04,0.1875);
			\draw[<->,gray!60,thick] (0.95,0.3125) -- (1.04,0.3125);
			\draw[<->,gray!60,thick] (0.95,0.4375) -- (1.04,0.4375);
			\draw[<->,gray!35,thick,densely dotted] (b-k1) -- (b-k2);
			\draw[<->,gray!35,thick,densely dotted] (b-k3) -- (b-k4);
			\draw[<->,gray!35,thick,densely dotted] (b-k1) -- (b-k3);
			\draw[<->,gray!35,thick,densely dotted] (b-k2) -- (b-k4);
			\node[font=\scriptsize] at (0.625,-0.1) {(b)};
		\end{scope}

		\draw[->,thick] (1.18,0.35) -- (1.41,0.35) node[midway,above,font=\tiny] {refine};
		\draw[->,thick] (1.41,0.15) -- (1.18,0.15) node[midway,below,font=\tiny] {coarsen};
	\end{tikzpicture}
	\caption{Connectivity update for refinement (\cref{alg:connec-refinement}) and coarsening (\cref{alg:connec-coarsening}).
		(a)~Before refinement: $V_\lambda$ (green, $\ell{=}0$) connected to $V_1$ ($\ell{=}0$), $V_2$, $V_3$ ($\ell{=}1$), and $V_4$--$V_7$ ($\ell{=}2$), covering all three neighbour configurations relative to the children.
		(b)~After refining $V_\lambda$ into $V_{\kappa_1}$--$V_{\kappa_4}$ ($\ell{=}1$); sibling connections (dotted) are added separately.
		Reading (b)$\to$(a) illustrates coarsening.}
	\label{fig:connec-update}
\end{figure}

\subsection{Distributed memory parallelisation via MPI}
\label{sec:distr-memory-parall}
With the core data structures established, we now describe how \textsl{MultiWave} scales to distributed memory. Large-scale simulations of hyperbolic balance laws typically require thousands of cores to reach the necessary precision with feasible runtimes, making scalability a critical requirement. Due to the local structure of the DG scheme we target distributed-memory parallelisation via the \emph{Message Passing Interface} (MPI). Domain decomposition is performed using a \emph{space-filling curve} (SFC): each cell is assigned a one-dimensional index derived from its LMI, and the resulting linear order induces a partition into contiguous subdomains distributed across MPI ranks. While this is straightforward on uniform grids, dynamic grid adaptation introduces additional complexity, as the grid changes at every time step and requires periodic repartitioning and ghost cell communication at subdomain boundaries. In the MRA setting, the challenge is compounded by the two-scale transformations, which may require data from cells across MPI boundaries temporarily without permanently transferring ownership~\cite{Brix2009,Brix2011}. We adapt the approach of~\cite{Brix2009,Brix2011} to a Morton-type SFC \cite{Bader2013}.

\subsubsection*{\textbf{SFC-based domain decomposition.}}
\label{sec:SFC-domain-decomp}
Space-filling curves play a central role in domain decomposition: the idea is to map the higher-dimensional adaptive grid to a linear ordering, enabling a contiguous partition of cells across MPI ranks. While the \emph{Hilbert curve} used in~\cite{Brix2009,Brix2011} offers better spatial locality, \textsl{MultiWave} uses a Morton-type SFC, which can be computed more efficiently via bit interleaving of the LMI components~\cite{Bader2013}. An example of the Morton ordering applied to the grid of \cref{fig:lmi-example} is shown in \cref{fig:morton-example}.

\begin{figure}[h]
	\noindent
	\begin{minipage}[t]{0.45\linewidth}
		\centering
		\begin{tikzpicture}[scale=5.5]
			\draw [thick] (0,0) -- (1,0) -- (1,1) -- (0,1) -- (0,0);
			\draw [thick] (0,0.5) -- (1,0.5);
			\draw [thick] (0.5,0) -- (0.5,1);
			\draw [thick] (0.5,0.75) -- (1,0.75);
			\draw [thick] (0.75,0.5) -- (0.75,1);
			\draw [thick] (0.625,0.5) -- (0.625,0.75);
			\draw [thick] (0.5,0.625) -- (0.75,0.625);

			\draw[->,ultra thick,blue] (0.25,0.25) -- (0.75,0.25) -- (0.25,0.75) -- (0.5625,0.5625) -- (0.6875,0.5625) -- (0.5625,0.6875) -- (0.6875,0.6875) -- (0.875,0.625) -- (0.625,0.875) -- (0.875,0.875);

			\node[font=\footnotesize] at (0.07,  0.42)   {$0$};
			\node[font=\footnotesize] at (0.93,  0.42)   {$1$};
			\node[font=\footnotesize] at (0.07,  0.92)   {$2$};
			\node[font=\footnotesize] at (0.5625,0.6)  {$3$};
			\node[font=\footnotesize] at (0.6875,0.6)  {$4$};
			\node[font=\footnotesize] at (0.5625,0.725)  {$5$};
			\node[font=\footnotesize] at (0.6875,0.725)  {$6$};
			\node[font=\footnotesize] at (0.93,  0.72)   {$7$};
			\node[font=\footnotesize] at (0.535, 0.965)  {$8$};
			\node[font=\footnotesize] at (0.93, 0.965)  {$9$};
		\end{tikzpicture}
	\end{minipage}~
	\begin{minipage}[t]{0.45\linewidth}
		\centering
		\vspace{-12.25em}
		\begin{tikzpicture}[scale=0.85,
				lf/.style={circle, draw, fill=white,   minimum size=5.5mm, inner sep=0pt, font=\tiny},
				nd/.style={circle, draw, fill=gray!30, minimum size=5.5mm, inner sep=0pt}]
			\node[anchor=east, font=\scriptsize] at (-0.4, -1.5) {$\ell=0$};
			\node[anchor=east, font=\scriptsize] at (-0.4, -3.0) {$\ell=1$};
			\node[anchor=east, font=\scriptsize] at (-0.4, -4.5) {$\ell=2$};
			\node[lf] (n0)  at (0.0, -1.5) {$0$};
			\node[lf] (n1)  at (1.3, -1.5) {$1$};
			\node[lf] (n2)  at (2.6, -1.5) {$2$};
			\node[nd] (r11) at (5.6, -1.5) {};
			\node[nd] (n22) at (3.55,-3.0) {};
			\node[lf] (n7)  at (5.3, -3.0) {$7$};
			\node[lf] (n8)  at (6.3, -3.0) {$8$};
			\node[lf] (n9)  at (7.3, -3.0) {$9$};
			\draw (r11) -- (n22);
			\draw (r11) -- (n7);
			\draw (r11) -- (n8);
			\draw (r11) -- (n9);
			\node[lf] (n3) at (2.05,-4.5) {$3$};
			\node[lf] (n4) at (3.05,-4.5) {$4$};
			\node[lf] (n5) at (4.05,-4.5) {$5$};
			\node[lf] (n6) at (5.05,-4.5) {$6$};
			\draw (n22) -- (n3);
			\draw (n22) -- (n4);
			\draw (n22) -- (n5);
			\draw (n22) -- (n6);
			\draw[->,ultra thick,blue] (n0) -- (n1) -- (n2) -- (n3) -- (n4) -- (n5) -- (n6) -- (n7) -- (n8) -- (n9);
		\end{tikzpicture}
	\end{minipage}
	\caption{(Left) Morton SFC traversal on the adaptive grid from \cref{fig:lmi-example}; numbers indicate Morton order. (Right) Corresponding quadtree; grey nodes are internal (refined), white leaves are numbered by Morton order.}
	\label{fig:morton-example}
\end{figure}

\subsubsection*{\textbf{Ghost halo communication.}}
\label{sec:ghost-layer-comm}
Flux computations require data from face-adjacent neighbours that may reside on a different MPI rank. For this purpose, each process maintains a single-cell ghost layer, which must be kept up to date.  Two kinds of updates are necessary: a \emph{structural update}, performed after every refinement and coarsening step, which determines the set of ghost cells and their corresponding ranks; and a \emph{data update}, performed e.g., before every surface integral evaluation, which refreshes the solution coefficients stored in the ghost cells. The corresponding rank of any ghost cell is determined by its Morton index via a separator list of the SFC, which encodes the index range assigned to each rank. Both updates are performed via sparse point-to-point communication. To reduce the communication volume, we represent each cell by its 64-bit LMI encoding (cf. \cref{sec:adaptive-solution-representation}) and transmit the raw byte stream of the coefficient matrix.

In the context of SSPRK-DG, communication can be made more efficient by overlapping ghost exchange with local computation. Since no ghost data is needed during the computation of the volume integral, the ghost exchange can be initiated immediately after the previous stage and completed before the surface integral begins. A single time step with asynchronous ghost updates is presented in \cref{alg:dg-rhs-with-async-ghosts}.

\begin{algorithm}[H]
	\caption{\texttt{Modal\_dg::rhs(slab) with asynchronous ghost exchange}}
	\label{alg:dg-rhs-with-async-ghosts}
	\begin{algorithmic}[1]
		\State  $\mathbf{R} \gets $ empty map $\lambda \mapsto \vec{u}_\lambda$
		\State {\color{red}$\texttt{MPI\_handle} \gets \texttt{MPI::start\_ghost\_exchange()}$}
		\For{each cell $V_{\lambda}$}
		\State $\mathbf{R}_\lambda \gets \texttt{Divergence::vol\_coeff}(\vec{u}_\lambda)
			+ \texttt{Source::vol\_coeff}(\vec{u}_\lambda)$
		\EndFor
		\State{\color{red} $\texttt{MPI::complete\_ghost\_exchange(MPI\_handle)}$}
		\For{each interior edge $e$ between cells $V_{\lambda_1}$, $V_{\lambda_2}$}
		\State $(\vec{w}_{L}, \vec{w}_{R}) \gets
			\texttt{Div::surf\_coeff}(e, \vec{u}_{\lambda_1}, \vec{u}_{\lambda_2})
			+ \texttt{Source::surf\_coeff}(e, \vec{u}_{\lambda_1}, \vec{u}_{\lambda_2})$
		\State $\mathbf{R}_{\lambda_1} \mathrel{-}= \vec{w}_{L} $
		\State $\mathbf{R}_{\lambda_2} \mathrel{+}= \vec{w}_{R}$
		\EndFor
		\State $\texttt{bc.apply}{(t,\,\vec{u},\,\mathbf{R})}$
		\State \Return $\mathbf{R}$
	\end{algorithmic}
\end{algorithm}

The key idea is that \texttt{MPI::start\_ghost\_exchange} posts non-blocking \texttt{MPI\_Isend}/\texttt{Irecv} calls and returns immediately, so that the volume integrals, which depend only on local data, proceed in parallel with the network transfer. Ghost values are unpacked only when \texttt{MPI::complete\_ghost\_exchange} is called, just before they are required by \textsc{surface\_integrals}, effectively hiding the communication latency behind volume integral computation. Furthermore, since the grid does not change between stages of a single time step, the ghost cell neighbourhood can be cached, reducing repeated MPI setup overhead.

After each refinement or coarsening step, the ghost layer must be updated: former ghost cells may no longer exist, new boundary cells may have appeared, and subdomain connectivity must be revised. Each rank exchanges the LMIs of its boundary cells with neighbouring ranks and receives the LMIs of new potential ghost cells. Each received cell is connected to the local grid by searching, in each outward direction, for a direct neighbour, a coarse ancestor, or a finer descendant in the received set, following the same logic as \cref{alg:connec-coarsening} and \cref{alg:connec-refinement}. The resulting complexity is $\mathcal{O}(d\cdot L\cdot|\mathcal{G}|)$, where $\mathcal{G}$ is the ghost cell set. The full procedure is given as \cref{alg:ghost-neigh-update} in \cref{sec:ghost-neigh-update}. \Cref{fig:ghost-update} illustrates the three connection cases: for each new ghost, the algorithm computes the hypothetical same-level neighbour via \textsc{next\_index} and projects candidates to a common level using \textsc{get\_parent}. In panel~(b), $V_{\mu_2}$ matches $V_{\kappa_2}$ at the same level; $V_{\mu_3}$ is projected to the level of the coarser $V_{\kappa_1}$, which matches. In panel~(d), the finer local cells $V_{\mu_3}$ and $V_{\mu_5}$ are projected to the level of $V_{\kappa_1}$'s same-level neighbour index, identifying them as descendants.

\subsubsection*{\textbf{MRA two-scale communication.}}
\label{sec:mra-comm}
Computing the detail coefficients of a cell $V_\lambda$ in the two-scale transform requires the DG coefficients of its $2^d$ children. Under a standard SFC-based load-balancing strategy, siblings may be distributed across different MPI ranks, necessitating data communication at every level of the transform hierarchy. Since the two-scale computation requires data not only from active grid cells but also from virtual cells in the MRA hierarchy, the communication overhead and associated bookkeeping become substantial.

To avoid this, we require that the complete sub-tree of each level-0 cell resides on a single MPI rank. This is achieved by assigning each cell the Morton index of its level-0 ancestor as its partition key, rather than its own Morton index.
\revAK{While this does not guarantee perfect load balance, the number of cells per rank deviates from the mean $N/Q$, where $N$ denotes the total number of cells and $Q$ the number of MPI ranks, by at most $\max_{\lambda\in\mathcal I_0} |\mathrm{subtree}_\lambda|$. Hence the relative imbalance is bounded by $Q \max_{\lambda\in\mathcal I_0} |\mathrm{subtree}_\lambda| / N$, which becomes negligible when $|\mathcal I_0| \gg Q$. Using this approach, MPI communication is performed only once during grid refinement and grid coarsening, respectively.}


\begin{figure}[h]
	\centering
	\begin{tikzpicture}[scale=3.0]

		\begin{scope}[shift={(0,0)}]
			\fill[blue!12] (0,0)     rectangle (0.5,0.5);   \draw[thick] (0,0)     rectangle (0.5,0.5);
			\fill[blue!12] (0,0.5)   rectangle (0.25,0.75);  \draw[thick] (0,0.5)   rectangle (0.25,0.75);
			\fill[blue!12] (0.25,0.5) rectangle (0.5,0.75);  \draw[thick] (0.25,0.5) rectangle (0.5,0.75);
			\fill[blue!12] (0,0.75)  rectangle (0.25,1.0);   \draw[thick] (0,0.75)  rectangle (0.25,1.0);
			\fill[blue!12] (0.25,0.75) rectangle (0.5,1.0);  \draw[thick] (0.25,0.75) rectangle (0.5,1.0);
			\fill[orange!10] (0.5,0)   rectangle (1.0,0.5);  \draw[orange!70!black,thick,dashed] (0.5,0)   rectangle (1.0,0.5);
			\fill[orange!10] (0.5,0.5) rectangle (1.0,1.0);  \draw[orange!70!black,thick,dashed] (0.5,0.5) rectangle (1.0,1.0);
			\draw[very thick,red!75!black] (0.5,0) -- (0.5,1.0);
			\node[font=\tiny,fill=blue!12,inner sep=0.5pt] (a-k1) at (0.25,0.25)   {$V_{\kappa_1}$};
			\node[font=\tiny,fill=blue!12,inner sep=0.5pt]        at (0.125,0.625) {$V_{\kappa_3}$};
			\node[font=\tiny,fill=blue!12,inner sep=0.5pt] (a-k4) at (0.375,0.625) {$V_{\kappa_4}$};
			\node[font=\tiny,fill=blue!12,inner sep=0.5pt]        at (0.125,0.875) {$V_{\kappa_5}$};
			\node[font=\tiny,fill=blue!12,inner sep=0.5pt] (a-k6) at (0.375,0.875) {$V_{\kappa_6}$};
			\node[font=\tiny,fill=orange!10,inner sep=0.5pt] (a-m1) at (0.75,0.25)   {$V_{\mu_1}$};
			\node[font=\tiny,fill=orange!10,inner sep=0.5pt] (a-m2) at (0.75,0.75)   {$V_{\mu_2}$};
			\draw[<->,gray!60,thick] (a-k1) -- (a-m1);
			\draw[<->,gray!60,thick] (a-k4) to[bend right=15] (a-m2);
			\draw[<->,gray!60,thick] (a-k6) to[bend left=15]  (a-m2);
			\node[font=\scriptsize] at (0.5,-0.13) {(a)};
		\end{scope}

		\begin{scope}[shift={(1.3,0)}]
			\fill[blue!12] (0,0)   rectangle (0.5,0.5);  \draw[thick] (0,0)   rectangle (0.5,0.5);
			\fill[blue!12] (0,0.5) rectangle (0.5,1.0);  \draw[thick] (0,0.5) rectangle (0.5,1.0);
			\fill[gray!12] (0.75,0)    rectangle (1.0,0.25);  \draw[gray!45,semithick] (0.75,0)    rectangle (1.0,0.25);
			\fill[gray!12] (0.75,0.25) rectangle (1.0,0.5);   \draw[gray!45,semithick] (0.75,0.25) rectangle (1.0,0.5);
			\fill[orange!10] (0.5,0)    rectangle (0.75,0.25); \draw[orange!70!black,thick,dashed] (0.5,0)    rectangle (0.75,0.25);
			\fill[orange!10] (0.5,0.25) rectangle (0.75,0.5);  \draw[orange!70!black,thick,dashed] (0.5,0.25) rectangle (0.75,0.5);
			\fill[orange!10] (0.5,0.5)  rectangle (1.0,1.0);   \draw[orange!70!black,thick,dashed] (0.5,0.5)  rectangle (1.0,1.0);
			\draw[very thick,red!75!black] (0.5,0) -- (0.5,1.0);
			\node[font=\tiny,fill=blue!12,inner sep=0.5pt] (b-k1) at (0.25,0.25)   {$V_{\kappa_1}$};
			\node[font=\tiny,fill=blue!12,inner sep=0.5pt] (b-k2) at (0.25,0.75)   {$V_{\kappa_2}$};
			\node[font=\tiny,fill=orange!10,inner sep=0.5pt] (b-m3) at (0.625,0.125) {$V_{\mu_3}$};
			\node[font=\tiny,fill=orange!10,inner sep=0.5pt] (b-m5) at (0.625,0.375) {$V_{\mu_5}$};
			\node[font=\tiny,fill=orange!10,inner sep=0.5pt] (b-m2) at (0.75,0.75)   {$V_{\mu_2}$};
			\draw[<->,gray!60,thick] (b-k1) to[bend right=15] (b-m3);
			\draw[<->,gray!60,thick] (b-k1) to[bend left=15]  (b-m5);
			\draw[<->,gray!60,thick] (b-k2) -- (b-m2);
			\node[font=\scriptsize] at (0.5,-0.13) {(b)};
		\end{scope}

		\begin{scope}[shift={(0,-1.5)}]
			\fill[gray!12] (0,0.5)  rectangle (0.25,0.75);  \draw[gray!45,semithick] (0,0.5)  rectangle (0.25,0.75);
			\fill[gray!12] (0,0.75) rectangle (0.25,1.0);   \draw[gray!45,semithick] (0,0.75) rectangle (0.25,1.0);
			\fill[blue!8] (0,0)       rectangle (0.5,0.5);   \draw[blue!70!black,thick,dashed] (0,0)       rectangle (0.5,0.5);
			\fill[blue!8] (0.25,0.5)  rectangle (0.5,0.75);  \draw[blue!70!black,thick,dashed] (0.25,0.5)  rectangle (0.5,0.75);
			\fill[blue!8] (0.25,0.75) rectangle (0.5,1.0);   \draw[blue!70!black,thick,dashed] (0.25,0.75) rectangle (0.5,1.0);
			\fill[orange!12] (0.5,0)   rectangle (1.0,0.5);  \draw[thick] (0.5,0)   rectangle (1.0,0.5);
			\fill[orange!12] (0.5,0.5) rectangle (1.0,1.0);  \draw[thick] (0.5,0.5) rectangle (1.0,1.0);
			\draw[very thick,red!75!black] (0.5,0) -- (0.5,1.0);
			\node[font=\tiny,fill=blue!8,inner sep=0.5pt] (c-k1) at (0.25,0.25)   {$V_{\kappa_1}$};
			\node[font=\tiny,fill=blue!8,inner sep=0.5pt] (c-k4) at (0.375,0.625) {$V_{\kappa_4}$};
			\node[font=\tiny,fill=blue!8,inner sep=0.5pt] (c-k6) at (0.375,0.875) {$V_{\kappa_6}$};
			\node[font=\tiny,fill=orange!12,inner sep=0.5pt] (c-m1) at (0.75,0.25)   {$V_{\mu_1}$};
			\node[font=\tiny,fill=orange!12,inner sep=0.5pt] (c-m2) at (0.75,0.75)   {$V_{\mu_2}$};
			\draw[<->,gray!60,thick] (c-k1) -- (c-m1);
			\draw[<->,gray!60,thick] (c-k4) to[bend right=15] (c-m2);
			\draw[<->,gray!60,thick] (c-k6) to[bend left=15]  (c-m2);
			\node[font=\scriptsize] at (0.5,-0.13) {(c)};
		\end{scope}

		\begin{scope}[shift={(1.3,-1.5)}]
			\fill[blue!8] (0,0)   rectangle (0.5,0.5);  \draw[blue!70!black,thick,dashed] (0,0)   rectangle (0.5,0.5);
			\fill[blue!8] (0,0.5) rectangle (0.5,1.0);  \draw[blue!70!black,thick,dashed] (0,0.5) rectangle (0.5,1.0);
			\fill[orange!12] (0.5,0)    rectangle (0.75,0.25); \draw[thick] (0.5,0)    rectangle (0.75,0.25);
			\fill[orange!12] (0.75,0)   rectangle (1.0,0.25);  \draw[thick] (0.75,0)   rectangle (1.0,0.25);
			\fill[orange!12] (0.5,0.25) rectangle (0.75,0.5);  \draw[thick] (0.5,0.25) rectangle (0.75,0.5);
			\fill[orange!12] (0.75,0.25) rectangle (1.0,0.5);  \draw[thick] (0.75,0.25) rectangle (1.0,0.5);
			\fill[orange!12] (0.5,0.5)  rectangle (1.0,1.0);   \draw[thick] (0.5,0.5)  rectangle (1.0,1.0);
			\draw[very thick,red!75!black] (0.5,0) -- (0.5,1.0);
			\node[font=\tiny,fill=blue!8,inner sep=0.5pt] (d-k1) at (0.25,0.25)   {$V_{\kappa_1}$};
			\node[font=\tiny,fill=blue!8,inner sep=0.5pt] (d-k2) at (0.25,0.75)   {$V_{\kappa_2}$};
			\node[font=\tiny,fill=orange!12,inner sep=0.5pt] (d-m3) at (0.625,0.125) {$V_{\mu_3}$};
			\node[font=\tiny,fill=orange!12,inner sep=0.5pt]        at (0.875,0.125) {$V_{\mu_4}$};
			\node[font=\tiny,fill=orange!12,inner sep=0.5pt] (d-m5) at (0.625,0.375) {$V_{\mu_5}$};
			\node[font=\tiny,fill=orange!12,inner sep=0.5pt]        at (0.875,0.375) {$V_{\mu_6}$};
			\node[font=\tiny,fill=orange!12,inner sep=0.5pt] (d-m2) at (0.75,0.75)   {$V_{\mu_2}$};
			\draw[<->,gray!60,thick] (d-k1) to[bend right=15] (d-m3);
			\draw[<->,gray!60,thick] (d-k1) to[bend left=15]  (d-m5);
			\draw[<->,gray!60,thick] (d-k2) -- (d-m2);
			\node[font=\scriptsize] at (0.5,-0.13) {(d)};
		\end{scope}

		\node[font=\small] at (0.5,  1.12) {before};
		\node[font=\small] at (1.8,  1.12) {after};
		\node[font=\small,rotate=90] at (-0.14,  0.5)  {rank~0};
		\node[font=\small,rotate=90] at (-0.14, -1.0)  {rank~1};
		\draw[->,thick] (1.04, 0.5)  -- (1.26, 0.5)  node[midway,above,font=\tiny] {update};
		\draw[->,thick] (1.04, -1.0) -- (1.26, -1.0) node[midway,above,font=\tiny] {update};

	\end{tikzpicture}
	\caption{Ghost neighbour update (\cref{alg:ghost-neigh-update}) from both ranks' perspectives. Blue solid: rank~0 local; orange solid: rank~1 local; dashed: ghost copies;
		faded: non-ghost cells of the remote rank. Before: rank~0's top cell $V_{\kappa_2}$ is refined (children $V_{\kappa_3}$--$V_{\kappa_6}$), rank~1 is uniform. After:
		$V_{\kappa_2}$'s children are coarsened and $V_{\mu_1}$ is refined. (a,c)~Valid connectivity before the grid changes. (b)~Rank~0's ghost layer updated: $V_{\mu_1}$
		replaced by adjacent children $V_{\mu_3}$, $V_{\mu_5}$; $V_{\kappa_2}$ now connects to $V_{\mu_2}$. (d)~Rank~1's connectivity updated: ghost $V_{\kappa_2}$ replaces
		$V_{\kappa_4}$, $V_{\kappa_6}$, and new local cells $V_{\mu_3}$, $V_{\mu_5}$ connect to $V_{\kappa_1}$.}
	\label{fig:ghost-update}
\end{figure}

\subsubsection*{\textbf{Load rebalancing.}}
\label{sec:load-rebalancing}
To monitor load balance, \textsl{MultiWave} checks for imbalance after each refinement or coarsening step and triggers rebalancing when the relative imbalance \revakk{$(\max_q |\mathcal{U}_q| - \min_q |\mathcal{U}_q|)/N$, where $\mathcal U_q$ denotes the cells owned by rank $q$,} exceeds a user-defined threshold $\tau_\mathrm{lb}\in[0,1]$. \revakk{The optimal choice of $\tau_\mathrm{lb}$ balances the cost of redistributing cells against the idle time caused by unequal rank loads and is therefore problem-dependent and determined empirically.} The rebalancing procedure, shown in \cref{alg:load-rebalancing}, proceeds in three phases: the SFC separators are updated to reflect the new target partition, cells whose assigned rank has changed are packed together with their solution data and neighbour lists and migrated via sparse point-to-point communication, and finally a ghost exchange refreshes the ghost layer on all ranks. Since rebalancing is triggered only when the threshold is exceeded, its overhead is amortised over a sufficiently large number of time steps.

\begin{algorithm}[ht]
	\caption{Load rebalancing}
	\label{alg:load-rebalancing}
	\begin{algorithmic}[1]
		\Require{distributed solution $\mathcal{U} = \bigcup_{q=0}^{Q-1} \mathcal{U}_q$ over $Q$ MPI ranks,
			adjacency list $\mathcal{A}$, load balance threshold $\tau_{\mathrm{lb}} > 0$}
		\If{$(\max_q |\mathcal{U}_q| - \min_q |\mathcal{U}_q|)/\sum_q |\mathcal{U}_q| > \tau_{\mathrm{lb}}$}
		\State $\textsc{update\_separators}(\mathcal{U})$
		\State $\mathcal{M} \gets \{V_\lambda \in \mathcal{U}_q :
			\textsc{get\_process}(V_\lambda) \neq q\}$
		\Comment{cells to migrate}
		\For{each $V_\lambda \in \mathcal{M}$}
		\State pack $(V_\lambda,\; \mathcal{U}_q[V_\lambda],\; \mathcal{N}_\lambda)$ into send buffer
		\Comment{cell data + neighbour list}
		\State erase $V_\lambda$ from $\mathcal{U}_q$
		\EndFor
		\State exchange send buffers via sparse P2P
		\For{each received $(V_\lambda,\; \mathbf{u},\; \mathcal{N}_\lambda)$}
		\State insert $V_\lambda$ into $\mathcal{U}_q$ with data $\mathbf{u}$
		\State update connectivity with $\mathcal{N}_\lambda$
		\EndFor
		\State $\textsc{ghost\_exchange}(\mathcal{U})$
		\EndIf
	\end{algorithmic}
\end{algorithm}

The following section presents numerical results that validate both the parallel performance and the modular design of \textsl{MultiWave}.

\section{Numerical results}
\label{sec:appl-fram}
This section presents  numerical studies that demonstrate the capabilities of \textsl{MultiWave}. \revakk{First, single-core studies on smooth and discontinuous solutions inspect the computational costs   and validate the correctness of the adaptive numerical method. Second, grid-scaling studies on the Taylor--Green vortex confirm the parallel efficiency of the framework on both uniform and adaptive grids. Finally, two applications are presented:  random conservation laws illustrate the weighted MRA extension and a fluid-structure coupling example showcases the modularity of the framework.} 

For all numerical examples we use a local Lax--Friedrichs flux~\cite{toroRiemannSolversNumerical2009}:
\begin{align*}
\widehat{\vec f}\left(\vec u^+_h, \vec u^-_h, \vec n_\lambda\right)
= \frac{1}{2}\left(\vec f\left(\vec u^+_h\right) + \vec f\left(\vec u^-_h\right)\right) \cdot \vec n_\lambda
- \frac{\lambda_\text{max}\left(\vec u^+_h, \vec u^-_h, \vec n_\lambda\right)}{2}
  \left(\vec u^+_h - \vec u^-_h\right),
\end{align*}
where
\begin{align*}
\lambda_\text{max}\left(\vec u^+_h, \vec u^-_h, \vec n_\lambda\right)
= \max_{\vec v \in \{\vec u^+_h,\, \vec u^-_h\}}
  \max_{i=1,\dots,m} \left| \lambda_i\!\left(\mathbf{A}(\vec v, \vec n_\lambda)\right) \right|
\end{align*}
denotes the largest absolute eigenvalue of $\mathbf{A}$, see~\eqref{eq:flux-jacobian}.
For time stepping in the RKDG scheme we apply an order-dependent CFL condition
\begin{align*}
\Delta t \leq \text{CFL}\cdot\frac{C}{2p-1}
\left(\sum_{j=1}^{d}\frac{\lambda_{\text{max},j}(\vec u_h)}{\Delta x_j}\right)^{-1},
\end{align*}
where $\Delta x_j$ is the smallest diameter in the $j$-th coordinate direction over all cells $V_\lambda$ in the grid $\mathcal G_{L,\varepsilon}$, $\lambda_{\text{max},j}(\vec u_h)$ is the magnitude of the maximum characteristic speed in the $j$-th direction, $C>0$ is the SSP coefficient in the time integration and $\text{CFL}>0$, see~\cite{cockburnRungeKuttaDiscontinuous2001}.

Beyond these, prototype versions of \textsl{MultiWave} have been applied in a variety of further works, including uncertainty quantification~\cite{gersterHyperbolicStochasticGalerkin2019,herty2022higher,herty2024novel,herty2024multiresolution,Kolb2023}, model adaptation~\cite{giesselmann2024model,Joshi2023}, boundary control~\cite{Herty2025}, data compression~\cite{Khosrawi2026}, \revAK{coupled Riemann problems~\cite{Du2026a}}, and a posteriori error estimates~\cite{giesselmann2025posteriori}. The framework has further been benchmarked on radially symmetric ultra-relativistic Euler equations~\cite{kunik2024radially}.

\subsection{Validation and single-core performance}
\label{sec:single-core}
\revakk{We first verify the performance of \textsl{MultiWave} on a single core. Here we consider several performance criteria which not only show preservation of the accuracy of the reference scheme, cf.~\eqref{eq:error-sparse-approx}, but also quantify adaptation costs. To this end, we consider the compressible Euler equations for a perfect gas with no source term and conserved variables $\vec u = (\rho, \rho \vec v, \rho E)^\top$,
\begin{equation}
	\label{eq:compressible-euler}
	\frac{\partial}{\partial t}\begin{pmatrix}\rho        \\
		\rho\vec{v} \\
		\rho E
	\end{pmatrix}
	+ \nabla\cdot \begin{pmatrix}
		\rho\vec{v}                               \\
		\rho\vec{v}\otimes\vec{v} + P\,\mathbf{I} \\
		\vec{v}(\rho E + P)
	\end{pmatrix} = \vec{0},
\end{equation}
where $\rho > 0$ is the density, $\vec{v} \in \mathbb{R}^d$ the velocity, $E$ the specific total energy, $P$ the pressure and $\mathbf{I}\in\R^{d\times d}$ the identity matrix. The system is closed by the equation of state (EoS) for a perfect gas,
\begin{equation}\label{eq:perfect-gas-eos}
	P = (\kappa - 1)\!\left(\rho E - \tfrac{1}{2}\rho|\vec{v}|^2\right),
\end{equation}
with adiabatic index $\kappa > 1$. We consider four test cases for the compressible Euler equations~\eqref{eq:compressible-euler}, two with smooth solutions and two containing discontinuities. Each pair consists of a one-dimensional and a two-dimensional case. The smooth test cases verify the order of the underlying RKDG scheme with grid adaptation. The discontinuous test cases  show the compression capabilities of the multiresolution-based grid adaptation.}

\subsubsection*{\textbf{1. One-dimensional smooth case}:}
\revakk{A Ricker wavelet  in the density is advected at a constant velocity and constant pressure on the domain $\Omega = [-5,10]$ with periodic boundary conditions and final time $T = 0.2$,
\begin{equation}
	\label{eq:hat-initial}
	\rho_0(x) = 1 + \psi(x), \quad v_0(x) \equiv 30, \quad P_0(x) \equiv 1, \qquad
	\psi(x) = \frac{2}{\sqrt{3}\,\pi^{1/4}}\left(1 - x^2\right)e^{-x^2/2}.
\end{equation}
Since the pressure and the velocity are constant, the density profile is transported and the exact solution is the periodic translation $\vec u = \vec u_0(x-v_0t)$ for all $t \geq 0$. Since the solution is smooth, we do not perform limiting and we choose $\gamma = p$ in the local threshold value~\eqref{eq:local-threshold-value}.}

\subsubsection*{\textbf{2. One-dimensional discontinuous case}:}
\revakk{Here, we consider \emph{Sod's shock tube} problem~\cite{sod1978survey} on $\Omega = [0,1]$ and final time $T=0.2$,
\begin{equation}
	\label{eq:sod-initial}
	(\rho_0, v_0, P_0)(x) = \begin{cases}
		(1,\ 0,\ 1),       & x < \tfrac12, \\
		(0.125,\ 0,\ 0.1), & x > \tfrac12.
	\end{cases}
\end{equation}
The exact solution is obtained from the exact Riemann solver for the Euler equations. It consists of  a rarefaction, a contact and a shock wave. Since the solution is discontinuous, we perform limiting on the finest resolution level. Due to the presence of discontinuities, we expect first-order convergence and choose $\gamma = 1$ in the local threshold value~\eqref{eq:local-threshold-value}.}

\subsubsection*{\textbf{3. Two-dimensional smooth case}:}
\revakk{We consider the \emph{isentropic Euler vortex}~\cite{spiegel2015survey} on $\Omega = [-10,10]^2$ with periodic boundary conditions and final time $T = 20/M_\infty$. With angular strength $\omega(\vec x) = \beta e^{-(x_0^2+x_1^2)/2}$ and temperature deviation $\delta T(\vec x) = -\tfrac12(\kappa-1)\omega(\vec x)^2$, the initial data read
\begin{equation}
	\label{eq:vortex-initial}
	\vec v_0(\vec x) = \begin{pmatrix} M_\infty - x_1\,\omega(\vec x) \\ x_0\,\omega(\vec x)\end{pmatrix}, \qquad
	\rho_0(\vec x) = \left(1 + \delta T(\vec x)\right)^{\frac{1}{\kappa-1}}, \qquad
	P_0(\vec x) = \frac{1}{\kappa}\left(1 + \delta T(\vec x)\right)^{\frac{\kappa}{\kappa-1}},
\end{equation}
with $M_\infty = \sqrt{2/\kappa}$ and $\beta = \tfrac{5}{2\pi}\sqrt{\mathrm{e}}$. Since the entropy is constant, the vortex is transported and the exact solution is the periodic translation $\vec u = \vec u_0\big(\vec x - (M_\infty t, 0)^\top\big)$ for all $t\geq0$, where the final time $T$ corresponds to one traversal of the domain. Since the solution is smooth, we do not perform limiting and we choose $\gamma = p$ in the local threshold value~\eqref{eq:local-threshold-value}.}

\subsubsection*{\textbf{4. Two-dimensional discontinuous case}:}
\revakk{We set the initial data of the one-dimensional Sod shock tube problem~\eqref{eq:sod-initial} on the circle with radius $r = \tfrac12$ on $\Omega = [-1,1]^2$ with final time $T = 0.2$, i.e.,
\begin{equation}
	\label{eq:radial-sod-initial}
	(\rho_0, \vec v_0, P_0)(\vec x) = \begin{cases}
		(1,\ \vec 0,\ 1),       & \|\vec x\|_2 < r, \\
		(0.125,\ \vec 0,\ 0.1), & \|\vec x\|_2 > r.
	\end{cases}
\end{equation}
A reference solution is obtained from a high-resolution one-dimensional computation of the radially symmetric system
\begin{align*}
	\frac{\partial}{\partial t} \vec U + \partial_r \vec F (\vec U) = -\frac1r\left(\rho v, \rho v^2, v(\rho E + P)\right)^\top
\end{align*}
where the singularity at $r\to 0$ is removable and, thus, the domain is clamped as $ r \geq 0.1 \Delta x$. For this test case we require a positivity-preserving limiter in addition to the Cockburn--Shu limiter to prevent non-physical states. Again, we expect first-order convergence and set $\gamma = 1$ in the local threshold value~\eqref{eq:local-threshold-value}.}

\revakk{For all four test cases, we set the order of the DG scheme to $p = 4$, the CFL number to $0.5$ and use a local Lax--Friedrichs flux together with a fourth-order, ten-stage SSPRK scheme. The initial grid at the coarsest refinement level consists of $15$ cells for the one-dimensional smooth case, $10$ cells for the one-dimensional discontinuous case and $10\times10$ cells for each of the two-dimensional cases. The one-dimensional cases run for a maximum refinement level $L = 3,\ldots,7$, the isentropic vortex for $L = 3,\ldots,6$ and the radial Sod problem for $L = 4,\ldots,6$, each for the three threshold constants $C_\text{thr}\in\{10^{-3},10^{-1},1\}$. All simulations are performed on one pinned core of an AMD EPYC 9555 processor (Zen~5, $4.4$~GHz sustained single-core clock, 1~MiB private L2 cache per core, 32~MiB shared L3 cache per core complex, 32 floating-point operations per core per cycle). Instruction counts and floating-point throughput are obtained from the hardware performance counters.}

\begin{figure}[htbp]
	\centering
	\begin{subfigure}{0.5\textwidth}
		\includegraphics[width=\textwidth]{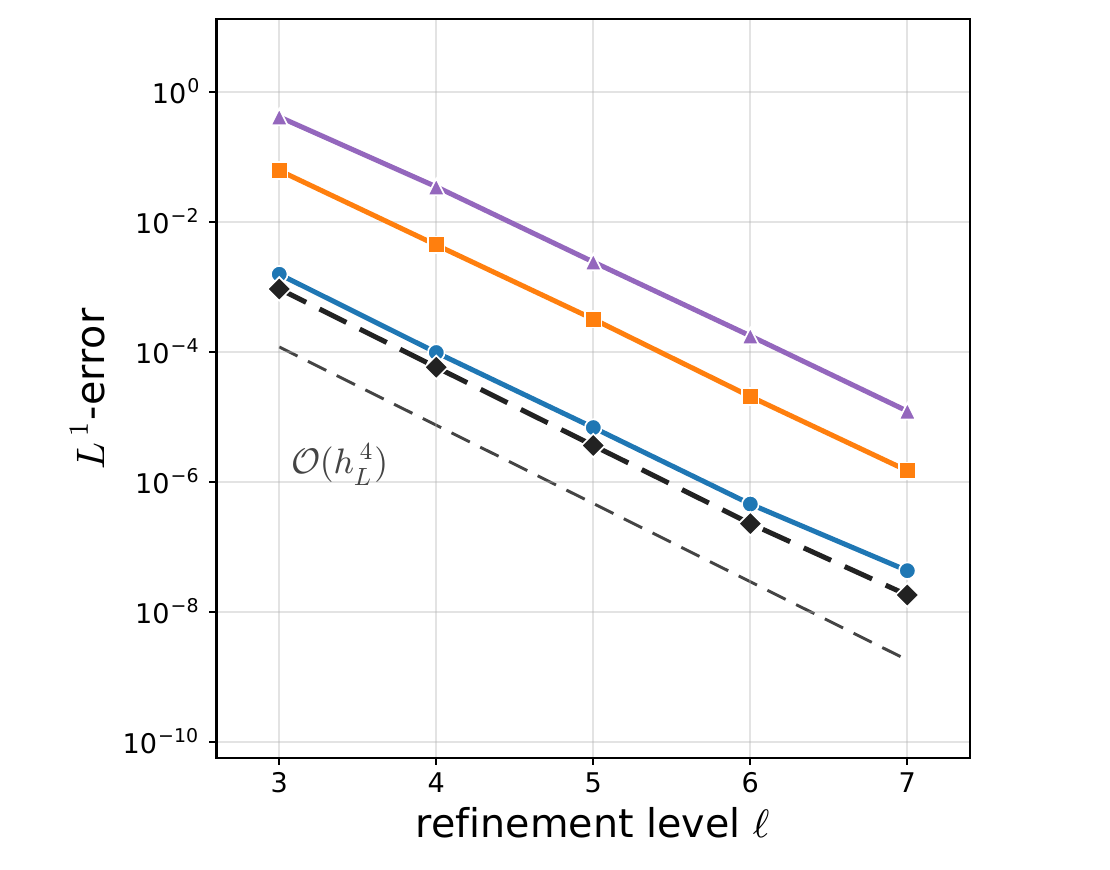}
		\caption{accuracy against refinement level}\label{fig:sc-hat-conv}
	\end{subfigure}\hfill
	\begin{subfigure}{0.5\textwidth}
		\includegraphics[width=\textwidth]{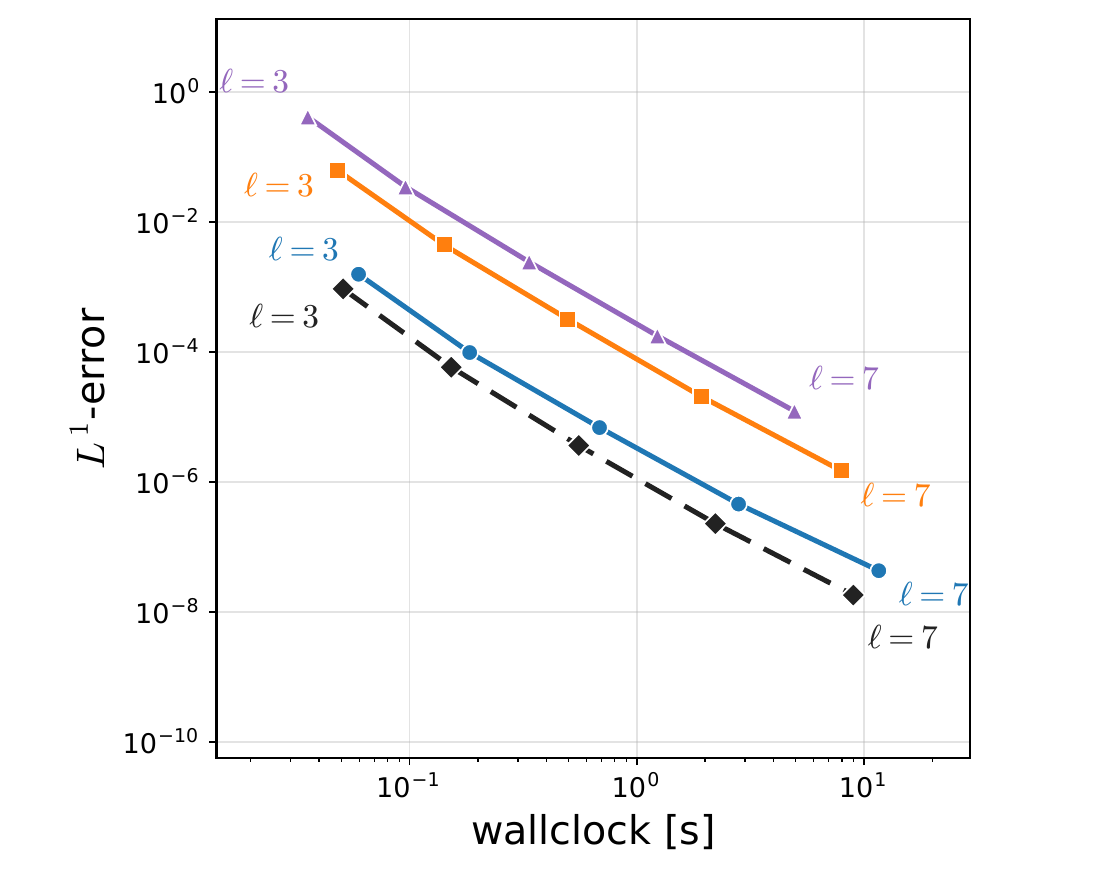}
		\caption{accuracy against runtime}\label{fig:sc-hat-cost}
	\end{subfigure}

	\vspace{2ex}

	\begin{subfigure}{0.5\textwidth}
		\includegraphics[width=\textwidth]{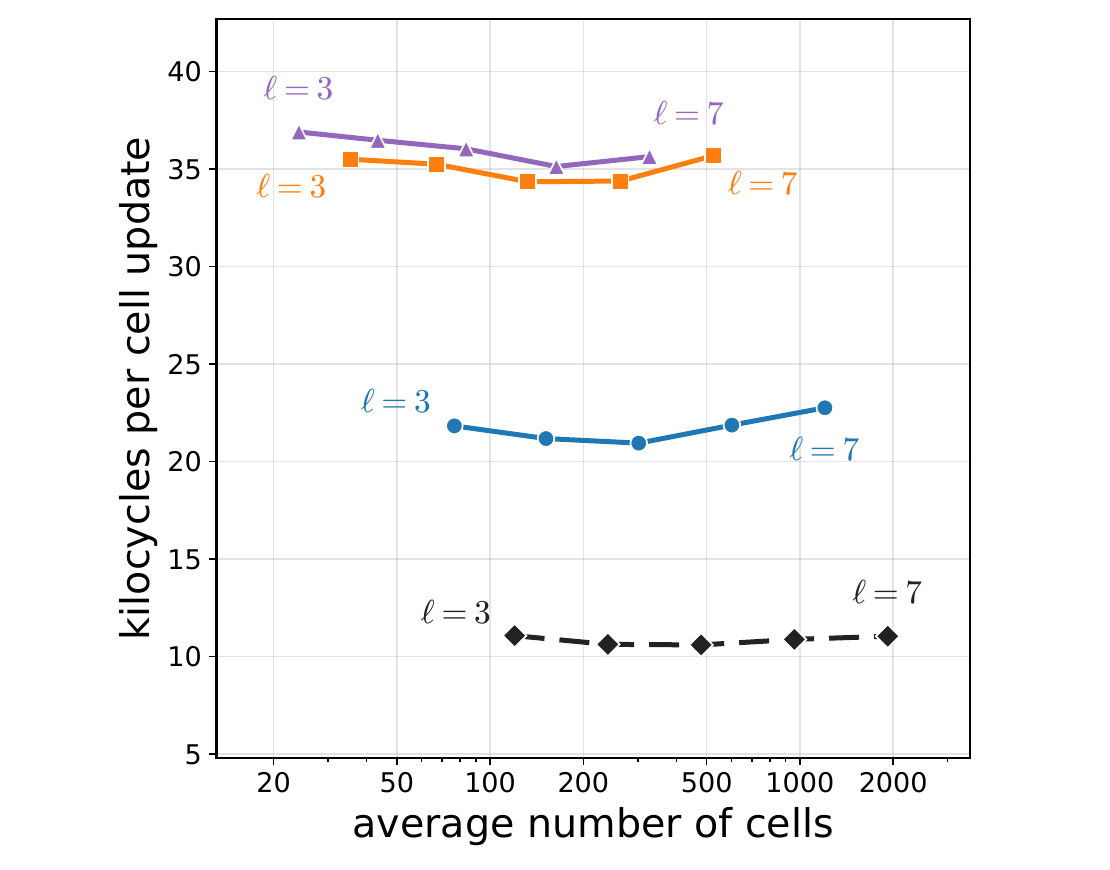}
		\caption{cost of one cell update}\label{fig:sc-hat-cell}
	\end{subfigure}\hfill
	\begin{subfigure}{0.5\textwidth}
		\includegraphics[width=\textwidth]{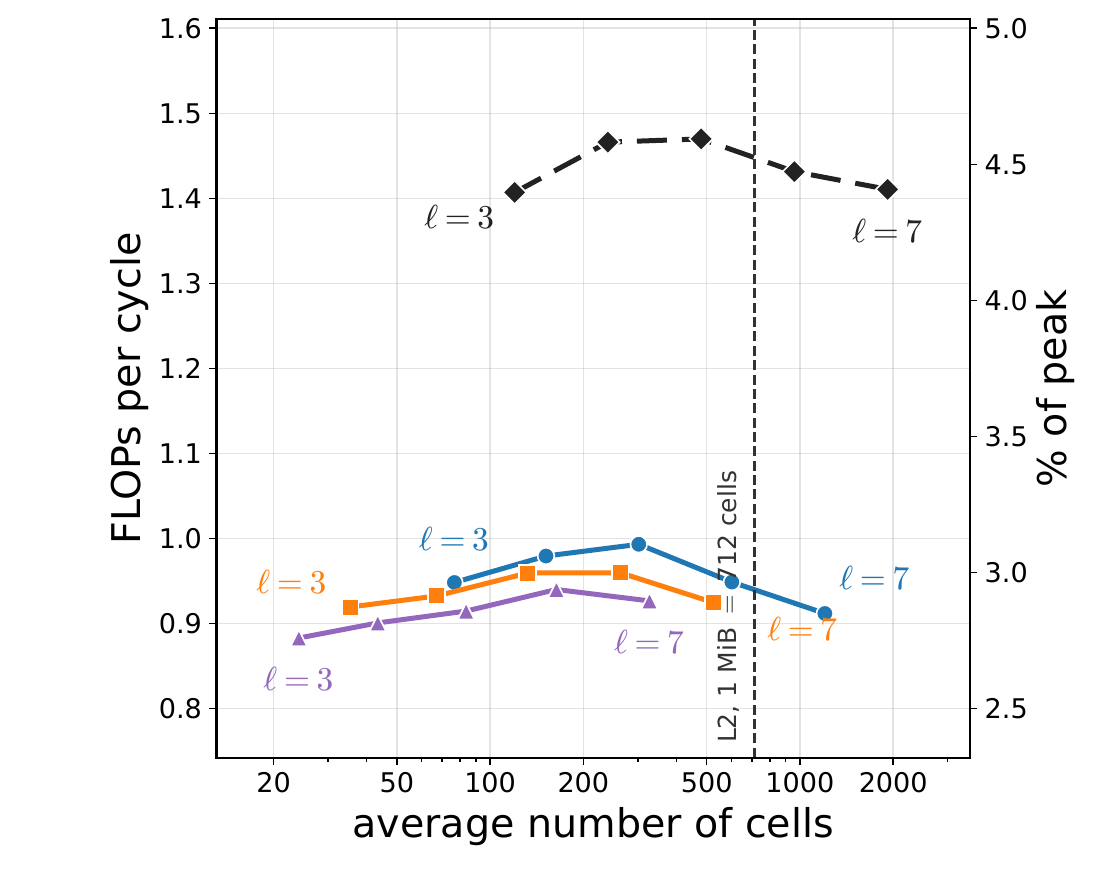}
		\caption{machine utilisation}\label{fig:sc-hat-util}
	\end{subfigure}

	\includegraphics[width=\textwidth]{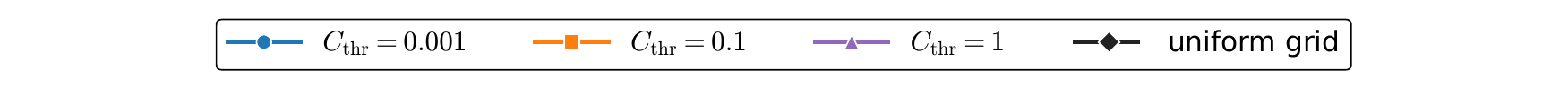}
	\caption{Advected Ricker wavelet~\eqref{eq:hat-initial} for the three threshold constants and the uniform reference grid.}
	\label{fig:sc-hat}
\end{figure}

\begin{figure}[htbp]
	\centering
	\begin{subfigure}{0.5\textwidth}
		\includegraphics[width=\textwidth]{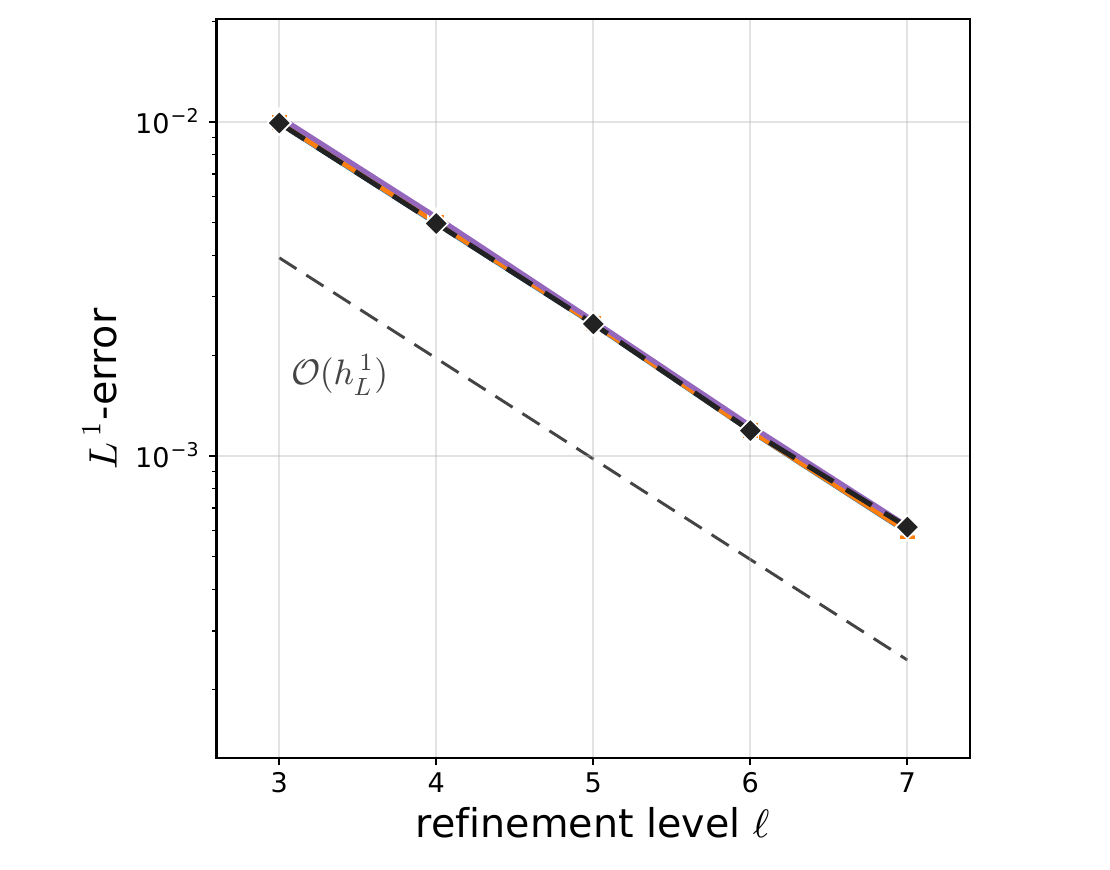}
		\caption{accuracy against refinement level}\label{fig:sc-sod-conv}
	\end{subfigure}\hfill
	\begin{subfigure}{0.5\textwidth}
		\includegraphics[width=\textwidth]{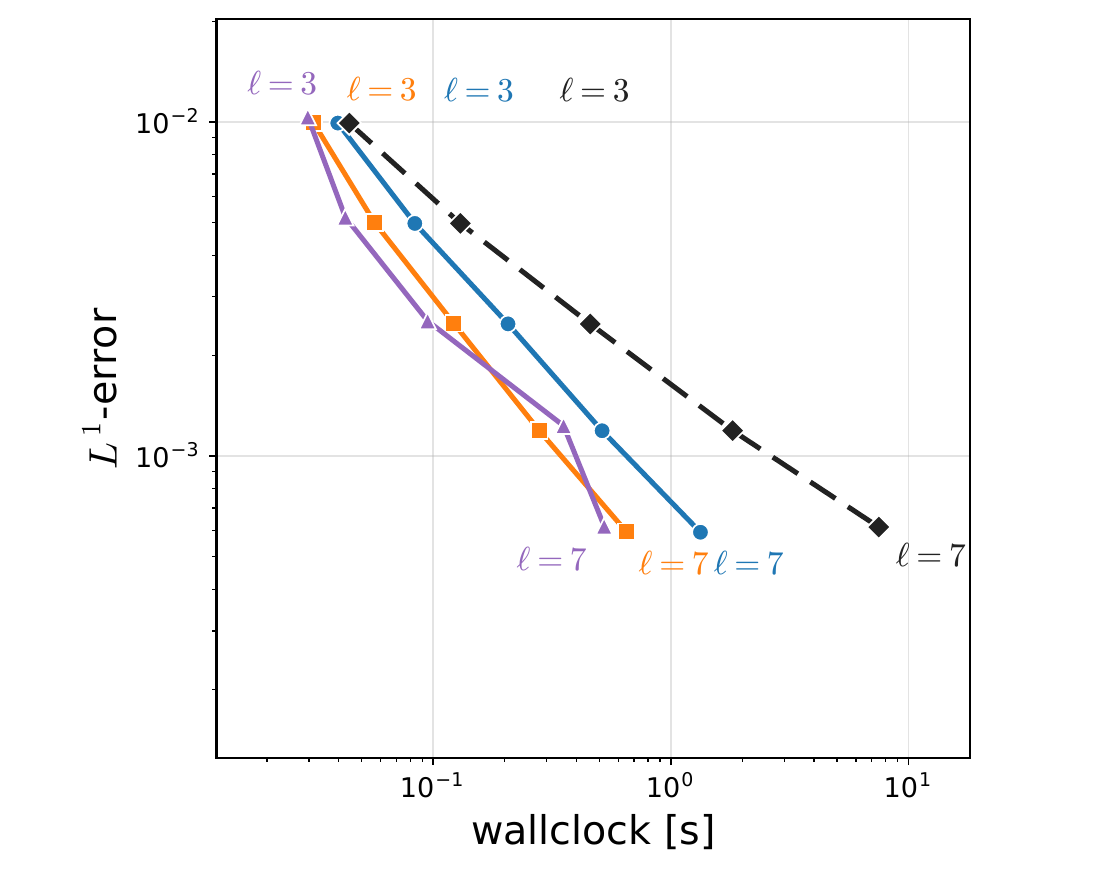}
		\caption{accuracy against runtime}\label{fig:sc-sod-cost}
	\end{subfigure}

	\vspace{2ex}

	\begin{subfigure}{0.5\textwidth}
		\includegraphics[width=\textwidth]{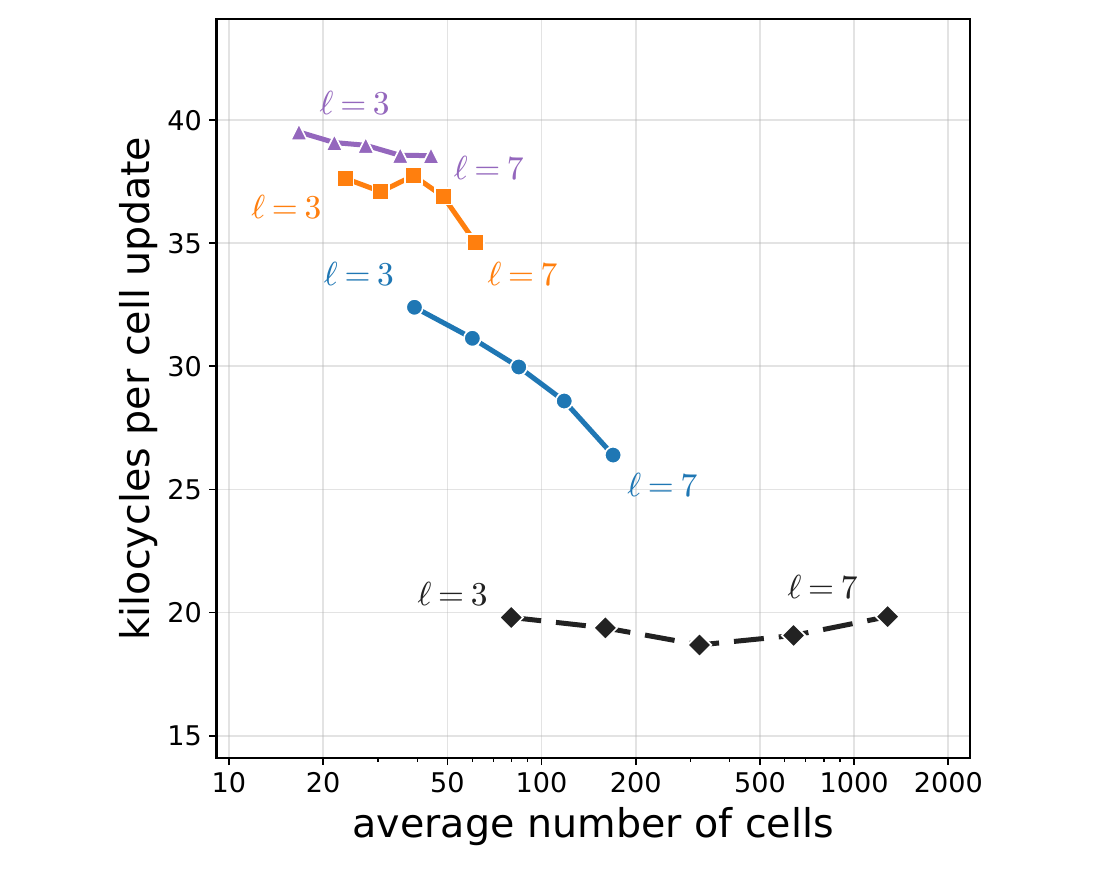}
		\caption{cost of one cell update}\label{fig:sc-sod-cell}
	\end{subfigure}\hfill
	\begin{subfigure}{0.5\textwidth}
		\includegraphics[width=\textwidth]{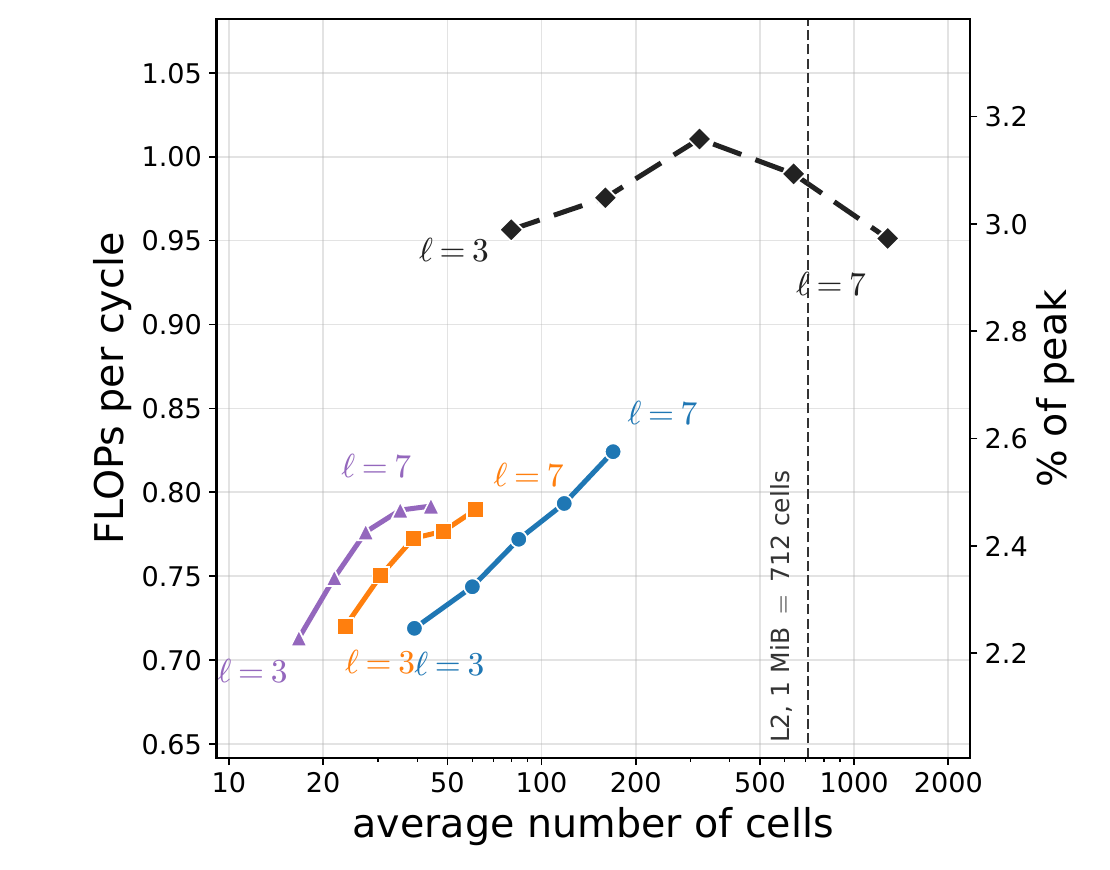}
		\caption{machine utilisation}\label{fig:sc-sod-util}
	\end{subfigure}

	\includegraphics[width=\textwidth]{figs/single_core_plots/single_core_legend.pdf}
	\caption{Sod's shock tube problem~\eqref{eq:sod-initial} for the three threshold constants and the uniform reference grid.}
	\label{fig:sc-sod}
\end{figure}

\begin{figure}[htbp]
	\centering
	\begin{subfigure}{0.5\textwidth}
		\includegraphics[width=\textwidth]{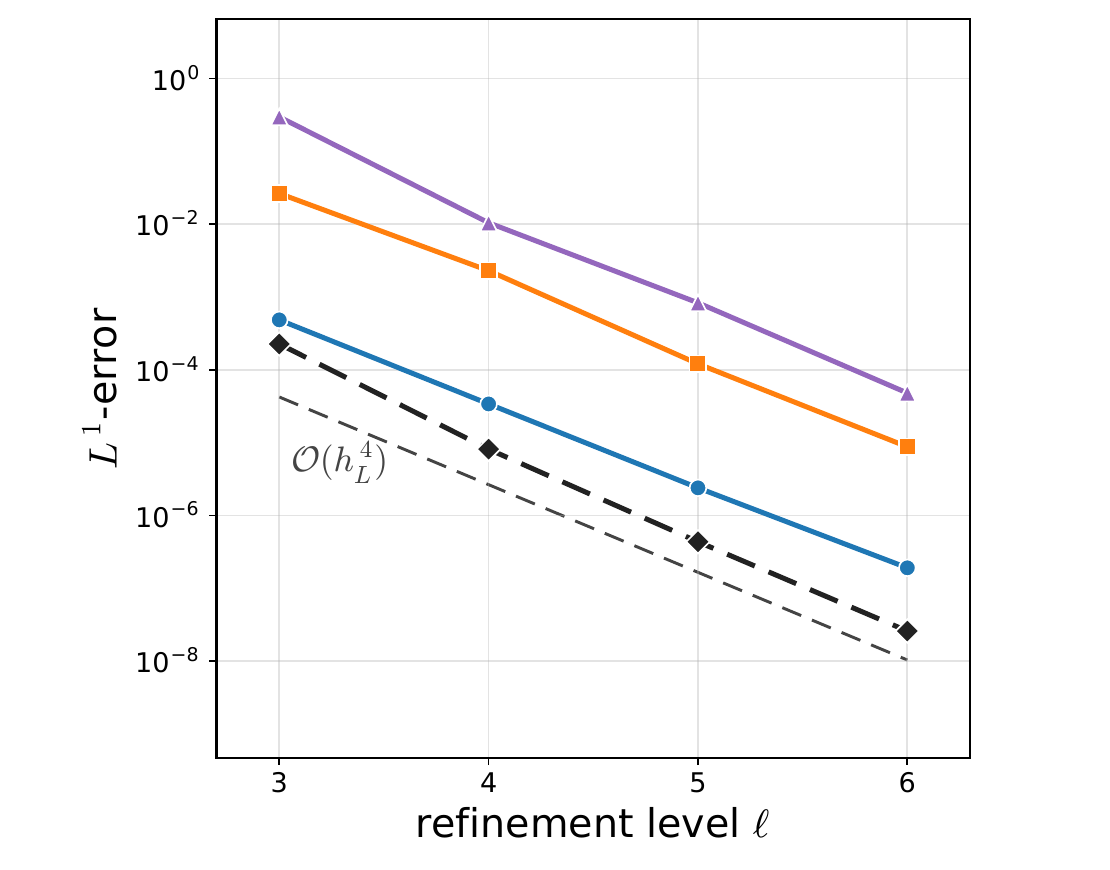}
		\caption{accuracy against refinement level}\label{fig:sc-vortex-conv}
	\end{subfigure}\hfill
	\begin{subfigure}{0.5\textwidth}
		\includegraphics[width=\textwidth]{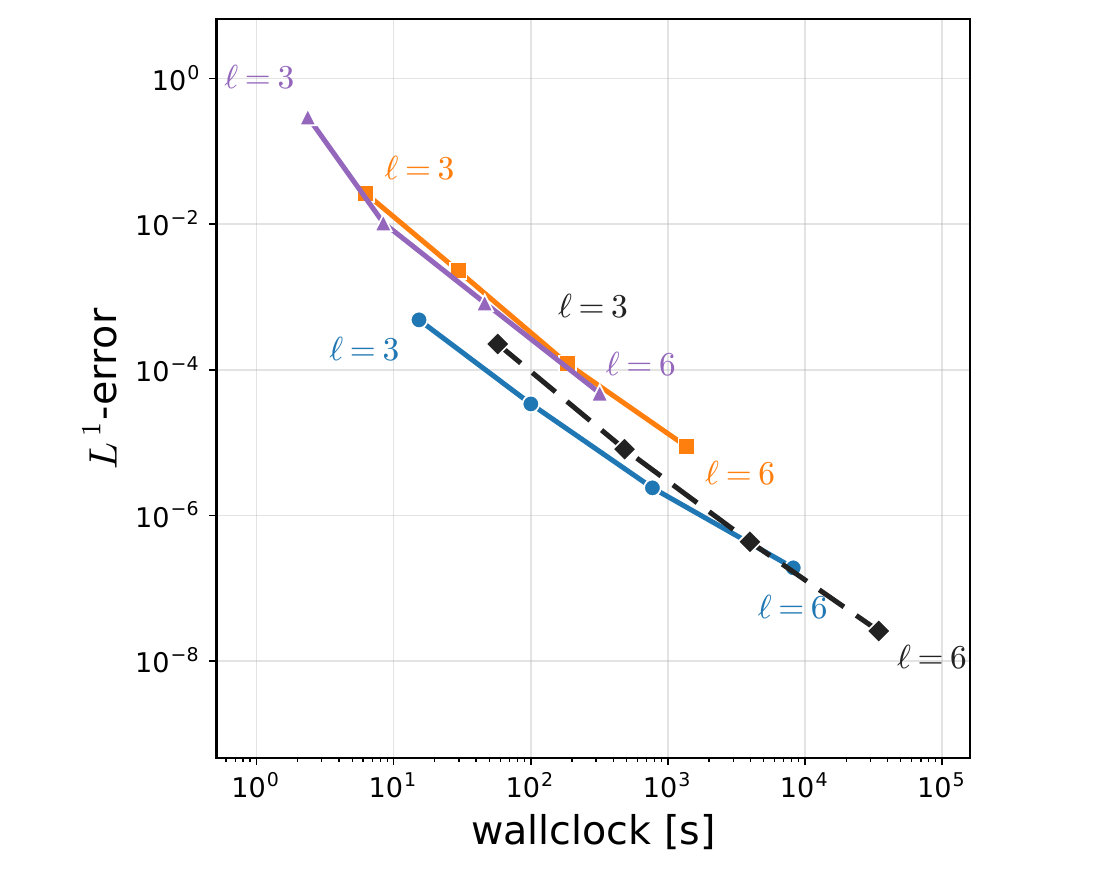}
		\caption{accuracy against runtime}\label{fig:sc-vortex-cost}
	\end{subfigure}

	\vspace{2ex}

	\begin{subfigure}{0.5\textwidth}
		\includegraphics[width=\textwidth]{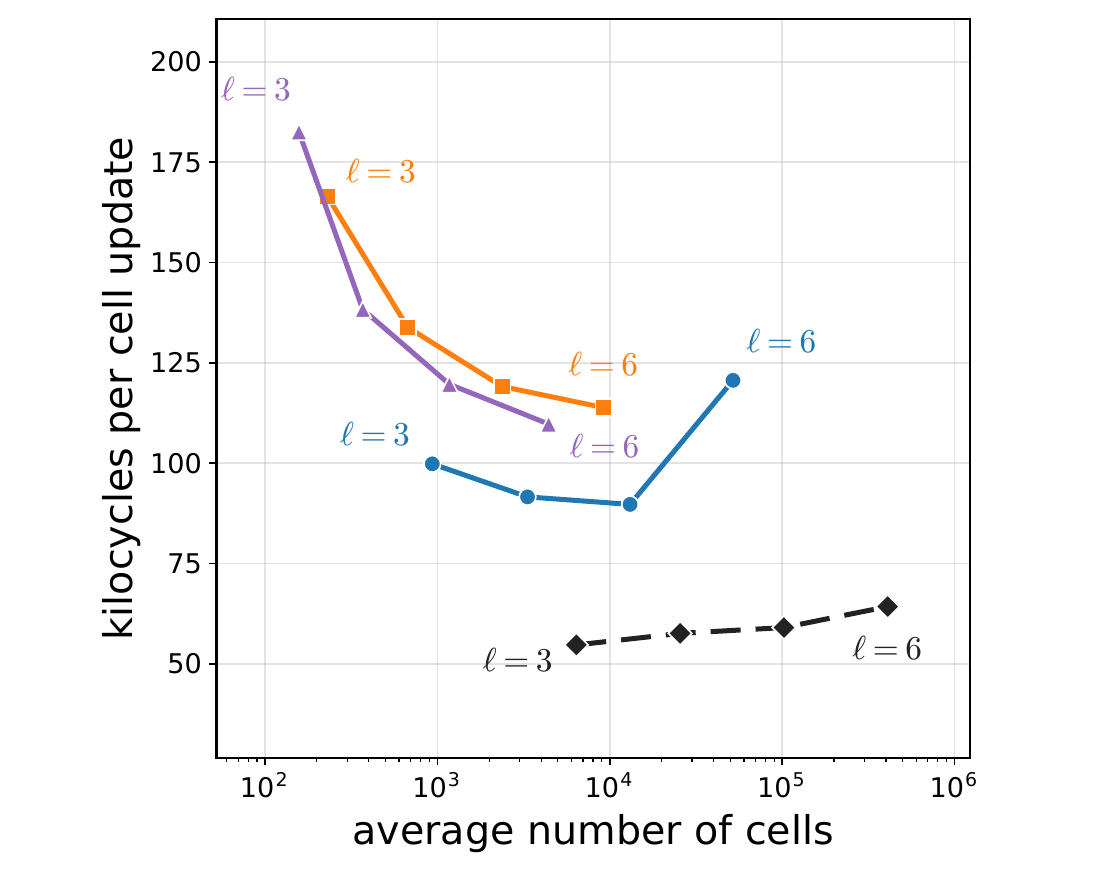}
		\caption{cost of one cell update}\label{fig:sc-vortex-cell}
	\end{subfigure}\hfill
	\begin{subfigure}{0.5\textwidth}
		\includegraphics[width=\textwidth]{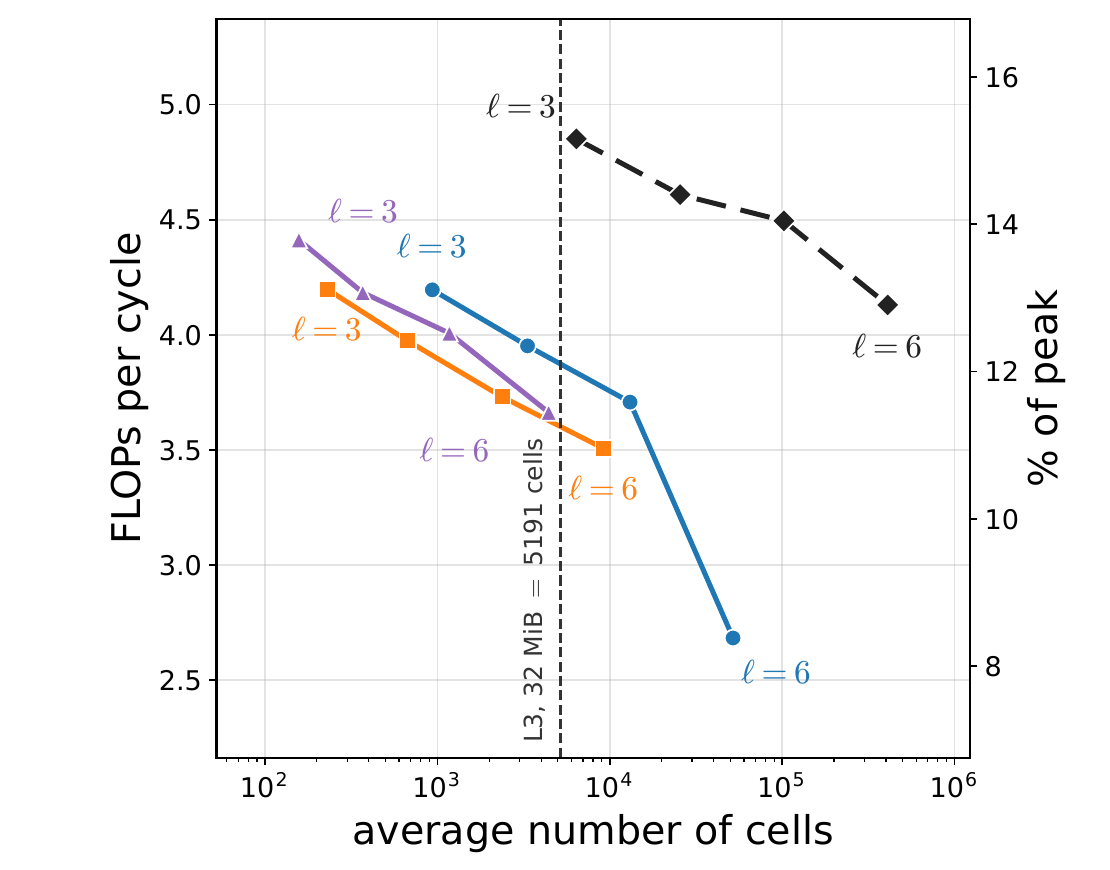}
		\caption{machine utilisation}\label{fig:sc-vortex-util}
	\end{subfigure}

	\includegraphics[width=\textwidth]{figs/single_core_plots/single_core_legend.pdf}
	\caption{Isentropic vortex~\eqref{eq:vortex-initial} for the three threshold constants and the uniform reference grid.}
	\label{fig:sc-vortex}
\end{figure}

\begin{figure}[htbp]
\centering
\begin{subfigure}{0.5\textwidth}
	\includegraphics[width=\textwidth]{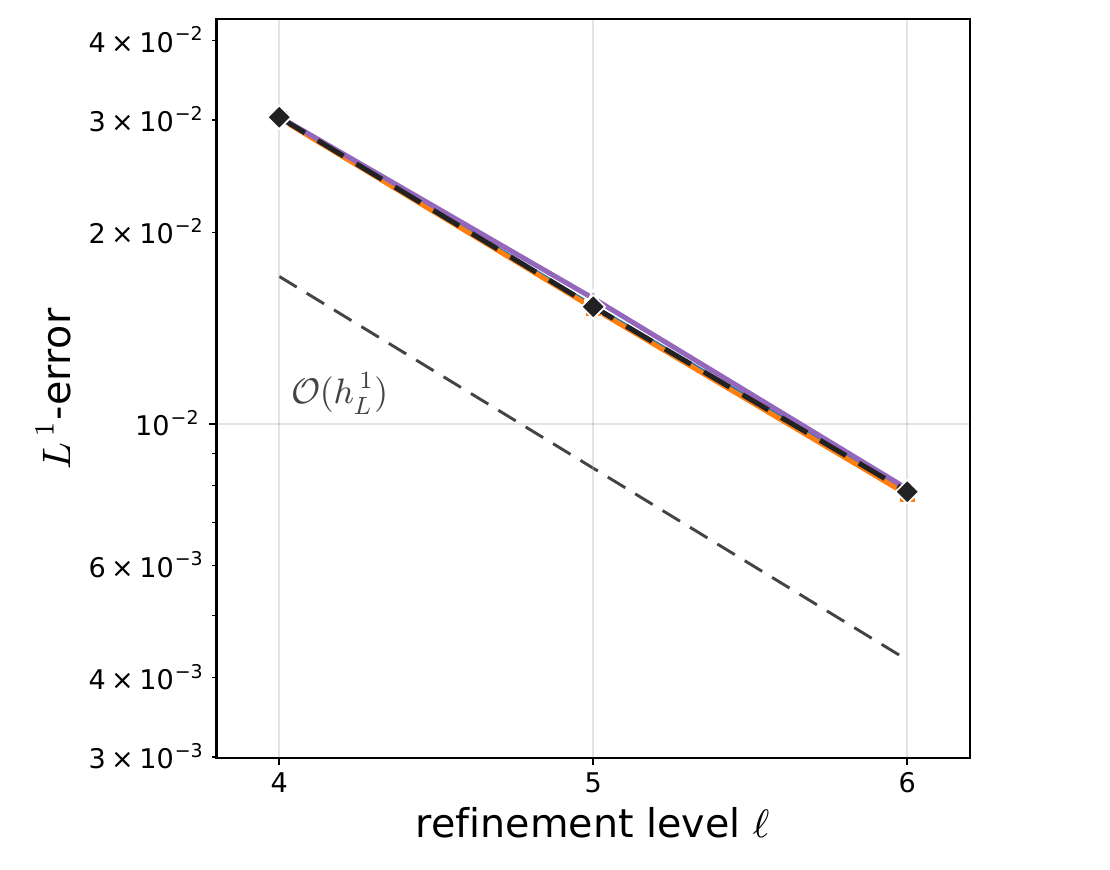}
	\caption{accuracy against refinement level}\label{fig:sc-radial-conv}
\end{subfigure}\hfill
\begin{subfigure}{0.5\textwidth}
	\includegraphics[width=\textwidth]{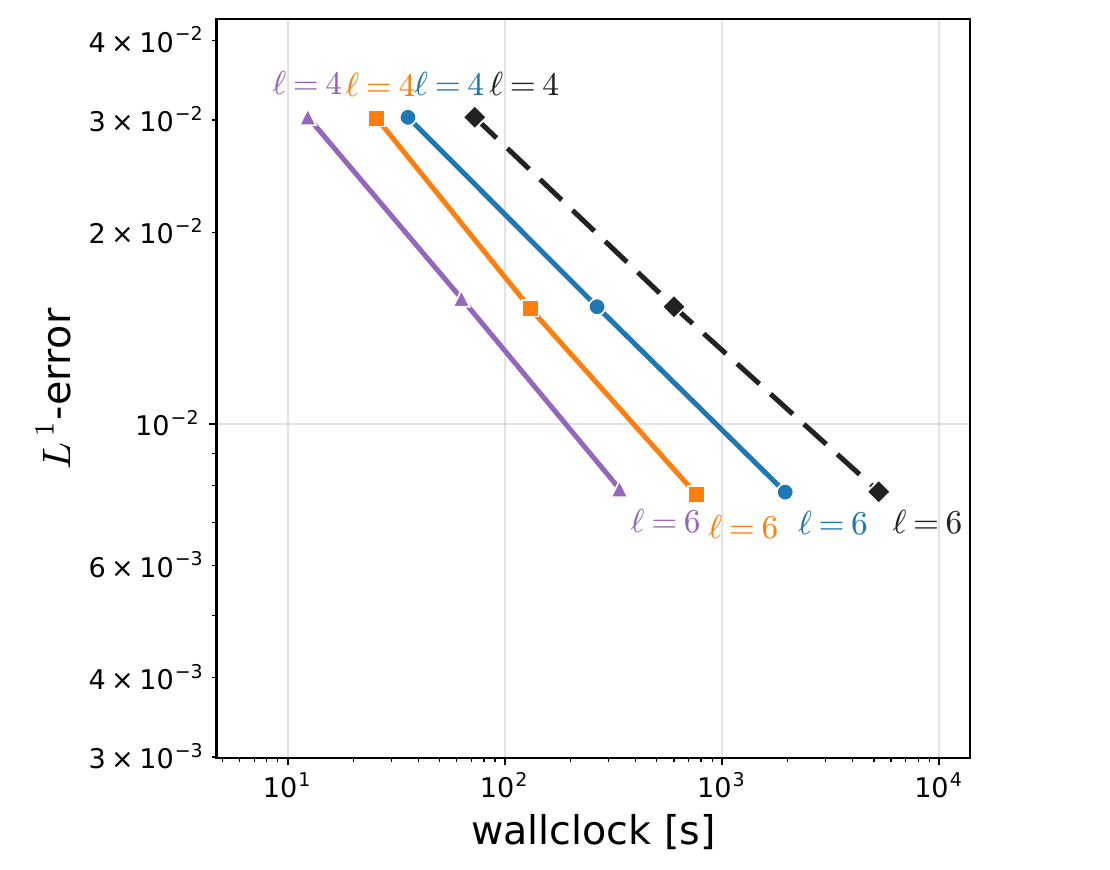}
	\caption{accuracy against runtime}\label{fig:sc-radial-cost}
\end{subfigure}

\vspace{2ex}

\begin{subfigure}{0.5\textwidth}
	\includegraphics[width=\textwidth]{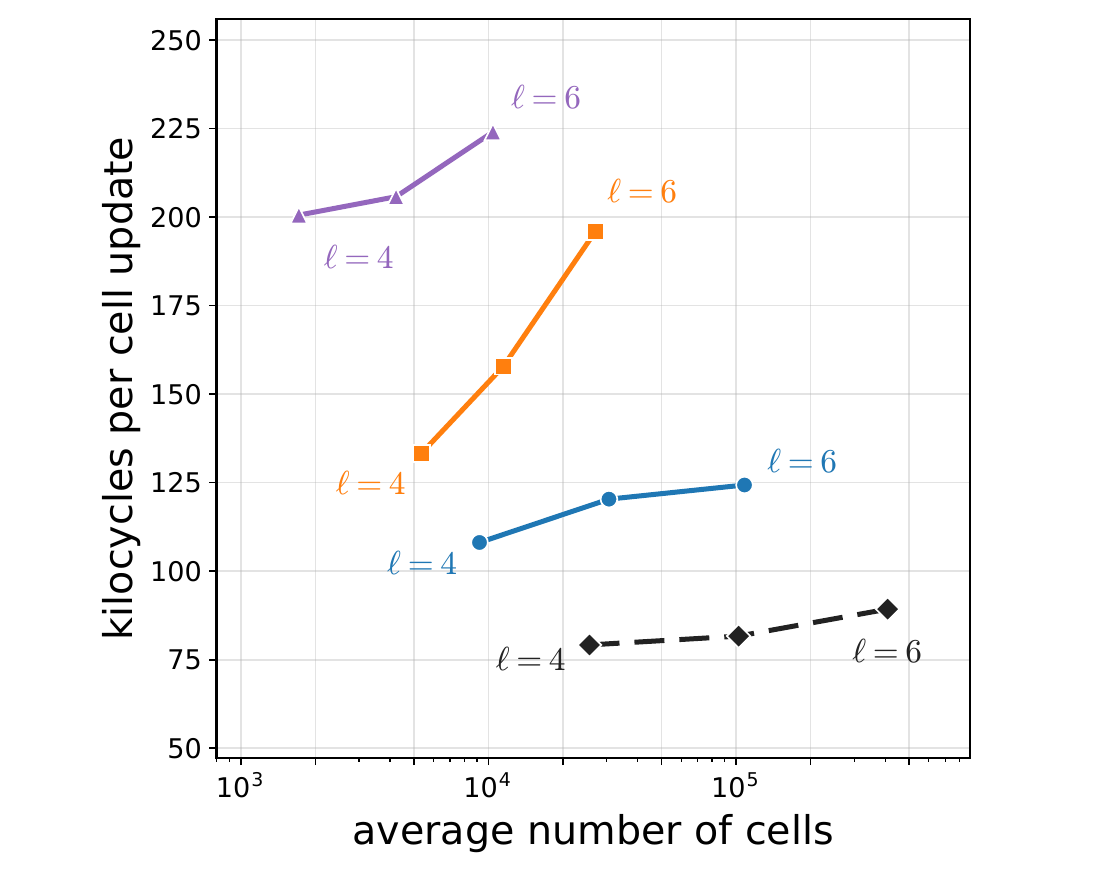}
	\caption{cost of one cell update}\label{fig:sc-radial-cell}
\end{subfigure}\hfill
\begin{subfigure}{0.5\textwidth}
	\includegraphics[width=\textwidth]{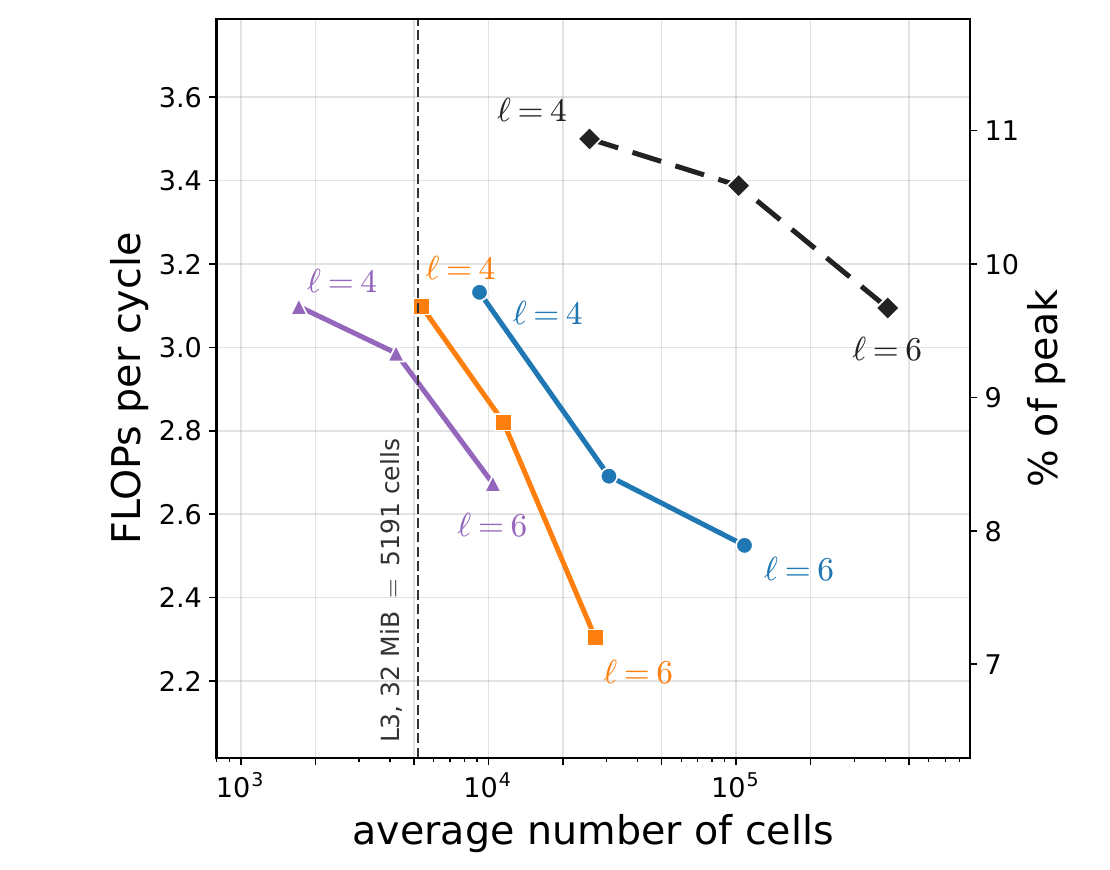}
	\caption{machine utilisation}\label{fig:sc-radial-util}
\end{subfigure}

\includegraphics[width=\textwidth]{figs/single_core_plots/single_core_legend.pdf}
\caption{Radially symmetric Sod problem~\eqref{eq:radial-sod-initial} for the three threshold constants and the uniform reference grid.}
\label{fig:sc-radial}
\end{figure}

\revakk{The results for each test case are presented in \Cref{fig:sc-hat,fig:sc-sod,fig:sc-vortex,fig:sc-radial}, one figure per case with the same four panels. Panel~(a) shows the $L^1$-error against the maximum refinement level and verifies that the grid adaptation preserves the order of the reference scheme. Panel~(b) shows the same error against wall-clock time and, therefore, the accuracy obtained per unit of computational time. 
Panel~(c) shows the cycles spent on a single cell update against the average number of cells and, thus, makes the cost of one update comparable between the uniform and the adaptive runs.
Panel~(d) shows the sustained floating-point throughput per cycle, i.e.\ the utilisation of the machine.}

\subsubsection*{\textbf{Discussion of the results}}

\revakk{Panel~(a) confirms that multiresolution-based grid adaptation preserves the order of accuracy of the reference scheme, where the dashed grey line indicates the expected rate $\mathcal{O}(h_L^p)$. As expected,  we obtain grid-convergence of order $p=4$ in the smooth cases and observe first-order convergence in the discontinuous cases. While this holds true for all tested threshold values, in the smooth cases a larger $C_\text{thr}$ leads to a coarser resolution and, thus, to a higher $L^1$-error. In the discontinuous cases we obtain almost the same error for all three threshold constants.}

\revakk{Panel~(b) shows the computational benefit of the grid adaptation.
  We plot the $L^1$-error of each adaptive and the uniformly resolved reference simulation  against its wall-clock time. For the two discontinuous cases we clearly benefit from the grid adaptation compared to the uniform simulations, i.e. we require substantially less computational time to obtain the same error. At the finest refinement level the grid consists of a few percent of the uniform grid. In the case of the Ricker wavelet the adaptation method utilizes about two thirds of the full grid, coarsening only in regions where the solution is almost constant. This is exactly the desired behaviour because otherwise the method could not maintain the order of convergence~\cite{devore1998nonlinear}. A larger threshold constant does reduce the number of cells, but the thresholding error  has to be compensated by an additional refinement level, which doubles the number of cells and halves the time step. For Sod's shock tube problem the behaviour is reversed. The thresholding error is of the same order as the discretisation error at the shock and never dominates it, so all three threshold constants reach the same error at a given refinement level. On the highest refinement level we thus significantly reduce the computational expense, resulting in a speed-up of about $14$ compared to the uniform simulation. This advantage grows with increasing refinement level, since local significant features, like discontinuities, are resolved on the same spatial width, whereas the remaining parts may be resolved on a coarser grid. Similar results are observed for the two-dimensional cases. As before, the smooth case profits less in computational effort, whilst for the discontinuous solution  a speed-up of about $16$ at the same $L^1$-error is recorded.}

\revakk{Panel~(c) shows the cost of a single cell update per time-step, i.e. spent processor cycles divided by the total number of cell updates. The simulation on a uniform grid is the cheapest in every case, as it does not include the computational overhead of the multiscale transformation. Although the per-cell cost is higher using multiresolution-based grid adaptation, the average number of cells per time step is significantly smaller than on the uniform grid, leading to reduced computational time in the discontinuous case, see Panel~(b).}

\revakk{Panel~(d) shows the sustained floating-point operations throughput. All configurations stay far below the peak of the core, so the computation is not limited by arithmetic. The grid adaptation lowers the throughput further, since the multiscale transformation and the grid bookkeeping add instructions that perform no floating-point work.
The dashed vertical line marks the number of cells at which the working set no longer fits into the cache available to a single core. For the isentropic vortex the working set of the smallest threshold constant exceeds this limit at the finest refinement level and the throughput drops noticeably, whereas the largest threshold constant stays below it.}

\subsection{Scaling studies on distributed parallel systems}
\label{sec:scaling-studies}
As a benchmark we investigate the Taylor--Green vortex for the two-dimensional Euler equations~\eqref{eq:compressible-euler} on $\Omega = [0,1]^2$ with $\kappa = 1.4$. The initial density, velocity and pressure are given by
\begin{equation}
	\label{eq:taylor-green-intial}
	\rho_0(\vec x) \equiv 1,\quad P_0(\vec x) \equiv \frac{1}{4},\text{ and } \vec v_0(\vec x) = \begin{pmatrix}
		\sin(2\pi x_0)\cos(2\pi x_1) \\
		-\cos(2\pi x_0)\sin(2\pi x_1)
	\end{pmatrix}.
\end{equation}
For the numerical simulations, we set the order of the scheme to $p = 3$ and the CFL number to $0.3$. The initial grid at the coarsest refinement level consists of $N_{\vec x_1} \times N_{\vec x_2} = \revAK{384 \times 384}$ cells. With maximum refinement level $L = 5$ a uniform grid contains \revAK{150\,994\,944} cells.
\revAK{
	The tests are performed on the JUWELS cluster system \cite{Krause2019} consisting of a total of 2271 compute nodes with two Intel Xeon Platinum 8168 CPUs, each with 24 cores at 2.7 GHz. We perform MPI parallelisation, see \cref{sec:distr-memory-parall}, using a high-speed, low-latency InfiniBand network.}

\begin{figure}[htbp]
	\centering
	\includegraphics[width=0.48\textwidth]{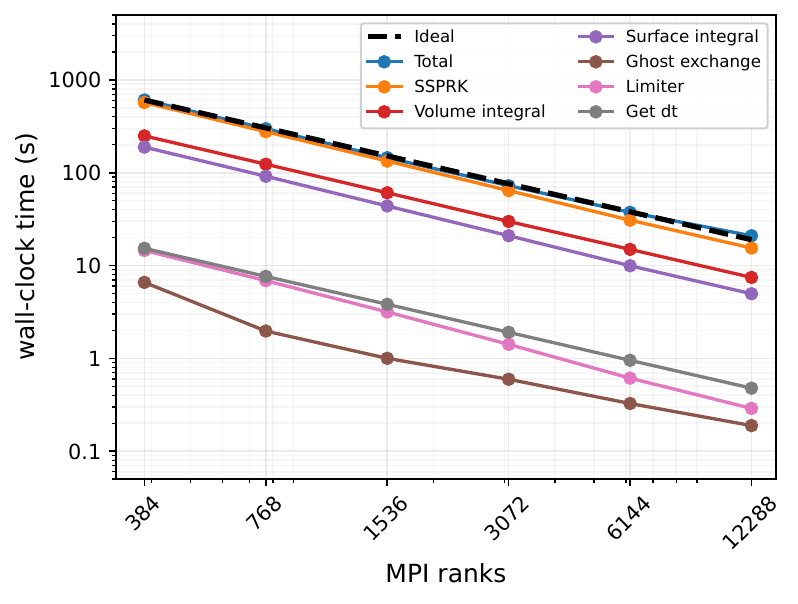}
	\hfill
	\includegraphics[width=0.48\textwidth]{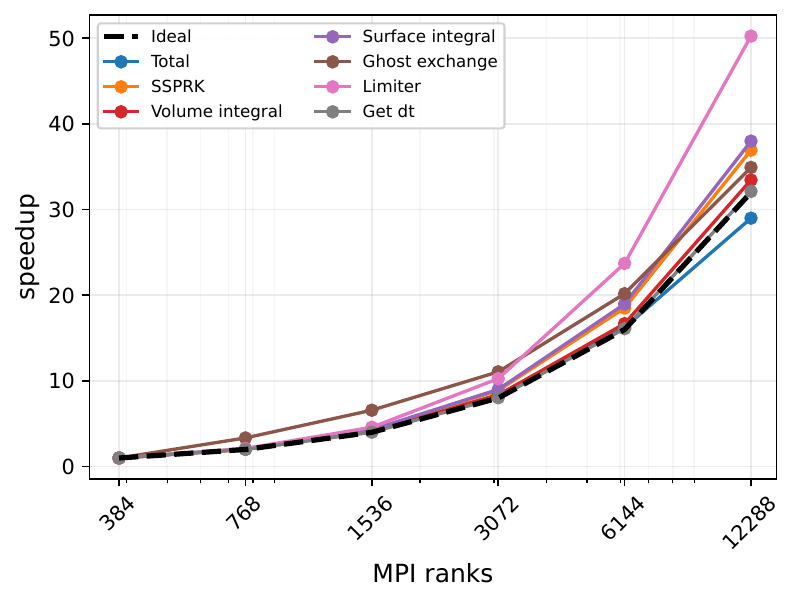}
	\caption{Strong scaling studies for the two-dimensional Euler equations on a uniform grid with up to \revAK{12\,288} MPI ranks. Left: Wall-clock time vs.\ number of ranks for the main components; Right: \revAK{Speedup relative to 384 ranks vs.\ number of ranks. On 12\,288 MPI ranks we obtain an efficiency of $90.6\%$.}}
	\label{fig:strong-scaling}
\end{figure}

\revAK{
	The scaling results for 50 time steps on a uniform grid are shown in \cref{fig:strong-scaling}. The SSPRK time-stepping dominates the total runtime, with the volume and surface integrals as the most expensive components, while the limiter and time-step computation account for only a small fraction of the total time. Since the ghost exchange is performed asynchronously, overlapping with the volume integral computation, its contribution to the overall runtime is minimised. The compute-bound operations scale slightly above ideal, which can be attributed to the local structure of the DG-SSPRK scheme. The ghost exchange shows a mild deviation from ideal scaling due to inter-node communication overhead. Overall, the strong scaling is nearly optimal across the full range of up to 12\,288 ranks.
}


\begin{figure}
	\centering
	\includegraphics[width=0.44\textwidth]{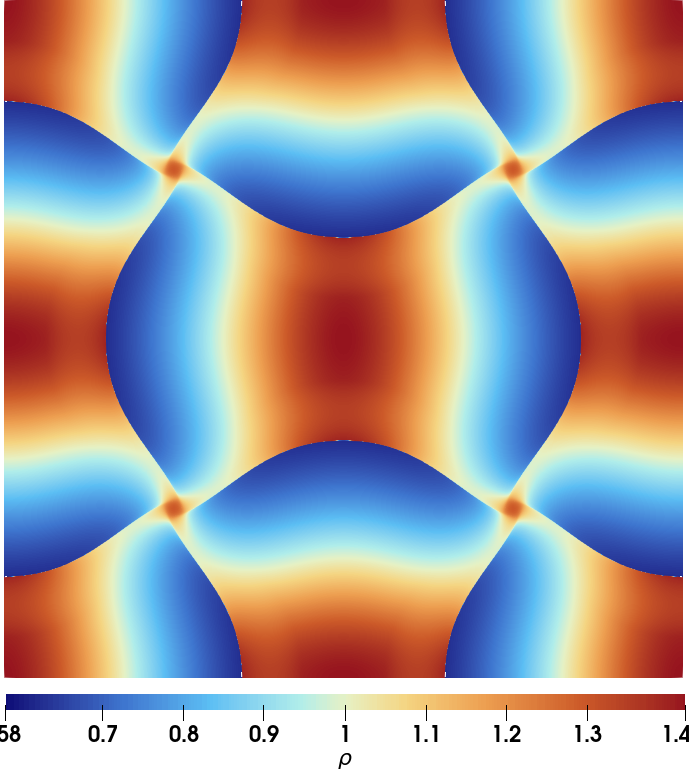}
	\hfill
	\includegraphics[width=0.44\textwidth]{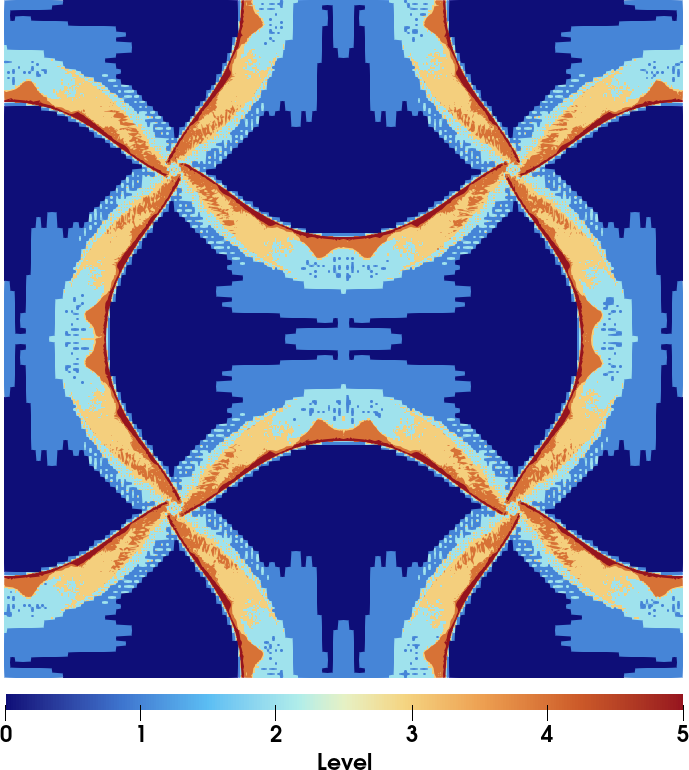}
	\caption{Solution for the two-dimensional Euler equations with initial data \eqref{eq:taylor-green-intial} at $T=0.35$. Left: density $\rho$; Right: refinement level.}
	\label{fig:taylor-green-mra}
\end{figure}

While the uniform grid results demonstrate good parallel scalability, computing on uniform grids at high refinement levels is computationally expensive. The multiresolution analysis described in \cref{sec:grid-adaptivity} allows \textsl{MultiWave} to significantly reduce the number of active cells while controlling the introduced perturbation error. \revAK{This is illustrated in \cref{fig:taylor-green-mra} for the compressible Euler equations~\eqref{eq:compressible-euler} with initial data~\eqref{eq:taylor-green-intial}.} A compression rate \revAK{larger than} $95\%$ is achieved while preserving the prescribed error tolerance, resulting in a significant reduction of computational time. The right plot of \cref{fig:taylor-green-mra} shows that the grid is only refined in regions of high variation such as discontinuities and steep gradients, while smooth areas remain at a coarser refinement level. For a detailed convergence analysis of multiresolution-based grid adaptation we refer to \cite{Gerhard2017} and the references therein.

\begin{figure}
	\centering
	\includegraphics[width=0.48\textwidth]{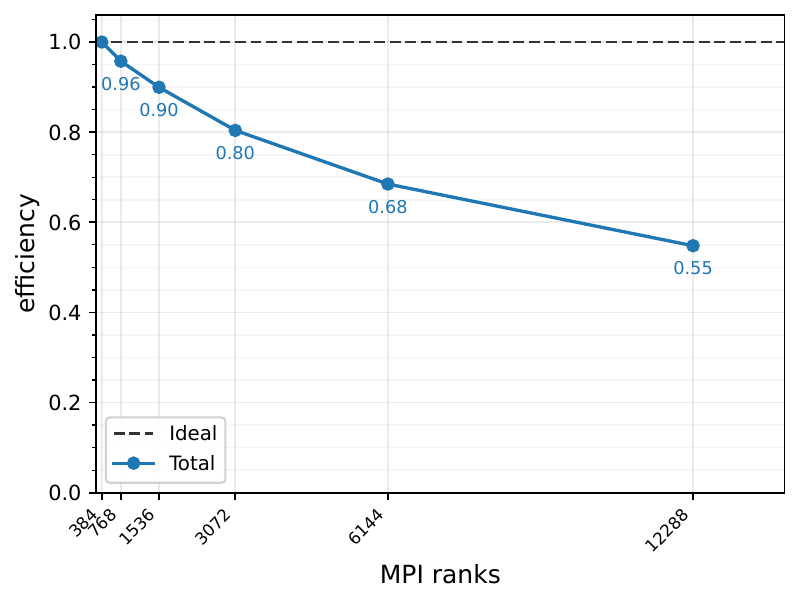}
	\hfill
	\includegraphics[width=0.48\textwidth]{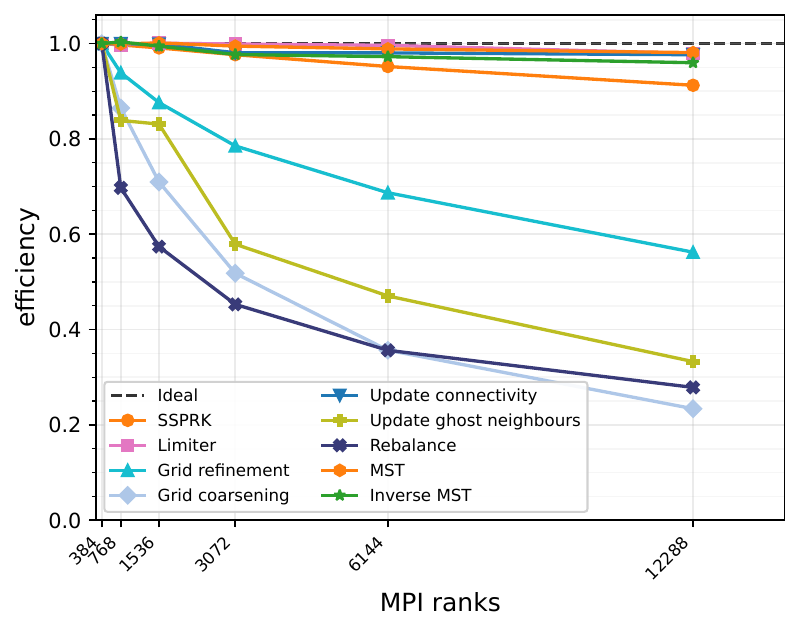}
	\caption{Weak-scaling studies for the two-dimensional Euler equations with initial data \eqref{eq:taylor-green-intial} on adaptive grids using MRA with up to \revAK{12\,288} ranks. Left: Efficiency of \textsl{MultiWave}; Right: Efficiency of main components.}
	\label{fig:weak-scaling}
\end{figure}

To demonstrate the scaling efficiency of adaptive solutions, we perform a weak-scaling study on the same test case with threshold $C_\text{thr} = \revAK{0.005}$. Since the solution is periodic, we extend the domain in the $\vec{x}_1$-direction proportionally to the number of processes, setting \revASS{$N_{\vec{x}_1} = 384 \cdot \text{Procs}/48$ with $\text{Procs} = \text{Nodes} \cdot 48$}, \revAK{so that the domain size grows proportionally with the number of MPI ranks.} The results for \revAK{5000} time steps are shown in \cref{fig:weak-scaling}.

Relative to the \revAK{384}-rank baseline, the overall efficiency stays above \revAK{55\%} up to \revAK{12\,288} ranks. The compute-bound time stepping retains \revAK{91\%}, consistent with the strong-scaling results, and the MRA operators retain \revAK{96\%} to \revAK{98\%} due to their local structure. The loss is concentrated in the adaptation pipeline. Ghost neighbour updates and load rebalancing drop to \revAK{33\%} and \revAK{28\%} as the number of inter-process cell boundaries grows with the domain, \revAK{and grid refinement and coarsening lose efficiency for the same reason, the latter most strongly at 23\%. Towards higher rank counts, the efficiency decline steepens, as the communication overhead of the adaptation pipeline grows relative to local computation.}

\begin{figure}[htbp]
	\centering
	\includegraphics[width=0.48\textwidth]{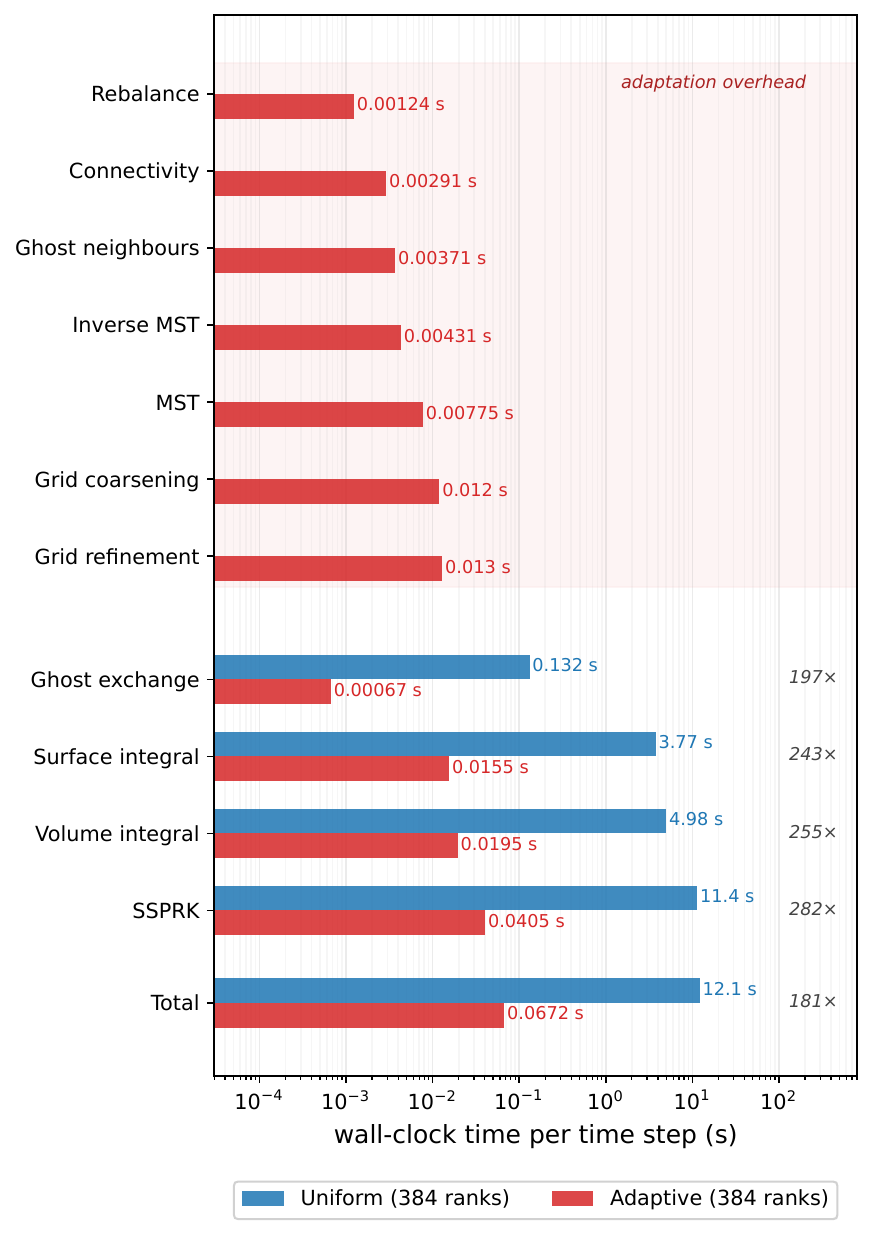}
	\hfill
	\includegraphics[width=0.48\textwidth]{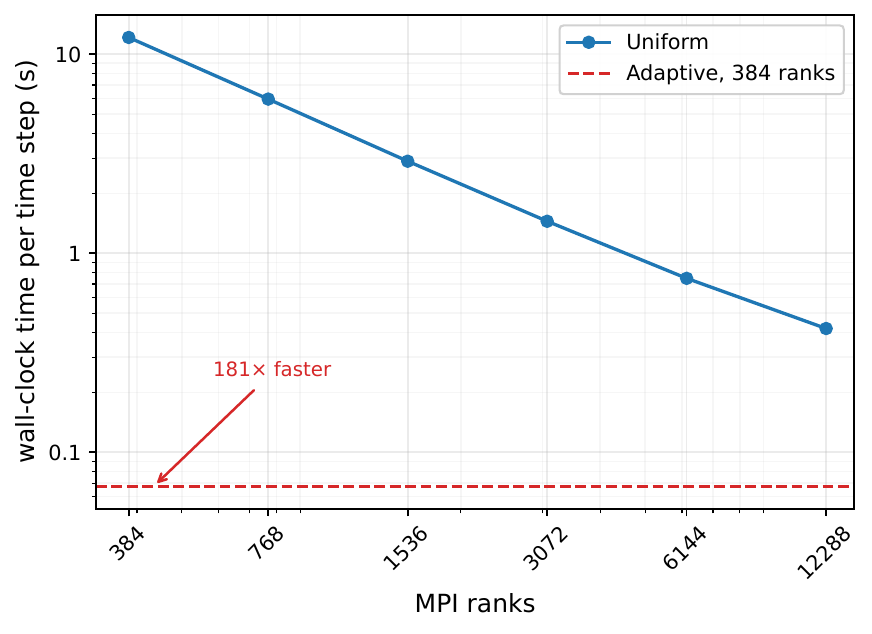}
	\caption{Adaptive vs.\ uniform scheme on \revAK{384} MPI ranks. Left: Per-component time for a single time step; Right: Total time for a single time step.}
	\label{fig:uniform-vs-adapt}
\end{figure}

\revAK{
	Finally, we quantify the benefit of multiresolution-based grid adaptation against the uniform scheme. \Cref{fig:uniform-vs-adapt} compares the per-time-step costs on 384 MPI ranks. The adaptive scheme is 181 times faster overall: ghost exchange and time stepping benefit most from the sparse cell representation, while the additional overhead introduced by the MRA operators is negligible by comparison. The right plot compares the adaptive scheme at 384 ranks with the uniform scheme running on up to 12\,288 ranks. Despite the order-of-magnitude difference in allocated resources, the uniform scheme does not reach the adaptive throughput. Since the multiresolution-based adaptation preserves \revakk{the order of accuracy} of the reference scheme, cf.~\cite{Gerhard2017}, the adaptive computation delivers \revakk{comparable} accuracy at a fraction of the cost, making high-resolution simulations feasible that would be prohibitively expensive on uniform grids.
}

\revAK{
	The scaling studies show that \textsl{MultiWave} scales well up to the tested rank counts. The DG core remains close to ideal scaling in both cases, while the overhead from dynamic grid adaptation stays moderate, even at 12\,288 ranks. Combined with the $181\times$ speedup over the uniform scheme, these results confirm that the adaptive framework is both scalable and cost-efficient.}

\subsection{Multiresolution analysis for random hyperbolic conservation laws}
\label{sec:weighted-mra}
\begin{figure}[htbp]
	\begin{center}
		\includegraphics[width=0.95\textwidth]{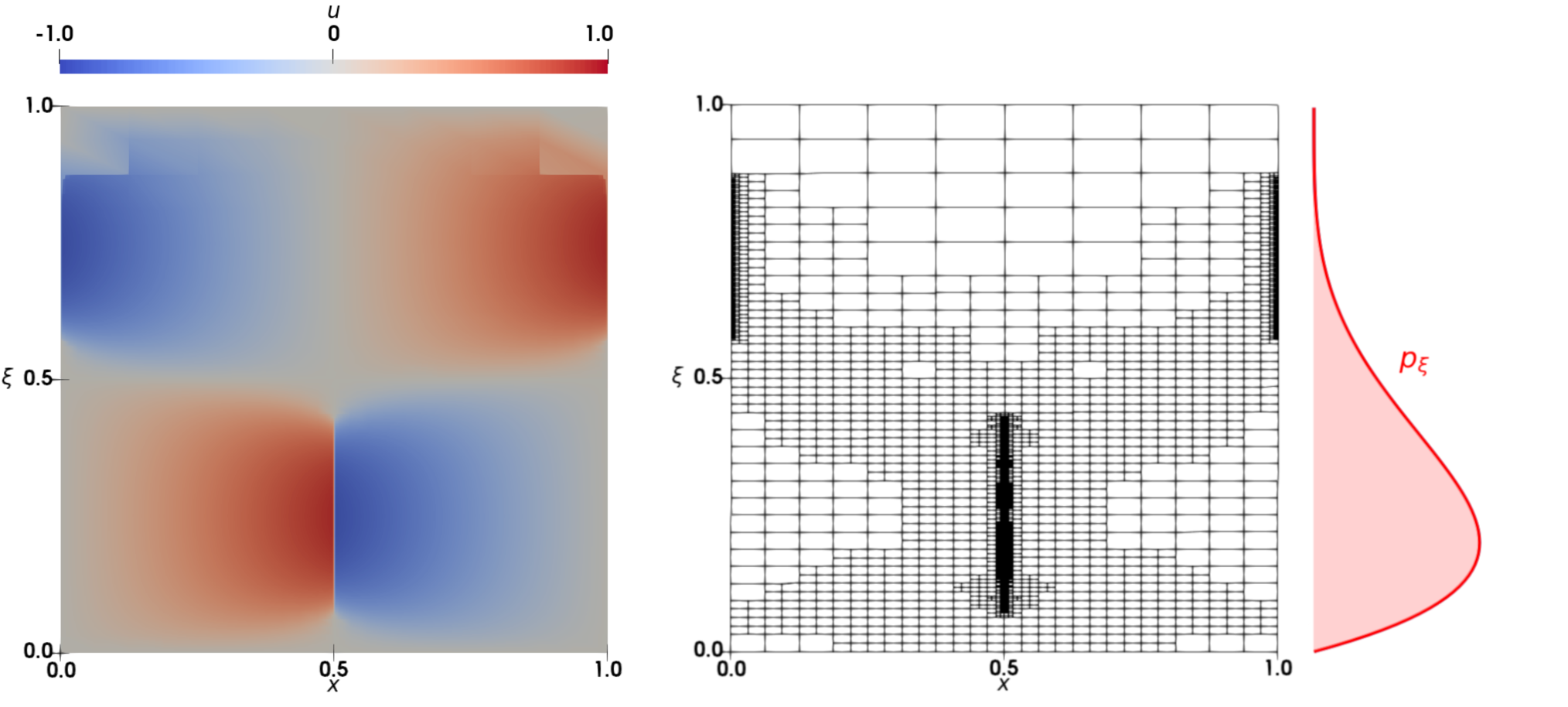}
	\end{center}
	\caption{Weighted MRA for Burgers' equation with Beta-distributed random variable \cite{herty2022higher}.}
	\label{fig:weighted-mra}
\end{figure}

A possible extension of \textsl{MultiWave} has been proposed in \cite{herty2022higher,Kolb2023,herty2024novel} in the context of random hyperbolic conservation laws of the form
\revASS{\begin{equation}
		\label{eq:random-hcl}
		\frac{\partial \vec{u}}{\partial t}(t,\vec{x}; \vec{\xi}(\omega)) + \nabla\cdot \vec{F}(\vec{u}(t,\vec x;\vec{\xi}(\omega)))= \vec{0},
	\end{equation}}
where $\vec{\xi}$ is a random variable. Since discontinuities in the spatial direction lead to discontinuities in the random space, efficient numerical solvers are necessary to obtain the stochastic moments. In \cite{herty2022higher}, the stochastic dimensions in \eqref{eq:random-hcl} are interpreted as additional spatial dimensions, leading to a higher-dimensional problem. This  allows us to investigate the interaction between the stochastic and spatial scales, and also enables weighted MRA based on the stochastic moments rather than on the solution itself. This is achieved by modifying the threshold value in \cref{eq:significant-detail} using the probability density $p_{\vec\xi}$ of the underlying random variable, leading to grids that are more refined in regions of high probability density and coarser in regions of low probability density. An example with a Beta-distributed random variable is shown in \cref{fig:weighted-mra}.
\revAK{In \textsl{MultiWave}, this is realised by replacing the local threshold value in~\eqref{eq:significant-detail} with a weighted threshold based on the local probability density of a cell $V_\lambda^{\vec \xi}$ in stochastic direction for a fixed spatial position $\varepsilon_{\lambda,L} = 2^{\ell-L}\|p_{\vec\xi}\|^{-1}_{L^\infty\left(V_\lambda^{\vec \xi}\right)}\cdot\varepsilon_L$. This choice controls the perturbation error in the stochastic moments. Note that this may lead to locally underresolved solutions as seen in \cref{fig:weighted-mra} but still guarantees convergence in the stochastic moments while reducing the number of needed cells.}

An extension of this formulation has been investigated in \cite{Kolb2023, herty2024novel} where multiple lower-dimensional functions have been used to represent the spatial and stochastic scales, respectively.
\revAK{In \textsl{MultiWave}, different refinement strategies are applied independently to the spatial and stochastic scales. A parent \textsl{MultiWave} object orchestrates all lower-dimensional instances, performing time stepping on all realisations while simultaneously updating both grids. Whenever the grid is adapted in one dimension, solutions in the other are inserted or removed accordingly. This is possible because each lower-dimensional solution is an independent object with its own configuration, including spatial dimension, polynomial order, and refinement strategy, and the template-based design allows them to be freely composed without modification to the surrounding framework.}

\subsection{Coupling}
\label{sec:coupling}

Fast moving objects in fluids may lead to a pressure drop below the evaporation
point such that hot vapour bubbles form. As the bubbles are surrounded by the
cold fluid they collapse under the ambient pressure. The waves that are emitted in
this process are able to generate extreme temperatures and pressures. These conditions
damage the material of the fast moving objects, a phenomenon known as cavitation
erosion, which is a major problem in the design of ship propellers, centrifugal
pumps and water turbines.
For this reason, it is essential to understand
the coupling of the fluid and the solid material.

In~\cite{sikstelAnalysisNumericalMethods2020,herty2021coupling} a fluid-structure coupling problem for steel and water in two
spatial dimensions has been simulated. The fluid part is governed by the Euler equations~(\ref{eq:compressible-euler}) with stiffened gas equations of state
\[
  P = (\kappa -1)\rho(e - Q) - \kappa\pi
\]
where $\kappa > 1$ denotes the adiabatic index, $e$ the internal energy, $Q$  the formation
enthalpy and $\pi$ the pressure stiffness. Setting $Q=pi=0$  results in the ideal gas EoS~\eqref{eq:perfect-gas-eos}.
The solid part is governed by a linear elastic system
\begin{align*}
  &\frac{\partial \vec{w}}{\partial t} - \frac{1}{\rho_s}\nabla \cdot \vec{\sigma} = \vec{0},\\
  &\frac{\partial \vec{\sigma}}{\partial t} - \lambda (\nabla \cdot \vec{w})\cdot \mathbf{I} - \mu (\nabla\vec{w} + \nabla\vec{w}^\top) = \vec{0},
\end{align*}
where conserved quantities are the deformation velocity $\vec{w}\in \mathbb{R}^d $ and the symmetric shear stress tensor $\vec{\sigma} = (\sigma_{i,j})^d_{i,j=1} = \vec{\sigma}^\top \in \mathbb{R}^{d\times d} $. Model parameters are the material density $\rho_s$ and the Lam\'e constants $\lambda, \mu > 0$ that relate to the dilatation wave velocities $c_1, c_2 \geq 0$ as $c_1^2 = \frac{2}{\rho_s}(\mu + \lambda)$ and $c_2^2 = \frac{\mu}{\rho_s}$.

\begin{figure}[!h]
	\centering
	\includegraphics[width=.8\textwidth]{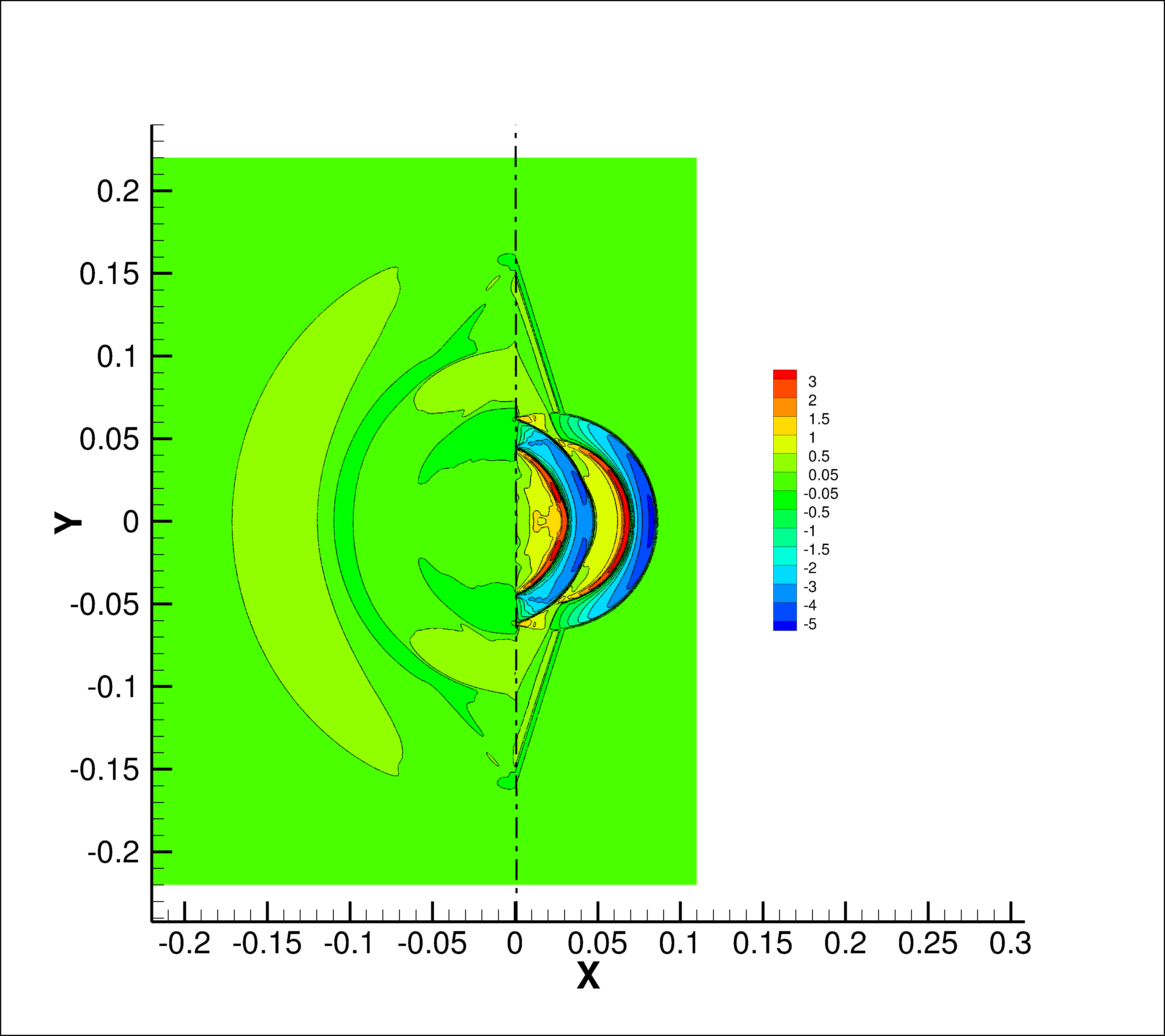}
	\caption{Contour plot of the velocity $\vec v_1$ of the 2D bubble simulation \cite{sikstelAnalysisNumericalMethods2020}.}
	\label{fig:fsi-pressure}
\end{figure}

The two systems are separated by a fixed interface at $x_1 = 0$ where a coupling condition is prescribed that is enforced by a Riemann-type solver that provides boundary information to each phase. The initial conditions of the solid are constant such that the coupling conditions are fulfilled and the  fluid part contains a hot bubble. \Cref{fig:fsi-pressure} shows the result at final time where the solution has produced a complex wave pattern including von Schmidt waves.

The coupling procedure has been realised by defining two different \textsl{MultiWave} objects. Since solutions of the coupling Riemann problems are required at each RK-stage, the two objects communicate their respective coupling-interface data during the boundary condition call. Only the ranks adjacent to the coupling interface participate in the MPI call. These ranks are known a priori, because complete subtrees of level-0 cells are assigned to a single MPI rank. After exchanging the interface data each object computes the solutions of the coupling Riemann problems and projects them to the boundary cell coefficients. Thus, the coupling procedure is implemented non-intrusively in a few concise steps.

\section{Conclusion}

\revASS{
	The central outcome of this work is that a high level of mathematical abstraction needs not be traded for computational performance. By encoding the complete numerical hierarchy as compile-time policies, \textsl{MultiWave} permits the PDE system, Riemann solver, spatial and temporal discretisations, basis representation, limiter, refinement strategy, and coupling mechanism to be exchanged independently. Thus, numerical methods can be expressed in a form that remains as close to their mathematical formulation as possible, while the compiler reduces each  configuration to a specialised executable without runtime polymorphism or abstraction overhead.}

\revASS{
	This design turns \textsl{MultiWave} into a computational laboratory in which theoretical concepts can be translated into scalable numerical experiments with little intervening infrastructure work. The close correspondence between the weak formulation, the multiresolution analysis, and their software representations preserves the structure required for conservation, stability, perturbation estimates, and rigorous error analysis. At the same time, new operators and application-specific methods can reuse the existing adaptive data structures, MPI parallelisation, load balancing, and communication machinery. The examples considered here, including viscous and non-conservative operators, well-balanced shallow-water schemes, stochastic adaptivity, and fluid--structure coupling, demonstrate that this modularity extends across substantially different mathematical models and algorithmic settings.}

\revASS{
	This generality does not compromise performance. The DG core exhibits nearly ideal strong scaling up to 12\,288 MPI ranks, attaining an efficiency of $90.6\%$, while the complete adaptive method maintains reasonable weak-scaling efficiency at the same scale. For the Taylor--Green benchmark, multiresolution-based adaptivity reduces the time per step by a factor of $181$ relative to the uniform-grid computation on 384 ranks while retaining \revakk{the order of accuracy} of the reference discretisation. These results show that the framework's abstractions disappear from the performance-critical execution path.}

\revASS{
	\textsl{MultiWave}, therefore, provides a foundation in which mathematical rigour, algorithmic flexibility, and large-scale performance coexist rather than being treated as competing objectives. Ongoing work on tensor-structured coefficients, $hp$-adaptivity, and rigorous a posteriori error control will make use of the framework's combination of abstraction, analytical structure, and computational efficiency. The resulting framework is intended not merely as a solver for a fixed class of PDEs, but as reusable infrastructure for turning new analytical ideas into efficient large-scale computations.
}




\begin{acks}
	This work was funded by the Deutsche Forschungsgemeinschaft (DFG, German Research Foundation) -- 320021702/GRK2326 -- Energy, Entropy and Dissipative Dynamics (EDDy) and by the BMFTR project ADAPTEX (FKZ: 16ME0670).
	\revAK{The authors acknowledge the computing time granted on the JUWELS supercomputer at the Jülich Supercomputing Centre, Forschungszentrum Jülich.}
	Both authors are very grateful to their dear colleagues who taught and supported them  during the journey of \textsl{MultiWave} implementation, in particular we thank Siegfried M\"uller, Michael Herty, Stephan Gerster, Igor Voulis and Frank Knoben.
	Furthermore, we thank all student assistants who have worked with us.
	Last but not least, the authors thank the W.H.I.P. \revAK{group} of the IGPM.
\end{acks}

\bibliographystyle{abbrv}
\bibliography{refs/MyLib}

\appendix

\newpage
\section{Ghost neighbour update}
\label{sec:ghost-neigh-update}

\begin{algorithm}[htbp]
	\caption{Ghost neighbour update after grid adaptation}
	\label{alg:ghost-neigh-update}
	\begin{algorithmic}[1]
		\Require{distributed solution $\mathcal{U} = \bigcup_{q=0}^{Q-1} \mathcal{U}_q$ over $Q$ MPI ranks, ghost cells $\mathcal{G}$, adjacency list $\mathcal{A}$}
		\Procedure{update\_ghost\_neighbours}{$\mathcal{U}_r,\, \mathcal{A}$} \Comment{executed on rank $r$}
		\State $\mathcal{S} \gets \varnothing$,\quad $\mathcal{V} \gets \varnothing$
		\For{each $V_\kappa \in \mathcal{G}$}
		\State $q \gets \textsc{get\_process}(V_\kappa)$
		\For{each $(V_\mu,\, \vec{n}) \in \textsc{get\_neighbours}(\mathcal{A},\, V_\kappa)$}
		\If{$V_\mu \in \mathcal{U}_r$}
		\State $\mathcal{S} \gets \mathcal{S} \cup \{(V_\mu,\, q)\}$
		\Comment{tell process $q$ about our boundary cell}
		\State $\mathcal{V} \gets \mathcal{V} \cup \{V_\mu\}$
		\Comment{candidate local neighbour for new ghosts}
		\EndIf
		\EndFor
		\State $\textsc{remove\_all\_neighbours}(\mathcal{A},\, V_\kappa)$
		\EndFor
		\State $\mathcal{G} \gets \varnothing$
		\State $\Lambda_\mathrm{recv} \gets \textsc{lmi\_exchange}(\mathcal{S})$
		\For{each $V_\lambda \in \Lambda_{\mathrm{recv}}$}
		\State add $V_\lambda$ to $\mathcal{G}$
		\State $\textsc{connect\_new\_ghost}(\mathcal{A},\, V_\lambda,\, \mathcal{V})$
		\EndFor
		\EndProcedure

		\Statex

		\Require{new ghost cell $V_\lambda$, adjacency list $\mathcal{A}$, candidate local neighbours $\mathcal{V}$}
		\Procedure{connect\_new\_ghost}{$\mathcal{A},\, V_\lambda,\, \mathcal{V}$}
		\For{each $\vec{n} \in \{\pm\vec{e}_1, \ldots, \pm\vec{e}_d\}$}
		\State $V_{\hat\kappa} \gets \textsc{next\_index}(V_\lambda,\, \vec{n})$
		\If{$\exists\, V_\mu \in \mathcal{V}$ with $\ell_\mu \leq \ell_{\hat\kappa}$
			and $\textsc{get\_parent}(V_{\hat\kappa},\, \ell_{\hat\kappa} - \ell_\mu) = V_\mu$}
		\State $\textsc{add\_neighbours}(\mathcal{A},\, V_\lambda,\, [(V_\mu,\, \vec{n})])$
		\Comment{same-level or coarser local neighbour}
		\Else
		\State $\mathcal{F} \gets \{V_\mu \in \mathcal{V} :
			\textsc{get\_parent}(V_\mu,\, \ell_\mu - \ell_{\hat\kappa}) = V_{\hat\kappa}\}$
		\Comment{finer local neighbours}
		\State $\textsc{add\_neighbours}(\mathcal{A},\, V_\lambda,\,
			[(V_\mu,\, \vec{n}) : V_\mu \in \mathcal{F}])$
		\EndIf
		\EndFor
		\EndProcedure
	\end{algorithmic}
\end{algorithm}

\end{document}